\documentclass[10pt,letterpaper]{amsart}
\usepackage{charter}
\usepackage[OT2,OT1]{fontenc}
\usepackage{textcomp}
\usepackage[latin9]{inputenc}
\usepackage{float}
\usepackage{units}
\usepackage{mathrsfs}
\usepackage{bm}
\usepackage{amsbsy}
\usepackage{amstext}
\usepackage{amsthm}
\usepackage{amssymb}
\usepackage{cancel}
\usepackage{stmaryrd}

\makeatletter

\pdfpageheight\paperheight
\pdfpagewidth\paperwidth

\providecommand{\tabularnewline}{\\}

\numberwithin{equation}{section}
\numberwithin{figure}{section}
\theoremstyle{plain}
\newtheorem{thm}{\protect\theoremname}[section]
\theoremstyle{definition}
\newtheorem{defn}[thm]{\protect\definitionname}
\theoremstyle{remark}
\newtheorem{rem}[thm]{\protect\remarkname}
\theoremstyle{plain}
\newtheorem{prop}[thm]{\protect\propositionname}
\theoremstyle{plain}
\newtheorem{cor}[thm]{\protect\corollaryname}
\theoremstyle{plain}
\newtheorem{lem}[thm]{\protect\lemmaname}
\theoremstyle{definition}
\newtheorem{example}[thm]{\protect\examplename}

\usepackage{amsfonts}
\usepackage{mathrsfs}
\usepackage{graphicx}
\usepackage{amscd}
\usepackage{enumerate}
\usepackage{url}
\usepackage{multirow,bm}
\usepackage{mathtools} 
\usepackage{xypic}
\usepackage{tikz-cd}

\newcommand{\cyr}{
\renewcommand\rmdefault{wncyr} \renewcommand\sfdefault{wncyss} \renewcommand\encodingdefault{OT2} \normalfont
\selectfont
}
\DeclareTextFontCommand{\textcyr}{\cyr}    

   \def\@settitle
{
\begin{center} \baselineskip14\p@\relax \LARGE\@title \end{center}
}

\usepackage[normalem]{ulem}
\usepackage[mathscr]{eucal}
\numberwithin{equation}{section}

\usepackage{verbatim}
\usepackage{amscd}
\usepackage{bm}\usepackage[mathscr]{eucal}
\usepackage[all]{xy}
\usepackage{enumerate}
\usepackage{color}
\usepackage{colortbl}

\input xyimport.tex

\title{{\bf  Key varieties for prime $\mQ$-Fano threefolds \\ defined by Jordan algebras of cubic forms. \\ Part I}}

\author{Hiromichi Takagi}

\address{Department of Mathematics, Gakushuin University, 
Mejiro, Toshima-ku, Tokyo 171-8588, Japan}
\email{hiromici@math.gakushuin.ac.jp}

\newcommand{\sO}{\mathcal{O}}

\newcommand{\mA}{\mathbb{A}}
\newcommand{\mC}{\mathbb{C}}

\newcommand{\mN}{\mathbb{N}}
\newcommand{\mP}{\mathbb{P}}
\newcommand{\mQ}{\mathbb{Q}}

\newcommand{\rank}{\mathrm{rank}\,}

\numberwithin{equation}{section}

 \newcounter{myparagraph}[subsection]

\makeatother

\providecommand{\corollaryname}{Corollary}
\providecommand{\definitionname}{Definition}
\providecommand{\examplename}{Example}
\providecommand{\lemmaname}{Lemma}
\providecommand{\propositionname}{Proposition}
\providecommand{\remarkname}{Remark}
\providecommand{\theoremname}{Theorem}

\begin{document}
\maketitle 
\begin{abstract}
The first aim of this paper is to construct a $13$-dimensional affine
variety $\mathscr{H}_{\mathbb{A}}^{13}$ related to $\mP^{2}\times\mP^{2}$-fibration
with relative Picard number $1$. It is well-known that the affine
cone over the Segre embedded $\mP^{2}\times\mP^{2}$ is defined as
the null locus of the so called $\sharp$-mapping of the 9-dimensional
nondegenerate quadratic Jordan algebra $J$ of a cubic form. Inspired
with this fact, we construct $\mathscr{H}_{\mathbb{A}}^{13}$ in the
same way coordinatizing $J$ with 8 parameters. We derive such a coordinatization
of $J$ using three fixed complementary primitive idempotents and
the associated Peirce decomposition of $J$. 

The second aim of this paper is to construct complex prime $\mathbb{Q}$-Fano
$3$-folds anticanonically embedded of codimension 4 as weighted complete
intersections in suitable weighted projectivizations of $\mathscr{H}_{\mathbb{A}}^{13}$,
its subvarieties or their weighted cones (possibly allowing some coordinates
of weight zero). The affine variety $\mathscr{H}_{\mathbb{A}}^{13}$
and such weighted projectivizations of $\mathscr{H}_{\mathbb{A}}^{13}$
are called key varieties for prime $\mathbb{Q}$-Fano 3-folds. We
show that a prime $\mathbb{Q}$-Fano 3-fold of genus $3$ with three
$\nicefrac{1}{2}(1,1,1)$-singularities of type No.$\,$5.4 as in
\cite{Tak1} is constructed as a linear section of a weighted projectivization
of $\mathscr{H}_{\mathbb{A}}^{13}$ such that all the coordinates
have positive weights. Conversely, we also show that any such prime
$\mathbb{Q}$-Fano 3-fold can be obtained in this way. Moreover, relating
$\mathscr{H}_{\mathbb{A}}^{13}$ to the $C_{2}$-cluster variety constructed
by Coughlan and Ducat \cite{CD1}, we show that weighted projectivizations
of $\mathscr{H}_{\mathbb{A}}^{13}$ or its subvarieties serve as key
varieties for prime $\mathbb{Q}$-Fano 3-folds belonging to 108 classes
in the online database \cite{GRDB}. 
\end{abstract}

\maketitle
\markboth{$\mQ$-Fano 3-fold and Jordan algebra }{Hiromichi Takagi}
{\small{\tableofcontents{}}}{\small\par}

2020\textit{ Mathematics subject classification}: 14J45, 14E30, 17C10,
17C50


\section{\textbf{Introduction\label{sec:Introduction}}}

\subsection{Background for classification of $\mQ$-Fano 3-folds\label{subsec:Background-for-classification}}

A complex projective variety is called a \textit{$\mQ$-Fano variety}
if it is a normal variety with only terminal singularities, and its
anticanonical divisor is ample. The classification of $\mQ$-Fano
3-folds is one of the central problems in Mori theory for projective
3-folds. Although the classification is far from completion, several
systematic classification results have been obtained so far. In the
online database \cite{GRDB}, a huge table of candidates for $\mQ$-Fano
3-folds is given. 

To proceed with the classification, we adopt the viewpoint of anticanonical
embeddings into weighted projective spaces. We say that a (not necessarily
prime) $\mQ$-Fano $3$-fold $X$ is \textit{anticanonically embedded
of codimension $c$} in a weighted projective space $\mP$ if the
following conditions are satisfied:

\begin{enumerate}[(1)]

\item $c=\dim\mP-\dim X$, 

\item the anticanonical sheaf $\omega_{X}^{[-1]}$ of $X$ coincides
with $\sO_{\mP}(1)|_{X}$, and

\item the homogeneous coordinate ring of $X$ is Gorenstein.

\end{enumerate} In particular, these conditions imply that the homogeneous
coordinate ring of $X$ coincides with its anticanonical ring $\bigoplus_{n\geq0}H^{0}(X,\omega_{X}^{[-n]})$
(cf.$\,$\cite[(5.16) (ii)]{GW}).

For prime $\mQ$-Fano 3-folds $X$ anticanonically embedded in codimension
$\leq2$, the classification has been completed, provided that $X$
is \textit{quasi-smooth}---that is, its affine cone in the affine
space associated with the anticanonical ring is smooth outside the
vertex (\cite{Fl}, \cite{CCC}). Examples of prime $\mQ$-Fano 3-folds
anticanonically embedded of codimension $3$ were constructed systematically
in \cite{Al} using $5\times5$ skew-symmetric matrices. 

\subsection{Key varieties for prime $\mQ$-Fano 3-folds anticanonically embedded
of codimension 4\label{subsec:Key-varieties-for}}

Under this circumstance, it is meaningful to construct examples of
prime $\mQ$-Fano 3-fold anticanonically embedded of codimension 4
systematically. In \cite{CD1,CD2}, Coughlan and Ducat do this as
weighted complete intersections in weighted projectivizations of the
cluster varieties of types $C_{2}$ and $G_{2}^{(4)}$, which are
related to $\mP^{2}\times\mP^{2}$- and $\mP^{1}\times\mP^{1}\times\mP^{1}$-fibrations,
respectively. This type of construction is modeled after the classification
of smooth prime Fano 3-folds by Gushel\'\ \cite{Gu1,Gu2} and Mukai
\cite{Mu2}. Their idea is to describe $\mQ$-Fano 3-folds $X$ in
appropriate projective varieties $\Sigma$ of larger dimensions such
that usually their affine cones have larger symmetries; in the case
of \cite{Mu2}, $\Sigma$ are rational homogeneous spaces, and in
the case of \cite{CD1,CD2}, $\Sigma$ are several weighted projectivizations
of the cluster varieties. Such $\Sigma$ or the affine cones over
$\Sigma$ are informally called \textit{key varieties} as in the title
of this paper. 

The purpose of this paper is to construct a new key variety $\mathscr{H}_{\mA}^{13}$,
related to $\mP^{2}\times\mP^{2}$-fibration with relative Picard
number $1$, for prime $\mQ$-Fano 3-folds $X$ anticanonically embedded
of codimension $4$. The variety $\mathscr{H}_{\mA}^{13}$ naturally
grows out of the theory of the quadratic Jordan algebra of a cubic
form.

\subsection{Quadratic Jordan algebra of a cubic form and its coordinatization\label{subsec:Jordan-algebra-of}}

In this subsection, we work over an algebraically closed field $\mathsf{k}$.
Let $J$ be the 9-dimensional nondegenerate quadratic Jordan algebra
of a cubic form over $\mathsf{k}$ with reasonable assumptions (for
details, we refer to Assumption (a)--(c) in Section \ref{Sect:Jordan-Algebra-of}).
The quadratic Jordan algebra $J$ is endowed with a quadratic mapping
called $\sharp$-mapping, and a bi-linear product called $\sharp$-product
defined by the $\sharp$-mapping (Section \ref{Sect:Jordan-Algebra-of}).
By the general theory of quadratic Jordan algebras of cubic forms,
$J$ has three complementary primitive idempotents $\upsilon_{1},\upsilon_{2},\upsilon_{3}$,
which give the so-called Peirce decomposition of $J$: 
\[
J=J_{11}\oplus J_{22}\oplus J_{33}\oplus J_{12}\oplus J_{13}\oplus J_{23},
\]
where $J_{ii}=\mathsf{k\upsilon}_{i}\,(1\leq i\leq3)$ and $J_{ij}\simeq\mathsf{k}^{2}\,(1\leq i<j\leq3)$.
Taking a suitable basis of each of the subspaces $J_{ij}$, we show
that the $\sharp$-products between elements of these bases can be
described by introducing $8$ parameters indexed as $p_{ijk}\,(i,j,k=1,2)$
(Proposition \ref{prop:paramH}). 

\vspace{5pt}

Hereafter in Section \ref{sec:Introduction}, we assume for simplicity
that $\mathsf{k}=\mC$, the complex number field. 

\subsection{Affine variety $\mathscr{H}_{\mA}^{13}$}

It is well-known that the affine cone $C(\mP^{2}\times\mP^{2})$ over
(the Segre embedded) $\mP^{2}\times\mP^{2}$ is defined as the null
locus of the $\sharp$-mapping of $J$. Inspired with this fact, we
construct the affine variety $\mathscr{H}_{\mA}^{13}$ in the same
way as $C(\mP^{2}\times\mP^{2})$ coordinatizing $J$ with 8 parameters
$p_{ijk}\,(i,j,k=1,2)$. 

Let $\mA_{\mathscr{H}}$ be the $17$-dimensional affine space extending
$J\simeq\mA^{9}$ by adding $p_{ijk}$ as additional $8$ coordinates.
Then we define $\mathscr{H}_{\mA}^{13}\subset\mA_{\mathscr{H}}^{17}$
as the null locus of the $\sharp$-mapping; we refer to Corollary
\ref{cor:HsharpN} and Definition \ref{def:H13} for the explicit
equations of $\mathscr{H}_{\mA}^{13}$. Here we summarize properties
of $\mathscr{H}_{\mA}^{13}$ as follows:
\begin{thm}
\label{thm:HUpPi} The affine variety $\mathscr{H}_{\mA}^{13}$ has
the following properties:

\begin{enumerate}[$(1)$]

\item There exists a $({\rm GL}_{2})^{3}\rtimes\mathfrak{S}_{3}$-action
on $\mathscr{H}_{\mA}^{13}$(Proposition \ref{prop:actionH}).

\item Let $\mA_{\mathsf{P}}$ be the affine space with coordinates
$p_{ijk}\,(i,j,k=1,2)$, and $\widehat{\mathscr{H}}$ the variety
in $\mP(J)\times\mA_{\mathsf{P}}$ with the same equations as $\mathscr{H}_{\mA}^{13}$.
The natural projection $\rho_{\mathscr{H}}\colon\widehat{\mathscr{H}}\to\mA_{\mathsf{P}}$
is equidimensional outside the origin of $\mA_{\mathsf{P}}$, and
a general $\rho_{\mathscr{H}}$-fiber is isomorphic to $\mP^{2}\times\mP^{2}$
(Proposition \ref{prop:Let--beP2P2fib}).

\item The variety $\mathscr{H}_{\mA}^{13}$ is a $13$-dimensional
irreducible normal variety with only factorial Gorenstein terminal
singularities (Propositions \ref{prop:The-singular-locus}, \ref{prop:9times16 H},
\ref{prop:UFD} and \ref{prop:HTerm}).

\item The ideal of $\mathscr{H}_{\mA}^{13}$ is generated by $9$
elements with $16$ relations (Proposition \ref{prop:9times16 H}).

\end{enumerate}
\end{thm}

In Section \ref{sec:Relatives-of-=002013AffineMS}, we also introduce
two subvarieties $\mathscr{M}_{\mA}^{8}$ and $\mathscr{S}_{\mA}^{6}$
of $\mathscr{H}_{\mA}^{13}$ with $({\rm GL}_{2})^{2}$- and ${\rm GL}_{3}$-action
respectively, and admitting $\mP^{2}\times\mP^{2}$-fibrations, which
turn out to be key varieties for many prime $\mQ$-Fano 3-folds (Section
\ref{sec:The--cluster-variety}).

\subsection{$\mathbb{Q}$-Fano $3$-fold of genus $3$ with three $\nicefrac{1}{2}(1,1,1)$-singularities\label{subsec:-Fano--fold-with}}
\begin{defn}
Let $\Sigma$ be a projective variety in a weighted projective space.
We say that a subvariety $X\subset\Sigma$ is a \textit{linear section
}of $\Sigma$ if $X$ is defined in $\Sigma$ by weighted homogeneous
hypersurfaces $H_{1},\dots,H_{c}$ of weight $1$, where $c$ is the
codimension of $X$ in $\Sigma$.
\end{defn}

Let $\mathscr{H}_{\mP}^{12}\subset\mP(1^{14},2^{3})$ be the weighted
projectivization of $\mathscr{H}_{\mA}^{13}$ such that the 3 weights
of the coordinates corresponding to $\upsilon_{i}\,(1\leq i\leq3)$
are $2$, and the other weights of coordinates are 1. We prove the
following result in Section \ref{sec:-Fano--fold-of 5.4}:
\begin{thm}
\label{thm:embthm} The following assertions hold:

\begin{enumerate}[$(1)$]

\item For general $9$ weighted homogeneous hypersurfaces $H_{1},\dots,H_{9}$
of weight $1$ in $\mP(1^{14},2^{3})$, $\mathscr{H}_{\mP}^{12}\cap H_{1}\cap\dots\cap H_{9}$
is a linear section of $\mathfrak{\mathscr{H}}_{\mP}^{12}$ of codimension
$9$ and has only three $\nicefrac{1}{2}(1,1,1)$-singularities. A
linear section of $\mathfrak{\mathscr{H}}_{\mP}^{12}$ of codimension
$9$ with only three $\nicefrac{1}{2}(1,1,1)$-singularities is a
prime $\mQ$-Fano $3$-fold, anticanonically embedded of codimension
$4$, and is of type No.$\,5.4$ as in \cite{Tak1}, regardless of
the singularity specified.

\item Conversely, any prime $\mathbb{Q}$-Fano $3$-fold with only
$\nicefrac{1}{2}(1,1,1)$-singularities of type No.$\,5.4$ for a
$\nicefrac{1}{2}(1,1,1)$-singularity is a codimension $9$ linear
section of $\mathscr{H}_{\mathbb{P}}^{12}$, and is anticanonically
embedded of codimension $4$. 

\end{enumerate}
\end{thm}

We will see in the course of the proof of Theorem \ref{thm:embthm}
(2) that the affine variety $\mathscr{H}_{\mathbb{A}}^{13}$ emerges
naturally also from birational geometry of the $\mQ$-Fano $3$-fold.
We mention that we obtain similar results for other prime $\mQ$-Fano
3-folds in \cite{Tak3}. 

\subsection{Relation to the $\mathsf{C}_{2}$-cluster variety\label{subsec:Relation-with-the}}

As is mentioned in Subsection \ref{subsec:Key-varieties-for}, several
weighted projectivizations of the $\mathsf{C}_{2}$-cluster variety
are key varieties for prime Fano 3-folds by \cite{CD1,CD2}. We show
that the $\mathsf{C}_{2}$-cluster variety is a subvariety of $\mathscr{H}_{\mA}^{13}$
(Proposition \ref{prop:C2}). Then this implies that, corresponding
to each of the $108$ classes of prime $\mQ$-Fano 3-folds constructed
by ibid., we can construct $\mQ$-Fano 3-folds from weighted projectivizations
of 4 subvarieties $\mathscr{H}_{\mA}^{12},\mathscr{H}_{\mA}^{11}$,
$\mathscr{M}_{\mA}^{8}$, $\mathscr{S}_{\mA}^{6}$ of $\mathscr{H}_{\mA}^{13}$,
which are more general than the ones constructed by ibid (Corollary
\ref{cor:example}). 

\subsection{Relation to our other works \label{subsec:Relation-with-our}}

We construct key varieties for prime $\mQ$-Fano $3$-folds in our
papers \cite{Tak4,Tak5,Tak6}. Here we clarify the relations of them
with the present paper. 

In the papers \cite{Tak5,Tak6}, we construct key varieties related
to $\mP^{2}\times\mP^{2}$-fibration; the varieties $\Upsilon_{\mA}^{14}$
and $\Pi_{\mA}^{15}$ are constructed in \cite{Tak5} via other coordinatizations
of the quadratic Jordan algebras of cubic forms (see also \cite{Tay}
as for $\Upsilon_{\mA}^{14}$, where it is constructed by the theory
of unprojection). The variety $\Sigma_{\mA}^{14}$ is constructed
in \cite{Tak6}. Examples of prime $\mQ$-Fano 3-folds are constructed
from $\Sigma_{\mA}^{14}$ and $\Pi_{\mA}^{15}$ in \cite{Tak7}.

In the paper \cite{Tak4}, we construct key varieties for prime $\mQ$-Fano
3-folds related to $\mP^{1}\times\mP^{1}\times\mP^{1}$-fibration,
and study the relation of them with the $G_{2}^{(4)}$-cluster variety.
Our constructions of them are based on the theory of Freudenthal triple
system. 

\subsection{Structure of the paper}

Sections \ref{Sect:Jordan-Algebra-of} and \ref{sec:CoordJ}, and
Subsection \ref{subsec:Nondegeneracy-of} are devoted to studying
the quadratic Jordan algebra of a cubic form, and, in the remaining
part, we are concerned with the classification of prime $\mQ$-Fano
$3$-folds. 

After we review general facts about the quadratic Jordan algebras
of cubic forms in Section \ref{Sect:Jordan-Algebra-of} and Subsection
\ref{subsec:Facts-about-three}, we obtain a coordinatization with
8 parameters of the 9-dimensional nondegenerate quadratic Jordan algebra
$J$ of a cubic form with some reasonable assumptions (Proposition
\ref{prop:paramH}). In Subsection \ref{subsec:-mapping-and-CubicH},
we write down explicitly the $\sharp$-mapping and the cubic form
of $J$ using $8$ parameters (Corollary \ref{cor:HsharpN}). This
leads us to our definition of the affine variety $\mathscr{H}_{\mA}^{13}$
(Definition \ref{def:H13}). It turns out that $\mathscr{H}_{\mA}^{13}$
has a natural $({\rm GL}_{2})^{3}\rtimes\mathfrak{S}_{3}$-action,
and then is related to $2\times2\times2$ hypermatrix (Subsections
\ref{subsec:-action-onH} and \ref{subsec:Cayley's-hyperdeterminant}).
In the remaining part of Section \ref{sec:Affine-varietyH}, we show
Theorem \ref{thm:HUpPi}. Section \ref{sec:-Fano--fold-of 5.4} is
occupied with the proof of Theorem \ref{thm:embthm}. In Section \ref{sec:Relatives-of-=002013AffineMS},
we introduce two subvarieties $\mathscr{M}_{\mA}^{8}$ and $\mathscr{S}_{\mA}^{6}$
of $\mathscr{H}_{\mA}^{13}$ with good group actions. In the final
section \ref{sec:The--cluster-variety}, we show that the $C_{2}$-cluster
variety $X_{\mathsf{C}_{2}}$ is a subvariety of $\mathscr{H}_{\mA}^{13}$.
Then, based on this fact, We also clarify relations of $X_{\mathsf{C}_{2}}$
with 4 subvarieties $\mathscr{H}_{\mA}^{12}$, $\mathscr{H}_{\mA}^{11}$,
$\mathscr{M}_{\mA}^{8}$, $\mathscr{S}_{\mA}^{6}$ of $\mathscr{H}_{\mA}^{13}$.
Then we show that these 4 varieties are key varieties of prime $\mQ$-Fano
3-folds belonging to 108 classes (Corollary \ref{cor:example}). 

\vspace{10pt}

\noindent\textbf{ Notation and Conventions:}

\begin{itemize}

\item $\mP(x_{1},\dots,x_{n})$ and $\mA(x_{1},\dots,x_{n})$: the
projective space and the affine space with coordinates $x_{1},\dots,x_{n}$,
respectively.

\item The $x$-point of a weighted projective space $\mP$$:=$the
point of $\mP$ whose only nonzero coordinate is $x$. 

\item \textit{Divisor and sheaf convention}: We denote a Weil divisor
on a normal variety by the same notation as the corresponding divisorial
sheaf. Although the notation feels somewhat awkward, we think that
it is better than unnecessarily increasing the number of symbols.

\end{itemize}

\vspace{3pt}

This paper contains several assertions which can be proved by straightforward
computations; we often omit such computation processes. Some computations
are difficult by hand but are easy within a software package. In our
computations, we use intensively the software systems Mathematica
\cite{W} and \textsc{Singular} \cite{DGPS}. 

\vspace{3pt}

\noindent\textbf{Acknowledgment:} I am very grateful to Professor
Shigeru Mukai; inspired by his article \cite{Mu1}, I was led to Jordan
algebra to describe key varieties. Moreover, I was helped to find
the key variety $\mathscr{H}_{\mA}^{13}$ by his inference as for
dimensions of key varieties. I owe many important calculations in
the paper to Professor Shinobu Hosono. I wish to thank him for his
generous cooperations. This work is supported in part by Grant-in
Aid for Scientific Research (C) 16K05090. 

\section{\textbf{Quadratic Jordan algebra $J$ of cubic forms}\label{Sect:Jordan-Algebra-of}}

In this section, we review the construction of the Jordan algebra
of a cubic form and collect some basic facts following \cite{R,Mc2}.

Let $\mathsf{k}$ be a field and $J$ a finite dimensional $\mathsf{k}$-vector
space and $N(x)\,(x\in J)$ a cubic form on $J$. We fix an element
$\mathfrak{1}\in J$ such that 
\[
N(\mathfrak{1})=1.
\]

We denote by $\partial_{y}N(x)$ the directional derivative of $N$
in the direction $y\in J$, evaluated at $x\in J$.

From the cubic form $N(x)$, the following several forms are derived:

\[
\begin{cases}
\text{Bilinear Trace Form} & :T(x,y):=-\partial_{x}\partial_{y}N(\mathfrak{1})+(\partial_{x}N(\mathfrak{1}))(\partial_{y}N(\mathfrak{1}))\\
\text{Linear Trace Form} & :T(x):=T(x,\mathfrak{1}).
\end{cases}
\]
Since $N(x)$ is a cubic form, we have $\partial_{1}N(\mathfrak{1})=3$
and $\partial_{x}\partial_{1}N(\mathfrak{1})=\partial_{1}\partial_{x}N(\mathfrak{1})=2\partial_{x}N(\mathfrak{1})$
by the Euler formula and the assumption that $N(\mathfrak{1})=1$.
Therefore we obtain
\[
T(x)=T(x,\mathfrak{1})=-2\partial_{x}N(\mathfrak{1})+3\partial_{x}N(\mathfrak{1})=\partial_{x}N(\mathfrak{1}),
\]
and
\[
T(\mathfrak{1})=\partial_{\mathfrak{1}}N(\mathfrak{1})=3.
\]
 Now $\sharp$-mapping and $\sharp$-product are introduced as follows,
which will guarantee that $J$ has the structure of quadratic Jordan
algebra defined from the cubic form $N(x)$: 
\[
\begin{cases}
\text{\ensuremath{\sharp}-mapping} & :x\mapsto x^{\sharp},\\
\sharp\text{-product} & :x\sharp y:=(x+y)^{\sharp}-x^{\sharp}-y^{\sharp}
\end{cases}
\]
satisfying the following three conditions: 

\vspace{5pt}

\noindent\textbf{Sharp Condition:}
\begin{align*}
(\sharp1) & \quad T(x^{\sharp},y)=\partial_{y}N(x),\\
(\sharp2) & \quad x^{\sharp\sharp}=N(x)x,\\
(\sharp3) & \mathfrak{\quad1}\sharp y=T(y)\mathfrak{1}-y.
\end{align*}
In this paper, the $\sharp$-product plays important roles. 

The following fundamental theorem is due to McCrimmon (\cite[Const. 4.2.2]{Mc2}):
\begin{thm}
\label{thm:Jordan}Assume that the cubic form $N(x)$ and $\mathfrak{1}\in J$
are endowed with the $\sharp$-mapping and the $\sharp$-product satisfying
the conditions $(\sharp1)$-$(\sharp3)$. Then $J$ has the structure
of quadratic Jordan algebra with $\mathfrak{1}$ as the unit element
and with the following $U$-operator:
\[
U_{x}y:=T(x,y)x-x^{\sharp}\sharp y.
\]
\end{thm}

We refer to \cite[Sect.4.3]{Mc2} and \cite{R} for the definition
of quadratic Jordan algebra.

\begin{defn}
We call the quadratic Jordan algebra defined as in Theorem \ref{thm:Jordan}
\textit{the quadratic Jordan algebra of the cubic form $N(x).$ }
\end{defn}

The introduction of the forms for $J$ continues:

\begin{align*}
\begin{cases}
\text{Quadratic Spur form} & :S(x):=T(x^{\sharp})\\
\text{Bilinear Spur form} & :S(x,y):=T(x\sharp y)=T(x,y)-T(x)T(y),
\end{cases}
\end{align*}
where the last equality of the second line holds by \cite[p.497, (16)]{Mc1}.
\begin{rem}
\label{rem:In-the-situationJprod}In the situation of Theorem \ref{thm:Jordan},
if the characteristic of $\mathsf{k}$ is not 2, then the usual Jordan
product $\empty\bullet\empty$ is defined as follows (\cite[Const. 4.2.2]{Mc2}):
\[
x\bullet y:=\frac{1}{2}(x\sharp y+T(x)y+T(y)x-S(x,y)\mathfrak{1}).
\]

The product $x\bullet x$ can be defined for any $\mathsf{k}$ as
follows:
\begin{equation}
x\bullet x:=x^{\sharp}+T(x)x-S(x)\mathfrak{1}.\label{eq:xx}
\end{equation}
\end{rem}

Now the following nondegeneracy condition for the quadratic Jordan
algebra $J$ of a cubic form is introduced, which will guarantee a
good structure theorem of $J$:
\begin{defn}
We set 
\begin{equation}
R_{J}:=\{\sigma\in J\mid\forall_{y\in J},\,U_{\sigma}y=0\}.\label{eq:nondeg}
\end{equation}
The quadratic Jordan algebra $J$ of the cubic form $N(x)$ is said
to be \textit{nondegenerate }if $R_{J}=\{0\}$. 
\end{defn}

We note that a nonzero element of $R_{J}$ is called in \cite{R}
\textit{an absolute zero divisor.} 

In \cite[Cor.4.1.58]{Pe}, the following effective paraphrase of nondegeneracy
is obtained:
\begin{thm}
The quadratic Jordan algebra $J$ of the cubic form $N(x)$ is nondegenerate
if and only if it holds that
\[
\{\sigma\in J\mid N(\sigma)=T(\sigma,J)=T(\sigma^{\sharp},J)=0\}=\{0\}.
\]
\end{thm}

\vspace{5pt}

Hereafter in this paper, we consider the quadratic Jordan algebra
$J$ of the cubic form $N(x)$ with the following assumptions:

\vspace{5pt}

\noindent \textbf{Assumptions.}

\begin{enumerate}[(a)]

\item $J$ is nondegenerate,

\item $N(x)\not=0$ for some $x\not=0$, and

\item $J\not\simeq\mathsf{k}\oplus J(Q,\mathfrak{1}'),$ where $J(Q,\mathfrak{1}')$
is the Jordan algebra of a quadratic form $Q$ with the identity $\mathfrak{1}'$.

\end{enumerate}

These assumptions originate in the paper \cite{R}, where Racine obtains
a coordinatization of the quadratic Jordan algebra $J$ of a cubic
form with the assumptions (a)--(c) ({[}ibid., Thm. 1.1{]}); he describe
$J$ using a $3\times3$ matrix over a composition algebra. We do
not state his result in detail since we do not need the result itself
in this paper. We, however, heavily use his arguments and several
formulas contained in the proof of {[}ibid., Thm. 1.1{]}, which we
will recall in Subsection \ref{subsec:Facts-about-three}.

\section{\textbf{Coordinatization of $J$ with $8$ parameters \label{sec:CoordJ}}}

\subsection{Facts about three complementary primitive idempotents\label{subsec:Facts-about-three}}

Let $\mathsf{k}$ be a field and $J$ be the quadratic Jordan algebra
of a cubic form with the assumptions (a)--(c) as in Section \ref{Sect:Jordan-Algebra-of}.
A nonzero element $e\in J$ is called an \textit{idempotent }if it
satisfies $e\bullet e=e$ (cf.~Remark \ref{rem:In-the-situationJprod}).
An idempotent is called\textit{ primitive }if it also satisfies $e^{\sharp}=0$.
For $x_{1},x_{2},y\in J$, we define $U_{x_{1},x_{2}}(y):=U_{\mathfrak{1}}(y)-U_{x_{1}}(y)-U_{x_{2}}(y).$

By \cite{R}, there exist three primitive idempotents $\upsilon_{1},\,\upsilon_{2},\,\upsilon_{3}$
such that $\upsilon_{1}+\upsilon_{2}+\upsilon{}_{3}=\mathfrak{1}$.
Such three primitive idempotents are called\textit{ complementary}.
They give the following decomposition of $J$, which is called the
\textit{Peirce decomposition} (associated to $\upsilon_{1},\,\upsilon_{2},\,\upsilon_{3}$):
\[
J=J_{11}\oplus J_{22}\oplus J_{33}\oplus J_{12}\oplus J_{13}\oplus J_{23},
\]
where 
\[
J_{ii}=\mathsf{k}\upsilon_{i},J_{ij}=U_{\upsilon_{i},\upsilon_{j}}(J)\,(1\leq i<j\leq3).
\]
If the characteristic of $\mathsf{k}$ is not 2, then it holds that
\[
J_{ij}=\left\{ \sigma\in J\mid\upsilon_{i}\bullet\sigma=\upsilon_{j}\bullet\sigma=\nicefrac{1}{2}\,\sigma\right\} \,(1\leq i<j\leq3).
\]
We define for convenience that $J_{12}=J_{21}$, $J_{13}=J_{31}$
and $J_{23}=J_{32}$. The subspaces $J_{ij}$ are called \textit{the
Peirce subspaces}. In \cite[p.99]{R} , it is proved that 
\begin{equation}
\dim J_{12}=\dim J_{13}=\dim J_{23}.\label{eq:J12J13J23}
\end{equation}

The following facts are obtained in \cite{R} in the course of the
proof of {[}ibid., Thm.1{]}:
\begin{thm}
\label{thm:The-Jordan-algebra racine} Let $i,j,k$ be three integers
such that $\{i,j,k\}=\{1,2,3\}.$ The Peirce spaces have the following
properties:

\begin{enumerate}

\item[$(\rm{i})$] For $\sigma\in J_{jk}$, it holds that $T(\sigma)=0$
$($\cite[(28)]{R}$)$.

\item[$(\rm{ii})$] $J_{ij}$ and $J_{ik}$ are orthogonal with respect
to the bilinear trace form $T$ $($\cite[(33)]{R}$)$.

\item[$(\rm{iii})$] For $\upsilon_{1},$ $\upsilon_{2}$, $\upsilon_{3}$,
and $\sigma\in J_{jk}$, it holds that 
\begin{align*}
\sigma^{\sharp}=S(\sigma)\upsilon_{i},\ \upsilon_{i}\sharp\upsilon_{j}=\upsilon_{k},\ \upsilon_{j}\sharp\sigma=\upsilon_{k}\sharp\sigma=0,\ \upsilon_{i}\sharp\sigma=-\sigma
\end{align*}
 $($\cite[(33)]{R}$)$.

\item[$(\rm{iv})$] $J_{jk}\not=\left\{ 0\right\} $ and the quadratic
form $S(\sigma)\,(\sigma\in J_{jk})$ is not identically zero on $J_{jk}$
$($\cite[p.98, 5th line from the bottom]{R}$)$. 

\item[$(\rm{v})$] For $\sigma_{a}\in J_{ij}$ and $\sigma_{b}\in J_{jk}$,
it holds that $\sigma_{a}\sharp\sigma_{b}\in J_{ik}$ $($\cite[(34)]{R}$)$.

\item[$(\rm{vi})$] For $\sigma_{a}\in J_{ij}$ and $\sigma_{b}\in J_{jk}$,
it holds that 
\begin{equation}
\sigma_{a}\sharp(\sigma_{a}\sharp\sigma_{b})=-S(\sigma_{a})\sigma_{b}\label{eq:xxy}
\end{equation}
$($\cite[(36)]{R}$)$.

\end{enumerate}

\vspace{0.5cm}
\end{thm}

Hereafter in this section, we keep considering the quadratic Jordan
algebra $J$ of a cubic form with the assumptions (a)--(c), and fixing
the three complementary primitive idempotents $\upsilon_{1},\upsilon_{2},\upsilon_{3}$.\textit{
Moreover, we assume that 
\[
\dim J=9.
\]
}Then, by (\ref{eq:J12J13J23}), it holds that \textit{
\[
\dim J_{12}=\dim J_{13}=\dim J_{23}=2.
\]
}

\subsection{Coordinatization of $J$ using $3$ complementary primitive idempotents\label{subsec:Coordinatization-of-thrre idem}}

We consider how to describe the cubic form $N$ from the Peirce decomposition.
For this, it suffices to compute $x^{\sharp}$ for any $x$ since
$x^{\sharp\sharp}=N(x)x.$ We can compute $x^{\sharp}$ for any $x$
if we know all the $\sharp$ products between the elements of the
bases of the Peirce spaces since the $\sharp$ product is bilinear.
By Theorem \ref{thm:The-Jordan-algebra racine} (iii), we already
know the $\sharp$ products among $\upsilon_{1},\upsilon_{2},\upsilon_{3}$,
and an element $x$ of a Peirce space. Hence it remains to know the
$\sharp$ products $x\sharp y$ when $x\in J_{ij}$ and $y\in J_{jk}$,
or $x,\,y\in J_{ij}$. In the following proposition, assuming that
$\mathsf{k}$ is algebraically closed, we show that everything is
determined by taking an appropriate bases for $J_{12},J_{13},J_{23}$
and giving 8 parameters that define the $\sharp$ products between
the bases of $J_{23}$ and $J_{13}$ as in (\ref{eq:J12P}) below. 
\begin{prop}
\label{prop:paramH} Assume that $\mathsf{k}$ is algebraically closed.
Then there exist bases $\chi_{1i}$, $\chi_{2i}$ of $J_{jk}$ for
each $i,j,k$ with $\{i,j,k\}=\{1,2,3\}$ and $p_{abc}\in\mathsf{k}\,(a,b,c=1,2)$
such that the following relations (\ref{eq:J12P})--(\ref{eq:JiiP})
hold: 
\begin{equation}
\begin{cases}
\chi_{11}\sharp\chi_{12}=p_{221}\chi_{13}+p_{222}\chi_{23}, & \chi_{11}\sharp\chi_{22}=-p_{211}\chi_{13}-p_{212}\chi_{23},\\
\chi_{21}\sharp\chi_{12}=-p_{121}\chi_{13}-p_{122}\chi_{23}, & \chi_{21}\sharp\chi_{22}=p_{111}\chi_{13}+p_{112}\chi_{23}.
\end{cases}\label{eq:J12P}
\end{equation}

\begin{equation}
\begin{cases}
\chi_{11}\sharp\chi_{13}=p_{212}\chi_{12}+p_{222}\chi_{22}, & \chi_{11}\sharp\chi_{23}=-p_{211}\chi_{12}-p_{221}\chi_{22},\\
\chi_{21}\sharp\chi_{13}=-p_{112}\chi_{12}-p_{122}\chi_{22}, & \chi_{21}\sharp\chi_{23}=p_{111}\chi_{12}+p_{121}\chi_{22}.
\end{cases}\label{eq:J13P}
\end{equation}
\begin{equation}
\begin{cases}
\chi_{12}\sharp\chi_{13}=p_{122}\chi_{11}+p_{222}\chi_{21}, & \chi_{12}\sharp\chi{}_{23}=-p_{121}\chi_{11}-p_{221}\chi_{21},\\
\chi_{22}\sharp\chi_{13}=-p_{112}\chi_{11}-p_{212}\chi_{21}, & \chi_{22}\sharp\chi_{23}=p_{111}\chi_{11}+p_{211}\chi_{21}.
\end{cases}\label{eq:J23P}
\end{equation}

\begin{equation}
\begin{cases}
\chi_{11}\sharp\chi_{11}=2(-p_{212}p_{221}+p_{211}p_{222})\upsilon_{1},\\
\chi_{11}\sharp\chi_{21}=(-p_{122}p_{211}+p_{121}p_{212}+p_{112}p_{221}-p_{111}p_{222})\upsilon_{1},\\
\chi_{21}\sharp\chi_{21}=2(-p_{112}p_{121}+p_{111}p_{122})\upsilon_{1},\\
\chi_{12}\sharp\chi_{12}=2(-p_{122}p_{221}+p_{121}p_{222})\upsilon_{2},\\
\chi_{12}\sharp\chi_{22}=(p_{122}p_{211}-p_{121}p_{212}+p_{112}p_{221}-p_{111}p_{222})\upsilon_{2}\\
\chi_{22}\sharp\chi_{22}=2(-p_{112}p_{211}+p_{111}p_{212})\upsilon_{2},\\
\chi_{13}\sharp\chi_{13}=2(-p_{122}p_{212}+p_{112}p_{222})\upsilon_{3},\\
\chi_{13}\sharp\chi_{23}=(p_{122}p_{211}+p_{121}p_{212}-p_{112}p_{221}-p_{111}p_{222})\upsilon_{3},\\
\chi_{23}\sharp\chi_{23}=2(-p_{121}p_{211}+p_{111}p_{221})\upsilon_{3}.
\end{cases}\label{eq:JiiP}
\end{equation}
\end{prop}

\begin{rem}
\label{rem:It-is-convenient}It is convenient to put $p_{abc}\in\mathsf{k}\,(a,b,c=1,2)$
on the vertices of the cube as follows and we regard it as a hypermatrix
$\mathsf{P}$:\[
\mathsf{P}:=\vcenter{\hbox{\xymatrix@!0{ & p_{122} \ar@{-}[rr]\ar@{-}'[d][dd] & & p_{222} \ar@{-}[dd] \\ p_{121} \ar@{-}[ur]\ar@{-}[rr]\ar@{-}[dd] & & p_{221} \ar@{-}[ur]\ar@{-}[dd] \\ & p_{112} \ar@{-}'[r]
[rr]& & p_{212} \\ p_{111} \ar@{-}[rr]\ar@{-}[ur] & & p_{211}. \ar@{-}[ur] }}}
\]
\end{rem}

\begin{proof}[Proof of Proposition \ref{prop:paramH}]
For $\sigma\in J_{ij}$, we denote by 
\[
f_{ij}^{\sigma}\colon J_{ik}\to J_{jk}\ \text{and\ }f_{ji}^{\sigma}\colon J_{jk}\to J_{ik}
\]
the linear maps defined by $\sigma\sharp$ (cf. Theorem \ref{thm:The-Jordan-algebra racine}
(v)). We choose a basis $\chi_{1i}$, $\chi_{2i}$ of $J_{jk}$ for
any mutually distinct $i,j,k$. Let $P_{ij}^{\sigma}$ be the representation
matrix of $f_{ij}^{\sigma}$ with respect to these bases. By (\ref{eq:xxy}),
we have 
\begin{equation}
P_{ij}^{\sigma}P_{ji}^{\sigma}=-S(\sigma)E,\label{eq:PijPji}
\end{equation}
where $E$ is the identity matrix. Let $(P_{ij}^{\sigma})^{\dagger}$
be the adjoint matrix of $P_{ij}^{\sigma}$. Multiplying the equality
(\ref{eq:PijPji}) with $(P_{ij}^{\sigma})^{\dagger}$, we have 
\begin{equation}
(\det P_{ij}^{\sigma})P_{ji}^{\sigma}=-S(\sigma)(P_{ij}^{\sigma})^{\dagger}.\label{eq:detPij}
\end{equation}
Note that entries of $P_{ij}^{\sigma}$ are linear forms and then
$\det P_{ij}^{\sigma}$ is a quadratic form with respect to the coordinates
of $\sigma\in J_{ij}$ since the $\sharp$-product is bi-linear and
$\dim J_{jk}=2$.

By Theorem \ref{thm:The-Jordan-algebra racine} (iv), $S(\sigma)\not=0$
for some $\sigma\in J_{ij},$ and then, by (\ref{eq:PijPji}), $\det P_{ij}^{\sigma}\not=0$
for such a $\sigma.$ Now we show that, for some $\alpha_{ij}\in\mathsf{k},$
$\det P_{ij}^{\sigma}=\alpha_{ij}S(\sigma)$ for any $\sigma\in J_{ij}$
using (\ref{eq:detPij}). Assume the contrary. Then, for a reason
of degree, the highest common factor of $S(\sigma)$ and $\det P_{ij}^{\sigma}$
is a linear form, say, $m(\sigma)$. We may write $S(\sigma)=m(\sigma)l(\sigma)$
with a linear form $l(\sigma)$ (note that $l(\sigma)$ does not divide
$\det P_{ij}^{\sigma}$). Then, by (\ref{eq:detPij}), there is a
constant $2\times2$ matrix $M$ such that $P_{ji}^{\sigma}=l(\sigma)M.$
Therefore, for a nonzero $\sigma_{0}\in J_{ij}$ such that $l(\sigma_{0})=0$,
we have $P_{ji}^{\sigma_{0}}=0$. This means $\sigma_{0}\sharp\tau=\bm{0}$
for any $\tau\in J_{jk}.$ For a $\tau_{0}\in J_{jk}$ with $S(\tau_{0})\not=0$
(cf.~Theorem \ref{thm:The-Jordan-algebra racine} (iv)), however,
the linear map $f_{jk}^{\tau_{0}}$ is an isomorphism, thus $\sigma_{0}\sharp\tau_{0}$
cannot be zero since $\sigma_{0}\not=\bm{0}$, a contradiction. Therefore,
for some $\alpha_{ij}\in\mathsf{k},$ $\det P_{ij}^{\sigma}=\alpha_{ij}S(\sigma)$
as desired. 

Note that, by replacing the basis $\chi_{1i}$, $\chi_{2i}$ of $J_{jk}$
with $\beta_{jk}\chi_{1i},\beta_{jk}\chi_{2i}$ for some $\beta_{jk}\in$$\mathsf{k^{*}}$
and the basis $\chi_{1j}$, $\chi_{2j}$ of $J_{ik}$ with $\beta_{ik}\chi_{1j}$,
$\beta_{ik}\chi_{2j}$ for some $\beta_{ik}\in\mathsf{k}^{*}$, the
matrix $P_{ij}^{\sigma}$ is multiplied with $\frac{\beta_{ik}}{\beta_{jk}}$.
Choosing $\beta_{12},\beta_{13},\beta_{23}$ such that $\frac{\beta_{12}}{\beta_{23}}=\sqrt{\alpha_{13}}$
and $\frac{\beta_{12}}{\beta_{13}}=\sqrt{\alpha_{23}}$ (\textit{here
we use the assumption that $\mathsf{k}$ is algebraically closed}),
and taking the corresponding base changes for $J_{12},J_{13},J_{23}$,
we may assume that 
\begin{equation}
\det P_{13}^{\sigma}=S(\sigma)\ \text{and\ }\det P_{23}^{\tau}=S(\tau)\label{eq:SxSy}
\end{equation}
for $\sigma\in J_{13}$ and $\tau\in J_{23}$. Then, by (\ref{eq:detPij}),
we have 
\begin{equation}
P_{31}^{\sigma}=-(P_{13}^{\sigma})^{\dagger}\ \text{and}\ P_{32}^{\tau}=-(P_{23}^{\tau})^{\dagger}\label{eq:P31P32}
\end{equation}
for $\sigma\in J_{13}$ and $\tau\in J_{23}$. 

Now, by setting 
\[
P_{32}^{\chi_{11}}=\left(\begin{array}{cc}
p_{221} & -p_{211}\\
p_{222} & -p_{212}
\end{array}\right)\ \text{and}\ P_{32}^{\chi_{21}}=\left(\begin{array}{cc}
-p_{121} & p_{111}\\
-p_{122} & p_{112}
\end{array}\right),
\]
the equality (\ref{eq:J12P}) holds. We have also the first three
equalities of (\ref{eq:JiiP}) by the first equality of Theorem \ref{thm:The-Jordan-algebra racine}
(iii), (\ref{eq:SxSy}), (\ref{eq:P31P32}), and the definition of
$\sharp$ product. Since $P_{23}^{\tau}=-(P_{32}^{\tau})^{\dagger}$
for $\tau\in J_{23}$ by (\ref{eq:P31P32}), we have (\ref{eq:J13P}).
By the commutativity of $\sharp$ product, (\ref{eq:J12P}) implies
\begin{align*}
\chi_{12}\sharp\chi_{11}=p_{221}\chi_{13}+p_{222}\chi_{23},\  & \chi_{12}\sharp\chi_{21}=-p_{121}\chi_{13}-p_{122}\chi_{23},\\
\chi_{22}\sharp\chi_{11}=-p_{211}\chi_{13}-p_{212}\chi_{23},\  & \chi_{22}\sharp\chi_{21}=p_{111}\chi_{13}+p_{112}\chi_{23}.
\end{align*}
Thus we have $P_{31}^{\chi_{12}}=\left(\begin{array}{cc}
p_{221} & -p_{121}\\
p_{222} & -p_{122}
\end{array}\right)$ and $P_{31}^{\chi_{22}}=\left(\begin{array}{cc}
-p_{211} & p_{111}\\
-p_{212} & p_{112}
\end{array}\right)$. Then, since $P_{13}^{\sigma}=-(P_{31}^{\sigma})^{\dagger}$ for
$\sigma\in J_{13}$, we obtain (\ref{eq:J23P}). Finally we obtain
the last six equalities of (\ref{eq:JiiP}) in a similar way to obtain
the first three equalities of (\ref{eq:JiiP}).
\end{proof}

\subsection{$\sharp$-mapping and cubic form\label{subsec:-mapping-and-CubicH}}

By Proposition \ref{prop:paramH} and direct computations, we obtain
the formulas of $\sigma^{\sharp}$ for any $\sigma\in J$ and the
explicit presentation of the cubic form of $J$ as follows:
\begin{cor}
\label{cor:HsharpN}Under the same assumptions and notation as in
Proposition \ref{prop:paramH}, we also fix the bases $\chi_{1i}$,
$\chi_{2i}$ of $J_{jk}$ for each $i,j,k$ with $\{i,j,k\}=\{1,2,3\}$
and $p_{abc}\in\mathsf{k}\,(a,b,c=1,2)$. Let $\sigma\in J$ be any
element. We write 
\begin{align*}
 & \sigma=\sum_{i=1}^{3}(x_{1i}\chi_{1i}+x_{2i}\chi_{2i})+\sum_{j=1}^{3}u_{j}\upsilon_{j},\\
 & \sigma^{\sharp}=\sum_{i=1}^{3}(x_{1i}^{\sharp}\chi_{1i}+x_{2i}^{\sharp}\chi_{2i})+\sum_{j=1}^{3}u_{j}^{\sharp}\upsilon_{j},
\end{align*}
and set 
\[
\bm{x}_{i}:=\left(\begin{array}{c}
x_{1i}\\
x_{2i}
\end{array}\right),\,\bm{x}_{i}^{\sharp}:=\left(\begin{array}{c}
x_{1i}^{\sharp}\\
x_{2i}^{\sharp}
\end{array}\right)\,(1\leq i\leq3).
\]

For $i,j=1,2$ and $k=1,2,3$, we define
\begin{align*}
D_{ij}^{(1)}(x) & :=\begin{vmatrix}p_{1ij} & x_{11}\\
p_{2ij} & x_{21}
\end{vmatrix},D_{ij}^{(2)}(x):=\begin{vmatrix}p_{i1j} & x_{12}\\
p_{i2j} & x_{22}
\end{vmatrix},D_{ij}^{(3)}(x):=\begin{vmatrix}p_{ij1} & x_{13}\\
p_{ij2} & x_{23}
\end{vmatrix},\\
D_{ij}^{(1)}(x^{\sharp}) & :=\begin{vmatrix}p_{1ij} & x_{11}^{\sharp}\\
p_{2ij} & x_{21}^{\sharp}
\end{vmatrix},D_{ij}^{(2)}(x^{\sharp}):=\begin{vmatrix}p_{i1j} & x_{12}^{\sharp}\\
p_{i2j} & x_{22}^{\sharp}
\end{vmatrix},D_{ij}^{(3)}(x^{\sharp}):=\begin{vmatrix}p_{ij1} & x_{13}^{\sharp}\\
p_{ij2} & x_{23}^{\sharp}
\end{vmatrix},
\end{align*}
\begin{align*}
D^{(k)}(x):= & \begin{pmatrix}-D_{12}^{(k)}(x) & D_{11}^{(k)}(x)\\
-D_{22}^{(k)}(x) & D_{21}^{(k)}(x)
\end{pmatrix},\\
D^{(k)}(x,x^{\sharp})\boldsymbol{:=} & \begin{pmatrix}-D_{12}^{(k)}(x) & D_{11}^{(k)}(x^{\sharp})\\
-D_{22}^{(k)}(x) & D_{21}^{(k)}(x^{\sharp})
\end{pmatrix},D^{(k)}(x^{\sharp},x):=\begin{pmatrix}-D_{12}^{(k)}(x^{\sharp}) & D_{11}^{(k)}(x)\\
-D_{22}^{(k)}(x^{\sharp}) & D_{21}^{(k)}(x)
\end{pmatrix},
\end{align*}

$ $The coordinates of $\sigma^{\sharp}$ and the cubic form defining
the quadratic Jordan structure of $J$ are written as follows:
\begin{equation}
\begin{cases}
\bm{x}_{1}^{\sharp} & =-u_{1}\bm{x}_{1}+D^{(3)}(x)\bm{x}_{2},\\
\bm{x}_{2}^{\sharp} & =-u_{2}\bm{x}_{2}+D^{(1)}(x)\bm{x}_{3},\\
\bm{x}_{3}^{\sharp} & =-u_{3}\bm{x}_{3}-(D^{(2)}(x))^{\dagger}\bm{x}_{1},\\
u_{1}^{\sharp} & =u_{2}u_{3}+|D^{(1)}(x)|,\\
u_{2}^{\sharp} & =u_{3}u_{1}+|D^{(2)}(x)|,\\
u_{3}^{\sharp} & =u_{1}u_{2}+|D^{(3)}(x)|,
\end{cases}\label{eq:HypMat}
\end{equation}
where $(D^{(2)}(x))^{\dagger}$ is the adjoint matrix of $D^{(2)}(x)$.

We denote the r.h.s.~of (\ref{eq:HypMat}) in order from the top
by $\bm{G}_{1},\bm{G}_{2},\bm{G}_{3}$, $\mathsf{G_{4}},\mathsf{G}_{5},\mathsf{G}_{6}$,
where $\bm{G}_{1},\bm{G}_{2},\bm{G}_{3}$ are $2$-dimensional vectors
with polynomial entries, and $\mathsf{G_{4}},\mathsf{G}_{5},\mathsf{G}_{6}$
are polynomials. 

\begin{align}
N_{\mathsf{P}}: & =\nicefrac{1}{3}\big\{ u_{1}(u_{1}^{\sharp})-|D^{(1)}(x,x^{\sharp})|-|D^{(1)}(x^{\sharp},x)|\label{eq:NH}\\
 & \,\,\,\,\,\,\,\,\,\,\,+u_{2}(u_{2}^{\sharp})-|D^{(2)}(x,x^{\sharp})|-|D^{(2)}(x^{\sharp},x)|\nonumber \\
 & \,\,\,\,\,\,\,\,\,\,\,+u_{3}(u_{3}^{\sharp})-|D^{(3)}(x,x^{\sharp})|-|D^{(3)}(x^{\sharp},x)|\},\nonumber 
\end{align}
 where the coefficient $\nicefrac{1}{3}$ in (\ref{eq:NH}) only appears
in the expression and is canceled such that r.h.s. of (\ref{eq:NH})
has only integer coefficients (thus, it works also when the characteristic
of $\mathsf{k}$ is $3$).
\end{cor}

We also have the following converse statement, for which we do not
need to assume that $\mathsf{k}$ is algebraically closed:
\begin{prop}
\label{Prop:reveseJ}Let $J$ be a $9$-dimensional $\mathsf{k}$-vector
space with coordinates $u_{1}$, $u_{2}$, $u_{3}$, $x_{ij}$ $(i=1,2,3,j=1,2)$.
Let $N_{P}$ be the cubic form (\ref{eq:NH}) on $J$ with some hypermatrix
$\mathsf{P}$. We also define $\sigma^{\sharp}$ for $\sigma\in J$
as in (\ref{eq:HypMat}). The vector space $J$ has the structure
of the quadratic Jordan algebra of the cubic form $N_{\mathsf{P}}$
with the $\sharp$-mapping $J\ni\sigma\mapsto\sigma^{\sharp}\in J$. 
\end{prop}

\begin{proof}
It is straightforward to check the sharp condition ($\sharp$1)--($\sharp$3)
for the cubic form $N_{\mathsf{P}}$ and the map $J\ni\sigma\mapsto\sigma^{\sharp}\in J$.
Then the assertion follows from Theorem \ref{thm:Jordan}.
\end{proof}
\begin{defn}
We denote by $J_{\mathsf{P}}$ the quadratic Jordan algebra of the
cubic form $N_{\mathsf{P}}$ constructed in Proposition \ref{Prop:reveseJ}. 
\end{defn}

$J_{\mathsf{P}}$ is not necessarily nondegenerate. We will revisit
the problem to determine nondegeneracy of $J_{\mathsf{P}}$ in Subsection
\ref{subsec:Nondegeneracy-of}.

\section{\textbf{Affine variety $\mathscr{H}_{\mA}^{13}$\label{sec:Affine-varietyH}}}

Let $\mathsf{k}$ be a field.

\subsection{Definition of $\mathscr{H}_{\mA}^{13}$\label{subsec:Definition-of-theH}}
\begin{defn}
\label{def:H13} Let $\mA_{\mathscr{H}}$ be the $\mathsf{k}$-affine
$17$-space with coordinates $u_{1},u_{2},u_{3},x_{ij}\,(i=1,2,3,j=1,2)$
and $p_{ijk}\,(i,j,k=1,2,3).$ We define a morphism $\mA_{\mathscr{H}}\to\mA_{\mathscr{H}}$
by $\mA_{\mathscr{H}}\ni\sigma\mapsto\sigma^{\sharp}\in\mA_{\mathscr{H}}$
as in Corollary \ref{cor:HsharpN}. In $\mA_{\mathscr{H}}$, we define
$\mathscr{H}_{\mA}^{13}$ to be the affine scheme defined by the condition
\[
\sigma^{\sharp}=0,
\]
namely, the vanishing of the r.h.s. of (\ref{eq:HypMat}). 
\end{defn}

Finally, we show that $\mathscr{H}_{\mA}^{13}$ is a variety, i.e.,
it is irreducible (Proposition \ref{prop:9times16 H} (3)).

\subsection{$({\rm GL}_{2})^{3}\rtimes\mathfrak{S}_{3}$-action on $\mathscr{H}_{\mathbb{A}}^{13}$\label{subsec:-action-onH}}

In this subsection, we show the following: 
\begin{prop}
\label{prop:actionH}$\mathscr{H}_{\mathbb{A}}^{13}$ has a $({\rm GL}_{2})^{3}\rtimes\mathfrak{S}_{3}$-action. 
\end{prop}

\begin{proof}
There exist three ${\rm GL_{2}}$-actions on $\mathscr{H}_{\mathbb{A}}^{13}$
defined by the following rules (1)--(3):

\begin{enumerate}[(1)]

\item For an element $g_{1}\in{\rm GL_{2}}$, we define the ${\rm GL}_{2}$-action
on the affine space $\mA_{\mathscr{H}}$ by the rule 
\begin{align*}
 & \bm{x}_{1}\mapsto g_{1}\bm{x}_{1},u_{2}\mapsto(\det g_{1})u_{2},u_{3}\mapsto(\det g_{1})u_{3},\\
 & \left(\begin{array}{cccc}
p_{111} & p_{112} & p_{121} & p_{122}\\
p_{211} & p_{212} & p_{221} & p_{222}
\end{array}\right)\mapsto g_{1}\left(\begin{array}{cccc}
p_{111} & p_{112} & p_{121} & p_{122}\\
p_{211} & p_{212} & p_{221} & p_{222}
\end{array}\right),
\end{align*}
with the other coordinates being unchanged. We can check that $\bm{G}_{1},\dots,\mathsf{G}_{6}$
are changed as follows: 
\begin{align*}
 & \bm{G}_{1}\mapsto g_{1}\bm{G}_{1},\bm{G}_{2}\mapsto(\det g_{1})\bm{G}_{2},\bm{G}_{3}\mapsto(\det g_{1})\bm{G}_{3},\\
 & \mathsf{G}_{4}\mapsto(\det g_{1})^{2}\mathsf{G}_{4},\mathsf{G}_{5}\mapsto(\det g_{1})\mathsf{G}_{5},\mathsf{G}_{6}\mapsto(\det g_{1})\mathsf{G}_{6}.
\end{align*}
Therefore this ${\rm GL}_{2}$-action on the affine space $\mA_{\mathscr{H}}$
induces the action on $\mathscr{H}_{\mathbb{A}}^{13}$.

\item For an element $g_{2}\in{\rm GL_{2}}$, we set 
\begin{align*}
 & \bm{x}_{2}\mapsto g_{2}\bm{x}_{2},u_{1}\mapsto(\det g_{2})u_{1},u_{3}\mapsto(\det g_{2})u_{3},\\
 & \left(\begin{array}{cccc}
p_{111} & p_{112} & p_{211} & p_{212}\\
p_{121} & p_{122} & p_{221} & p_{222}
\end{array}\right)\mapsto g_{2}\left(\begin{array}{cccc}
p_{111} & p_{112} & p_{211} & p_{212}\\
p_{121} & p_{122} & p_{221} & p_{222}
\end{array}\right),
\end{align*}
 and other coordinates being unchanged. The changes of the equations
of $\mathscr{H}_{\mathbb{A}}^{13}$ are similar to the case (1), so
we omit them. We can check that this ${\rm GL}_{2}$-action on the
affine space $\mA_{\mathscr{H}}$ induces the action on $\mathscr{H}_{\mathbb{A}}^{13}$.

\item For and element $g_{3}\in{\rm GL_{2}}$, we set $\bm{x}_{3}\mapsto g_{3}\bm{x}_{3}$,
$u_{1}\mapsto(\det g_{3})u_{1}$, $u_{2}\mapsto(\det g_{3})u_{2}$,
{\small{} 
\[
\left(\begin{array}{cccc}
p_{111} & p_{121} & p_{211} & p_{221}\\
p_{112} & p_{122} & p_{212} & p_{222}
\end{array}\right)\mapsto g_{3}\left(\begin{array}{cccc}
p_{111} & p_{121} & p_{211} & p_{221}\\
p_{112} & p_{122} & p_{212} & p_{222}
\end{array}\right),
\]
}and other coordinates being unchanged. The changes of the equations
of $\mathscr{H}_{\mathbb{A}}^{13}$ are also similar to the case (1),
so we omit them. We can check that this ${\rm GL}_{2}$-action on
the affine space $\mA_{\mathscr{H}}$ induces the action on $\mathscr{H}_{\mathbb{A}}^{13}$.

\end{enumerate}

We can also check that these three ${\rm GL_{2}}$-actions on $\mathscr{H}_{\mathbb{A}}^{13}$
mutually commute. Thus these actions induce a $({\rm GL_{2})^{3}}$-action
on $\mathscr{H}_{\mathbb{A}}^{13}$.

We may define a natural $\mathfrak{S}_{3}$-action by associating
the three numbers $1,2,3$ to the three coordinates $u_{1},u_{2},u_{3}$,
the three vectors $\bm{x}_{1},\bm{x}_{2},\bm{x}_{3}$, and the three
subscripts $i,j,k$ for $p_{ijk}$. For example, by $(12)\in\mathfrak{S}_{3}$,
$u_{1}$ and $u_{2}$, $\bm{x}_{1}$ and $\bm{x}_{2}$, $p_{ijk}$
and $p_{jik}$ are interchanged respectively. Also, the action of
$(123)\in\mathfrak{S}_{3}$ cyclically permutes $u_{1},u_{2},u_{3}$,
$\bm{x}_{1},\bm{x}_{2},\bm{x}_{3}$, and the indices of $p_{ijk}$;
explicitly
\[
u_{1}\to u_{2}\to u_{3}\to u_{1},\,\bm{x}_{1}\to\bm{x}_{2}\to\bm{x}_{3}\to\bm{x}_{1},\,p_{ijk}\to p_{kij}\to p_{jki}\to p_{ijk}.
\]

We may check that the $({\rm GL}_{2})^{3}$-action and the $\mathfrak{S}_{3}$-action
are not commutative. Moreover, in ${\rm Aut}\,\mathscr{H}_{\mathbb{A}}^{13}$,
the image of $({\rm GL}_{2})^{3}$ is normalized by the image of $\mathfrak{S}_{3}$.
Therefore we have a $({\rm GL}_{2})^{3}\rtimes\mathfrak{S}_{3}$-action
on $\mathscr{H}_{\mathbb{A}}^{13}$. 
\end{proof}

\subsection{Cayley's hyperdeterminant\label{subsec:Cayley's-hyperdeterminant}}

Let $\mA_{\mathsf{P}}$ be the $\mathsf{k}$-affine $8$-space with
coordinates $p_{ijk}\,(i,j,k=1,2)$. Note that the $({\rm GL}_{2})^{3}\rtimes\mathfrak{S}_{3}$-action
on $\mA_{\mathscr{H}}$ induces the one on $\mA_{\mathsf{P}}$ by
the rules as in the proof of Proposition \ref{prop:actionH}. Now
we consider the hypermatrix $\mathsf{P}$ as in Remark \ref{rem:It-is-convenient}. 

The induced $({\rm SL}_{2})^{3}$-action on $\mathbb{A}_{\mathsf{P}}$
is studied in \cite[Chap.14, Ex.4.5]{GKZ}; $\mathbb{A}_{\mathsf{P}}$
has seven orbits. By this result and by considering the induced $\mathfrak{S}_{3}$-action
on $\mathbb{A}_{\mathsf{P}}$ together, we see that $\mathbb{A}_{\mathsf{P}}$
has five ${\rm (GL_{2})^{3}\rtimes\mathfrak{S}_{3}}$-orbits when
$\mathsf{k}$ is algebraically closed. To state the result, we set
\begin{align*}
\mathsf{p}_{1} & :=\text{the}\ p_{111}\text{-point},\\
\mathsf{p}_{2} & :=\text{the point with }p_{111}=p_{221}=1\ \text{and the other coordinates being zero},\\
\mathsf{p}_{3} & :=\text{the point with\ }p_{111}=p_{122}=p_{212}=1\ \text{and the other coordinates being zero,}\\
\mathsf{p}_{4} & :=\text{the point with\ }p_{111}=p_{222}=1\ \text{and the other coordinates being zero.}
\end{align*}
 Moreover, we set
\begin{align}
D_{\mathscr{H}}(\mathsf{P}) & :=p_{111}^{2}p_{222}^{2}+p_{112}^{2}p_{221}^{2}+p_{121}^{2}p_{212}^{2}+p_{122}^{2}p_{211}^{2}\label{eq:DH}\\
 & -2(p_{111}p_{122}p_{211}p_{222}+p_{111}p_{121}p_{212}p_{222}+p_{111}p_{112}p_{221}p_{222}\nonumber \\
 & +p_{121}p_{122}p_{211}p_{212}+p_{112}p_{122}p_{211}p_{221}+p_{112}p_{121}p_{212}p_{221})\nonumber \\
 & +4(p_{111}p_{122}p_{212}p_{221}+p_{112}p_{121}p_{211}p_{222}),\nonumber 
\end{align}
which is the hyperdeterminant of the hypermatrix $\mathsf{{P}}$,
classically known as the Cayley hyperdeterminant. It is known that
$\{D_{\mathscr{H}}=0\}$ is stable under the $({\rm GL}_{2})^{3}\rtimes\mathfrak{S}_{3}$-action
on $\mathbb{A}_{\mathsf{P}}$. 
\begin{prop}
\label{prop:orbits}Assume that $\mathsf{k}$ is algebraically closed.
Then the affine space $\mA_{\mathsf{P}}$ has the following five ${\rm (GL_{2})^{3}\rtimes\mathfrak{S}_{3}}$-orbits:
\end{prop}

\begin{enumerate}[(1)]

\item The origin of $\mathbb{A}_{\mathsf{P}}$.

\item The $4$-dimensional orbit of $\mathsf{p}_{1}$, which is the
complement $O_{1}$ of the origin in the affine cone over the Segre
variety $\mathbb{P}^{1}\times\mathbb{P}^{1}\times\mathbb{P}^{1}$.

\item The $5$-dimensional orbit of $\mathsf{p}_{2}$, which is the
complement $O_{2}$ of the affine cone over the Segre variety $\mathbb{P}^{1}\times\mathbb{P}^{1}\times\mathbb{P}^{1}$
in the union of three copies of the affine cone over the Segre variety
$\mathbb{P}^{1}\times\mathbb{P}^{3}$.

\item The $7$-dimensional orbit of $\mathsf{p}_{3}$, which is the
complement $O_{3}$ of the union of three copies of the affine cone
over the Segre variety $\mathbb{P}^{1}\times\mathbb{P}^{3}$ in $\{D_{\mathscr{H}}=0\}$.

\item The open orbit of $\mathsf{p}_{4}$, which is the complement
$O_{4}$ of $\{D_{\mathscr{H}}=0\}$ in $\mA_{\mathsf{P}}$.

\end{enumerate}
\begin{proof}
In \cite[Chap.14, Ex.4.5]{GKZ}, the result is stated without proof
for $\mathsf{k=\mC}$ (see \cite{Pa} for a proof). We give a quick
proof for convenience. 

We may identify $\mA_{\mathsf{P}}$ with the affine space associated
to the vector space $V_{1}\otimes V_{2}\otimes V_{3}$ with $V_{i}\simeq\mathsf{k}^{2}\,(1\leq i\leq3)$,
where ${\rm GL}_{2}$ acts on $V_{i}$ in a standard way. We may identify
$V_{2}\otimes V_{3}$ with the space of $2\times2$ $\mathsf{k}$-matrices
as ${\rm GL}_{2}\times{\rm GL}_{2}$-representation, where $(G_{1},G_{2})\in{\rm GL}_{2}\times{\rm GL}_{2}$
acts on a $2\times2$ $\mathsf{k}$-matrix $A$ as $A\mapsto G_{1}AG_{2}^{-1}$.
Having this, we write an element $\mathsf{p}$ of $V_{1}\otimes V_{2}\otimes V_{3}$
as $\bm{e}_{1}\otimes A_{1}+\bm{e}_{2}\otimes A_{2}$, where $\bm{e}_{1},\bm{e}_{2}$
form a basis of $V_{1}$, and $A_{1},A_{2}$ are $2\times2$ $\mathsf{k}$-matrices.
We also take bases $\bm{g}_{1},\bm{g}_{2}$ of $V_{2}$ and $\bm{h}_{1},\bm{h}_{2}$
of $V_{3}$. We determine the orbit of $\mathsf{p}$. We may assume
that $\bm{e}_{1}\otimes A_{1}+\bm{e}_{2}\otimes A_{2}\not=0$, i.e.,
$A_{1}\not=O$ or $A_{2}\not=O$. By the ${\rm GL}_{2}$-action on
$V_{1}$, we may assume that $A_{1}\not=O$.

Assume that $\rank A_{1}=2$. Then, by the ${\rm GL}_{2}\times{\rm GL}_{2}$-action,
we may assume that $A_{1}$ is the identity matrix $I$. By a further
transformation of the form $(I,A_{2})\mapsto(I,GA_{2}G^{-1})$, we
may assume that $A_{2}=\left(\begin{array}{cc}
\alpha & 0\\
0 & \beta
\end{array}\right)$ or $\left(\begin{array}{cc}
\alpha & 1\\
0 & \alpha
\end{array}\right)\,(\alpha,\beta\in\mathsf{k})$ (\textit{here we use the assumption $\mathsf{k}$ is algebraically
closed}). 

In the former case, the element $\mathsf{p}=\bm{e}_{1}\otimes I+\bm{e}_{2}\otimes A_{2}$
is identified with $(\bm{e}_{1}+\alpha\bm{e}_{2})\otimes\bm{g}_{1}\otimes\bm{h}_{1}+(\bm{e}_{1}+\beta\bm{e}_{2})\otimes\bm{g}_{2}\otimes\bm{h}_{2}$.
If $\alpha\not=\beta$, then, by the ${\rm GL}_{2}$-action on $V_{1}$,
we may assume that $\mathsf{p}=\bm{e}_{1}\otimes\bm{g}_{1}\otimes\bm{h}_{1}+\bm{e}_{2}\otimes\bm{g}_{2}\otimes\bm{h}_{2}$,
which is nothing but the point $\mathsf{p}_{4}$. If $\alpha=\beta$,
then, by the ${\rm GL}_{2}$-action on $V_{1}$, we may assume that
$\mathsf{p}=\bm{e}_{1}\otimes(\bm{g}_{1}\otimes\bm{h}_{1}+\bm{g}_{2}\otimes\bm{h}_{2})$,
which belongs to the orbit of the point $\mathsf{p}_{2}$. 

In the latter case, the element $\mathsf{p}=\bm{e}_{1}\otimes I+\bm{e}_{2}\otimes A_{2}$
is identified with $(\bm{e}_{1}+\alpha\bm{e}_{2})\otimes(\bm{g}_{1}\otimes\bm{h}_{1}+\bm{g}_{2}\otimes\bm{h}_{2})+\bm{e}_{2}\otimes\bm{g}_{1}\otimes\bm{h}_{2}$.
By the ${\rm GL}_{2}$-action on $V_{1}$, we may assume that $\mathsf{p}=\bm{e}_{1}\otimes(\bm{g}_{1}\otimes\bm{h}_{1}+\bm{g}_{2}\otimes\bm{h}_{2})+\bm{e}_{2}\otimes\bm{g}_{1}\otimes\bm{h}_{2}$,
which is nothing but the point $\mathsf{p}_{3}$. 

Assume that $\rank A_{1}=$1. If $\rank A_{2}=2$, then the above
discussion applies after interchanging $A_{1}$ and $A_{2}$ via the
${\rm GL}_{2}$-action on $V_{1}$ and then we see that $\mathsf{p}$
is in the orbit of $\mathsf{p}_{2}$. Therefore we may assume that
$A_{2}=O$ or $\rank A_{2}=1$. If $A_{2}=O$, then, by the ${\rm GL}_{2}\times{\rm GL}_{2}$-action
on $V_{2}\otimes V_{3}$, we may assume that $\mathsf{p}=\bm{e}_{1}\otimes\bm{g}_{1}\otimes\bm{h}_{1}$,
which is nothing but the point $\mathsf{p}_{1}$. Assume that $\rank A_{2}=1$.
Then by the ${\rm GL}_{2}\times{\rm GL}_{2}$-action on $V_{2}\otimes V_{3}$,
we may assume that $\mathsf{p}=\bm{e}_{1}\otimes\bm{g}_{1}\otimes\bm{h}_{1}+\bm{e}_{2}\otimes\bm{v}_{2}\otimes\bm{v}_{3},$
where $\bm{v}_{i}$ is a nonzero vector of $V_{i}\,(i=2,3)$. Now
it is easy to see the following: If $\bm{v}_{2}$ is parallel to $\bm{g}_{1}$
and $\bm{v}_{3}$ is parallel to $\bm{h}_{1}$, then $\mathsf{p}$
belongs to the orbit of $\mathsf{p}_{1}$. If $\bm{v}_{2}$ is parallel
to $\bm{g}_{1}$ and $\bm{v}_{3}$ is not parallel to $\bm{h}_{1}$,
then $\mathsf{p}$ belongs to the orbit of $\mathsf{p}_{2}$. If $\bm{v}_{2}$
is not parallel to $\bm{g}_{1}$ and $\bm{v}_{3}$ is not parallel
to $\bm{h}_{1}$, then $\mathsf{p}$ belongs to the orbit of $\mathsf{p}_{4}$. 

Now we have shown that all the orbits are those of $o$ and $\mathsf{p}_{i}\,(1\leq i\leq4)$.
It is easy to check that $\mathsf{p}_{i}\in O_{i}\,(1\leq i\leq4),$
and $O_{i}$ are $({\rm GL}_{2})^{3}\rtimes\mathfrak{S}_{3}$-stable
(note that the three copies of the Segre variety $\mathbb{P}^{1}\times\mathbb{P}^{3}$
as in the statement of (3) are nothing but $\mP(V_{1})\times\mP(V_{2}\otimes V_{3})$,
$\mP(V_{1}\otimes V_{2})\times\mP(V_{3})$, $\mP(V_{2})\times\mP(V_{1}\otimes V_{3})$).
Since $O_{i}$ are mutually disjoint, we conclude that $O_{i}\,(1\leq i\leq4)$
coincides with the orbit of $\mathsf{p}_{i}$. 
\end{proof}

\subsection{$\mathbb{P}^{2}\times\mathbb{P}^{2}$-fibration associated to $\mathscr{H}_{\mathbb{A}}^{13}$\label{subsec:-fibration-associated-toP2P2}}

Note that any equation of $\text{\ensuremath{\mathscr{H}_{\mathbb{A}}^{13}}}$
is of degree $2$ if we regard the coordinates of $\mathbb{A_{\mathsf{P}}}$
as constants. Therefore, considering the coordinates $u_{j},x_{ij}\,(i=1,2,j=1,2,3)$
as projective coordinates, we obtain a quasi-projective variety with
the same equations as $\text{\ensuremath{\mathscr{H}_{\mathbb{A}}^{13}}}$,
which we denote by $\widehat{\mathscr{H}}$. We also denote by $\rho_{\mathscr{H}}\colon\widehat{\mathscr{H}}\to\mathbb{A_{\mathsf{P}}}$
the natural projection.
\begin{defn}
We set
\[
\mathsf{S}:=\left(\begin{array}{ccc}
s_{11} & s_{12} & s_{13}\\
s_{12} & s_{22} & s_{23}\\
s_{13} & s_{23} & s_{33}
\end{array}\right),\,\bm{{\sigma}}:=\left(\begin{array}{c}
\sigma_{1}\\
\sigma_{2}\\
\sigma_{3}
\end{array}\right),
\]
 and denote by $\mathsf{S}^{\dagger}$ the adjoint matrix of $\mathsf{S}$.
We define 
\[
\mathbb{P}^{2,2}:=\{(\mathsf{S},\bm{\sigma})\mid\mathsf{S}^{\dagger}=O,\mathsf{S}\bm{\sigma}=\bm{0}\}\subset\mathbb{P}^{8},
\]
which is introduced in \cite{Mu3} as a degeneration of $\mathbb{P}^{2}\times\mathbb{P}^{2}$
and is shown to be a singular del Pezzo 4-fold of degree six (see
also \cite{Fu}). 
\end{defn}

\begin{prop}
\label{prop:Let--beP2P2fib} Assume that $\mathsf{k}$ is algebraically
closed. We use the notation of Proposition \ref{prop:orbits}. The
$\rho_{\mathscr{H}}$-fibers over points in the ${\rm (GL_{2})^{3}\rtimes\mathfrak{S}_{3}}$
-orbit of $\mathsf{p}_{4}$ are isomorphic to $\mathbb{P}^{2}\times\mathbb{P}^{2}$,
and the $\rho_{\mathscr{H}}$-fibers over points in the ${\rm (GL_{2})^{3}\rtimes\mathfrak{S}_{3}}$
-orbit of $\mathsf{p}_{3}$ are isomorphic to $\mathbb{P}^{2,2}$.
Any $\rho_{\mathscr{H}}$-fibers except the fiber over the origin
are (not necessarily irreducible) $4$-folds of degree $6$ while
the fiber over the origin is $5$-dimensional. In particular, $\mathscr{H}_{\mA}^{13}$
is $13$-dimensional.
\end{prop}

\begin{proof}
The assertions immediately follows by giving descriptions of all the
$\rho_{\mathscr{H}}$-fibers. Note that the ${\rm (GL_{2})^{3}\rtimes\mathfrak{S}_{3}}$-action
on $\text{\text{\ensuremath{\mathscr{H}_{\mathbb{A}}^{13}}}}$ induces
that on $\widehat{\mathscr{H}}$. Then $\rho_{\mathscr{H}}$-fibers
are isomorphic to the fibers over the origin or the $\mathsf{p}_{i}$-points
for some $i=1,2,3,4$, which are easily described by the equations
of $\mathscr{H}_{\mathbb{A}}^{13}$ as follows:

\begin{enumerate}[$(1)$]

\item The fiber over the origin is the union o{\small f $\{u_{1}=u_{2}=u_{3}=0\}\simeq\mP^{5},$
and three $\mP^{4}$'s} $\{u_{1}=u_{2}=x_{13}=x_{23}=0\}$, $\{u_{1}=u_{3}=x_{12}=x_{22}=0\}$,
and $\{u_{2}=u_{3}=x_{11}=x_{21}=0\}.$

\item The fiber over the point $\mathsf{p}_{1}$ is the union of
three quadric $4$-folds of rank $4$; {\small
\begin{align*}
 & \{u_{1}=u_{2}=x_{23}=u_{3}x_{13}-x_{21}x_{22}=0\}\cup\{u_{1}=u_{3}=x_{22}=u_{2}x_{12}-x_{21}x_{23}=0\}\\
 & \cup\{u_{2}=u_{3}=x_{21}=u_{1}x_{11}-x_{22}x_{23}=0\}.
\end{align*}
}\item The fiber over the point $\mathsf{p}_{2}$ is the union of
a quadric $4$-fold of rank $6$ and the $4$-dimensional cone over
a hyperplane section of $\mathbb{P}^{1}\times\mathbb{P}^{3}$; {\small
\begin{align*}
 & \{u_{1}=u_{2}=x_{23}=u_{3}x_{13}-x_{11}x_{12}-x_{21}x_{22}=0\}\\
 & \cup\left\{ u_{3}=0,{\rm rank\,}\begin{pmatrix}x_{23} & u_{2} & x_{11} & x_{21}\\
u_{1} & -x_{23} & x_{22} & -x_{12}
\end{pmatrix}\leq1\right\} .
\end{align*}
}\item The fiber over the point $\mathsf{p}_{3}$: {\small
\[
\left\{ {\rm rank}\,\begin{pmatrix}u_{1} & x_{13} & x_{22}\\
x_{13} & u_{2} & -x_{21}\\
x_{22} & -x_{21} & -u_{3}
\end{pmatrix}\leq1,\begin{pmatrix}u_{1} & x_{13} & x_{22}\\
x_{13} & u_{2} & -x_{21}\\
x_{22} & -x_{21} & -u_{3}
\end{pmatrix}\begin{pmatrix}-x_{11}\\
x_{12}\\
x_{23}
\end{pmatrix}=\bm{o}\right\} ,
\]
}which is isomorphic to $\mathbb{P}^{2,2}$.

\item The fiber over the point $\mathsf{p}_{4}$: {\small
\[
\left\{ {\rm rank}\,\begin{pmatrix}u_{1} & x_{13} & x_{22}\\
x_{23} & u_{2} & x_{11}\\
x_{12} & x_{21} & u_{3}
\end{pmatrix}\leq1\right\} ,
\]
}which is isomorphic to $\mathbb{P}^{2}\times\mathbb{P}^{2}$.

\end{enumerate}

Therefore we have shown the assertions. 
\end{proof}
\begin{rem}
\begin{enumerate}[(1)]

\item By Proposition \ref{prop:UFD} below, the relative Picard number
of $\rho_{\mathscr{H}}$ is $1$. 

\item By Proposition \ref{prop:Let--beP2P2fib}, a general $\rho_{\mathscr{H}}$-fiber
is isomorphic to $\mP^{2}\times\mP^{2}$. Usually, in this situation,
$\rho_{\mathscr{H}}$ is said to be a generically $\mP^{2}\times\mP^{2}$-fibration,
but for simplicity we shall simply say that $\rho_{\mathscr{H}}$
is a $\mP^{2}\times\mP^{2}$-fibration.

\end{enumerate}
\end{rem}

\subsection{Nondegeneracy of $J_{\mathsf{P}}$\label{subsec:Nondegeneracy-of}}

Assume that $\mathsf{k}$ is algebraically closed. Now we revisit
the quadratic Jordan algebra $J_{\mathsf{P}}$ of the cubic form $N_{\mathsf{P}}$
for each value of $\mathsf{P\in\mA_{\mathsf{P}}}$ as in Proposition
\ref{Prop:reveseJ}; we study nondegeneracy of $J_{\mathsf{P}}$.
We may regard $\mA_{\mathscr{H}}\to\mA_{\mathsf{P}}$ as a family
of $J_{\mathsf{P}}.$ By the group action of $\mathscr{H}_{\mA}^{13}$
as in Subsection \ref{subsec:-action-onH}, we have only to consider
nondegeneracy over the origin, the $\mathsf{p}_{1}$-, $\mathsf{p}_{2}$-,
$\mathsf{p_{3}}$-, and $\mathsf{p}_{4}$-point. For these points,
we calculate the subvarieties $R_{J_{\mathsf{P}}}$ defined as (\ref{eq:nondeg}),
which measure how degenerate the quadratic Jordan algebras $J_{\mathsf{P}}$
are. The results as presented below follow from direct calculations:

\vspace{3pt}

\noindent \textbf{The origin:}
\begin{align*}
R_{J_{\mathsf{P}}}=\{u_{1}=u_{2}=u_{3}=0\},\ N_{\mathsf{P}}=u_{1}u_{2}u_{3}.
\end{align*}

\vspace{3pt}

\noindent \noindent \textbf{The $\mathsf{p}_{1}$-point : }
\begin{align*}
R_{J_{\mathsf{P}}} & =\{u_{1}=u_{2}=u_{3}=x_{22}x_{23}=x_{21}x_{23}=x_{21}x_{22}=0\},\\
N_{\mathsf{P}} & =u_{1}u_{2}u_{3}.
\end{align*}

\vspace{3pt}

\noindent \textbf{The $\mathsf{p}_{2}$-point: }
\begin{align*}
R_{J_{\mathsf{P}}} & =\{u_{1}=u_{2}=u_{3}=x_{23}=x_{11}x_{12}+x_{21}x_{22}=0\},\\
N_{\mathsf{P}} & =u_{3}(u_{1}u_{2}+x_{23}^{2}).
\end{align*}

\vspace{3pt}

\noindent \textbf{The $\mathsf{p}_{3}$-point: }
\begin{align*}
R_{J_{\mathsf{P}}} & =\{u_{1}=u_{2}=u_{3}=x_{21}=x_{22}=x_{13}=0\},\\
N_{\mathsf{P}} & =u_{1}u_{2}u_{3}-u_{3}x_{13}^{2}+u_{1}x_{21}^{2}+2x_{13}x_{21}x_{22}+u_{2}x_{22}^{2}.
\end{align*}

\vspace{3pt}

\noindent \textbf{The $\mathsf{p}_{4}$-point:} $R_{J_{\mathsf{P}}}=\{0\}$,
i.e., $R_{J_{\mathsf{P}}}$ is nondegenerate.

\vspace{3pt}


Consequently, we have the following:
\begin{cor}
$J_{\mathsf{P}}$ is nondegenerate if and only if $\mathsf{P}$ belongs
to the open ${\rm (GL_{2})^{3}\rtimes\mathfrak{S}_{3}}$-orbit, namely,
$D_{\mathscr{H}}(\mathsf{P})\not=0$.
\end{cor}

\vspace{5pt}

In the remaining part of Section \ref{sec:Affine-varietyH}, we obtain
several more properties of $\mathscr{H}_{\mA}^{13}$ which are needed
to construct prime $\mQ$-Fano 3-folds from weighted projectivizations
of $\mathscr{H}_{\mA}^{13}$.

\subsection{Charts of $\mathfrak{\mathscr{H}}_{\mA}^{13}$\label{subsec:Descriptrion-of-charts}}

For a coordinate $*$, we call the open subset of $\mathfrak{\mathscr{H}}_{\mA}^{13}$
with $*\not=0$ \textit{the $*$-chart}. We describe the $*$-chart
such that $*$ is an entry of $\bm{x}_{i}$, or $u_{i}\,(i=1,2,3)$.
By the $\mathfrak{S}_{3}$-action on $\mathscr{H}_{\mA}^{13}$ as
in Proposition \ref{prop:actionH}, we may swap $u_{1},u_{2},u_{3}$,
and $\bm{x}_{1},\bm{x}_{2},\bm{x}_{3}$. We can also check that, by
the action of {\small
\[
\left(\left(\begin{array}{cc}
0 & 1\\
1 & 0
\end{array}\right),\left(\begin{array}{cc}
0 & 1\\
1 & 0
\end{array}\right),\left(\begin{array}{cc}
0 & 1\\
1 & 0
\end{array}\right)\right)\in({\rm GL}_{2})^{3},
\]
}we may swap $x_{1i}$ and $x_{2i}$ for $i=1,2,3.$ Therefore it
suffices to treat only the $u_{1}$-, $x_{11}$-charts.

\vspace{3pt}

\noindent\underline{\bf $u_1$-chart}: We note that $\mathfrak{\mathscr{H}}_{\mA}^{13}\cap\{u_{1}\not=0\}$
is isomorphic to $(\mathfrak{\mathscr{H}}_{\mA}^{13}\cap\{u_{1}=1\})\times\mA^{1*}$
by the map $(\bm{x}_{1},\bm{x}_{2},\bm{x}_{3},u_{1},u_{2},u_{3},p_{ijk})\mapsto\left((u_{1}^{-1}\bm{x}_{1},u_{1}^{-1}\bm{x}_{2},u_{1}^{-1}\bm{x}_{3},u_{1}^{-1}u_{2},u_{1}^{-1}u_{3},p_{ijk}),u_{1}\right)$.
This is because all the equations of $\mathfrak{\mathscr{H}}_{\mA}^{13}$
are quadratic when we consider $p_{ijk}$ are constants. Therefore
it suffices to describe $\mathscr{H}_{\mA}^{13}\cap\{u_{1}=1\}$.
Solving regularly the $9$ equations of $\mathfrak{\mathscr{H}}_{\mA}^{13}$
setting $u_{1}=1$, we see that the $9$ equations are reduced to
the following 4 equations: 
\begin{align*}
\bm{x}_{1}=D^{(3)}\bm{x}_{2},\ u_{3}=-\det D^{(2)},\ u_{2}=-\det D^{(1)}.
\end{align*}

This description of the $u_{1}$-chart shows that the $u_{1}$-, $u_{2}$-,
and $u_{3}$-charts of $\mathfrak{\mathscr{H}}_{\mA}^{13}$ are isomorphic
to $\mA^{12}\times\mA^{1*}$. 

\vspace{3pt}

\noindent\underline{\bf $x_{11}$-chart}: Similarly to the $u_{1}$-chart,
it suffices to describe $\mathfrak{\mathscr{H}}_{\mA}^{14}\cap\{x_{11}=1\}$.
When $x_{11}=1$, we may eliminate the coordinate $u_{1}$ by $u_{1}\bm{x}_{1}-D^{(3)}\bm{x}_{2}=0$
with $x_{11}=1$. Then we may verify that $\mathfrak{\mathscr{H}}_{\mA}^{14}\cap\{x_{11}=1\}$
is defined by the five $4\times4$ Pfaffians of the following skew-symmetric
matrix: 

\noindent
\begin{equation}
M_{{\rm T}}:=\left(\begin{array}{ccccc}
0 & x_{13} & x_{23} & -x_{22} & x_{12}\\
 & 0 & -u_{2} & D_{21}^{(1)} & -D_{11}^{(1)}\\
 &  & 0 & D_{22}^{(1)} & -D_{12}^{(1)}\\
 &  &  & 0 & u_{3}\\
 &  &  &  & 0
\end{array}\right)\label{eq:MT}
\end{equation}
with $x_{11}=1$. 

\vspace{10pt}

Based on the above results, we describe the singularities of $\mathscr{H}_{\mA}^{13}\setminus\mA_{\mathsf{P}}$.

We set 

\begin{align*}
\Delta_{1} & :=\{u_{1}=u_{2}=u_{3}=0,\bm{x}_{2}=\bm{x}_{3}=\bm{o},D^{(1)}=O,\bm{x}_{1}\not=\bm{o}\},\\
\Delta_{2} & :=\{u_{1}=u_{2}=u_{3}=0,\bm{x}_{1}=\bm{x}_{3}=\bm{o},D^{(2)}=O,\bm{x}_{2}\not=\bm{o}\},\\
\Delta_{3} & :=\{u_{1}=u_{2}=u_{3}=0,\bm{x}_{1}=\bm{x}_{2}=\bm{o},D^{(3)}=O,\bm{x}_{3}\not=\bm{o}\},
\end{align*}
and 
\begin{equation}
\Delta:=\Delta_{1}\cup\Delta_{2}\cup\Delta_{3}.\label{eq:Delta}
\end{equation}
We have $\Delta\cap\{x_{ij}\not=0\}=\Delta_{j}\cap\{x_{ij}\not=0\}\,(i=1,2,j=1,2,3).$
Identifying $\mA_{\mathsf{P}}$ with the subspace $\{\bm{x}_{1}=\bm{x}_{2}=\bm{x}_{3}=\bm{o},u_{1}=u_{2}=u_{3}=0\}$
of $\mA_{\mathscr{H}}$, we consider $\mA_{\mathsf{P}}\subset\mathscr{H}_{\mA}^{13}$.

By the above descriptions of the charts and the ${\rm (GL_{2})^{3}\rtimes\mathfrak{S}_{3}}$-action
on $\mathscr{H}_{\mA}^{13}$, we have the following:
\begin{prop}
\label{prop:The-singular-locus} The open subset $\mathscr{H}_{\mA}^{13}\setminus\mA_{\mathsf{P}}$
of $\mathscr{H}_{\mA}^{13}$ is $13$-dimensional and irreducible.
Its singular locus is equal to $\Delta$, and it has $c({\rm G}(2,5))$-singularities
along $\Delta$, where we call a singularity isomorphic to the vertex
of the cone over ${\rm G}(2,5)$ a \textup{$c({\rm G}(2,5))$-singularity. }
\end{prop}

We will show that $\mathscr{H}_{\mA}^{13}$ itself is 13-dimensional
and irreducible (Proposition \ref{prop:UFD}), and describe the singularities
of $\mathfrak{\mathscr{H}}_{\mA}^{13}$ along $\mA_{\mathsf{P}}$
in Subsection \ref{subsec:Singularties-ofH} when $\mathsf{k=\mC}$.

\vspace{1cm}

\textit{Hereafter till the end of the paper, we assume for simplicity
that $\mathsf{k}=\mC$.}

\subsection{Gorensteinness and $9\times16$ graded minimal free resolution of
the ideal of $\mathscr{H}_{\mA}^{13}$\label{subsec:Gorensteinness-and-9 16 H}}
\begin{prop}
\label{prop:9times16 H} Let $S_{\mathscr{H}}$ be the polynomial
ring over $\mC$ whose variables are $u_{1}$,$u_{2}$,$u_{3}$, $p_{ijk}\,(i,j,k=1,2)$
and the entries of $\bm{x}_{1},\bm{x}_{2},\bm{x}_{3}$. Let $I_{\mathscr{H}}$
be the ideal of the polynomial ring $S_{\mathscr{H}}$ generated by
the $9$ equations (\ref{eq:HypMat}) of $\mathscr{H}_{\mathbb{A}}^{13}$.
\begin{enumerate}[$(1)$]

\item We give nonnegative weights for coordinates of $S_{\mathscr{H}}$
such that all the equations of $\mathscr{H}_{\mA}^{13}$ are weighted
homogeneous, and we denote by $w(*)$ the weight of the monomial $*$.
We denote by $\mP$ the corresponding weighted projective space, and
by $\mathscr{H}_{\mP}\subset\mP$ the weighted projectivization of
$\mathscr{H}_{\mA}^{13}$, where we allow some coordinates being nonzero
constants (thus $\dim\mathscr{H}_{\mP}$ could be less than $12$).
We set
\[
c:=w(x_{11}x_{21}u_{1}),\,d:=w(u_{1}u_{2}u_{3}),\,\delta:=c+d.
\]

\begin{enumerate}[$({1}\text{-}1)$]

\item It holds that $\omega_{\mP}=\sO_{\mP}(c-4d)$.

\item The ideal $I_{\mathscr{H}}$ has the following graded minimal
$S_{\mathscr{H}}$-free resolution
\begin{equation}
0\leftarrow P_{0}\leftarrow P_{1}\leftarrow P_{2}\leftarrow P_{3}\leftarrow P_{4}\leftarrow0,\ \text{where}\label{eq:minfree}
\end{equation}
 
\begin{align*}
P_{0}= & S:=S_{\mathscr{H}},\\
P_{1}= & S(-(c-w(x_{11})))\oplus S(-(c-w(x_{21})))\oplus S(-(c-w(x_{12})))\oplus S(-(c-w(x_{22})))\\
\oplus & S(-(c-w(x_{13})))\oplus S(-(c-w(x_{23})))\\
\oplus & S(-(d-w(u_{1})))\oplus S(-(d-w(u_{2})))\oplus S(-(d-w(u_{3}))),\\
P_{2}= & S(-c)^{\oplus2}\oplus S(-d)^{\oplus2}\\
\oplus & S(-w(u_{1}u_{2}x_{11}))\oplus S(-w(u_{1}u_{2}x_{21}))\oplus S(-w(u_{1}u_{3}x_{11}))\oplus S(-w(u_{1}u_{3}x_{21}))\\
\oplus & S(-w(u_{1}u_{2}x_{12}))\oplus S(-w(u_{1}u_{2}x_{22}))\oplus S(-w(u_{2}u_{3}x_{12}))\oplus S(-w(u_{2}u_{3}x_{22}))\\
\oplus & S(-w(u_{1}u_{3}x_{13}))\oplus S(-w(u_{1}u_{3}x_{23}))\oplus S(-w(u_{2}u_{3}x_{13}))\oplus S(-w(u_{2}u_{3}x_{23})),\\
P_{3}= & S(-(d+w(x_{11})))\oplus S(-(d+w(x_{21})))\oplus S(-(d+w(x_{12})))\oplus S(-(d+w(x_{22})))\\
\oplus & S(-(d+w(x_{13})))\oplus S(-(d+w(x_{23})))\\
\oplus & S(-(c+w(u_{1})))\oplus S(-(c+w(u_{2})))\oplus S(-(c+w(u_{3}))),\\
P_{4}= & S(-\delta).
\end{align*}

\item It holds that 
\begin{equation}
\omega_{\mathscr{H}_{\mathbb{P}}}=\sO_{\mathscr{H}_{\mP}}(c-4d+\delta)=\sO_{\mathscr{H}_{\mP}}(2c-3d).\label{eq:omegaP}
\end{equation}

\end{enumerate}

\item $I_{\mathscr{H}}$ is a Gorenstein ideal of codimension $4$. 

\item $\mathscr{H}_{\mathbb{A}}^{13}$ is irreducible and reduced,
thus $I_{\mathscr{H}}$ is a prime ideal.

\item $\mathscr{H}_{\mathbb{A}}^{13}$ is normal. 

\end{enumerate}
\end{prop}

\begin{proof}
(1). We may compute the $S_{\mathscr{H}}$-free resolution (\ref{eq:minfree})
of $R_{\mathfrak{\mathscr{H}}}$ by \textsc{Singular} \cite{DGPS}.
Since $\mathscr{H}_{\mathbb{A}}^{13}$ is codimension 4 in $\mA_{\mathscr{H}}$
by Proposition \ref{prop:Let--beP2P2fib}, the length of an $S_{\mathscr{H}}$-free
resolution of $I_{\mathscr{H}}$ is at least 4. Therefore the minimality
of (\ref{eq:minfree}) follows. We have (1-1) by the conditions of
weights of coordinates, and (1-3) follows from (1-1) and (1-2). 

\noindent (2). By (1), $I_{\mathscr{H}}$ is Cohen-Macaulay. Moreover,
since the last term of the minimal free resolution is $S_{\mathscr{H}}(-\delta)$,
which is of rank 1, $I_{\mathscr{H}}$ is Gorenstein. 

\noindent (3). We use the notation of Proposition \ref{prop:orbits}.
Let $V$ be the union of the ${\rm (GL_{2})^{3}\rtimes\mathfrak{S}_{3}}$-orbits
of the points $\mathsf{p}_{3}$ and $\mathsf{p}_{4}$. Note that $V$
is an open subset of $\mA_{\mathsf{P}}$ such that $\dim(\mA_{\mathsf{P}}\setminus V)=6.$
By Proposition \ref{prop:Let--beP2P2fib}, the $\rho_{\mathscr{H}}$-inverse
image $U$ of $V$ is irreducible since any $\rho_{\mathscr{H}}$-fiber
in $V$ is isomorphic to $\mP^{2}\times\mP^{2}$ or $\mP^{2,2}$.
Moreover, by ibid., we can verify that $\widehat{\mathscr{H}}\setminus U$
is $10$-dimensional. Therefore, $\mathscr{H}_{\mA}^{13}$ has $13$-dimensional
irreducible open subset, and its complement is $11$-dimensional.
Since $\mathscr{H}_{\mA}^{13}$ is Gorenstein by (2), this implies
that $\mathscr{H}_{\mA}^{13}$ is irreducible and reduced by unmixedness
of a Gorenstein variety. 

\noindent (4). By Proposition \ref{prop:The-singular-locus}, $\mathscr{H}_{\mathbb{A}}^{13}$
is regular in codimension 1. By (2), $\mathscr{H}_{\mathbb{A}}^{13}$
satisfies $S_{2}$ condition. Therefore, $\mathscr{H}_{\mathbb{A}}^{13}$
is normal.
\end{proof}

\subsection{Factoriality of $\mathscr{H}_{\mathbb{A}}^{13}$\label{subsec:Factoriality}\protect 
}\begin{prop}
\label{prop:UFD} Let $S_{\mathscr{H}},\,I_{\mathscr{H}}$ be as in
Proposition \ref{prop:9times16 H}. The ring $R_{\mathscr{H}}:=S_{\mathscr{H}}/I_{\mathscr{H}}$
is a UFD. 
\end{prop}

\begin{proof}
\noindent Note that $R_{\mathscr{H}}$ is a domain by Proposition
\ref{prop:9times16 H} (3). Therefore, by Nagata's theorem \cite[Thm.~20.2]{Matsumura},
it suffices to show that the ring $(R_{\mathscr{H}})_{u_{1}}$ is
a UFD and $u_{1}$ is a prime element of $R_{\mathscr{H}}$.

By the description of the $u_{1}$-chart as in Subsection \ref{subsec:Descriptrion-of-charts},
we see that $(R_{\mathscr{H}})_{u_{1}}$ is a localization of a polynomial
ring, thus is a UFD.

We can show that $\mathscr{H}_{\mA}^{13}\cap\{u_{1}=0\}$ is irreducible
by a similar argument as in the proof of Proposition \ref{prop:9times16 H}
(3) using the fibration $\rho_{\mathscr{H}}$; note that $\mathscr{H}_{\mA}^{13}\cap\{u_{1}=0\}$
is also Gorenstein, and any fiber of $\rho_{\mathscr{H}}|_{\{u_{1}=0\}}$
contained in the open subset $U$ is irreducible since it is a hyperplane
section of $\mP^{2}\times\mP^{2}$ or $\mP^{2,2}$. Moreover, $(\widehat{\mathscr{H}}\setminus U)\cap\{u_{1}=0\}$
is at most $10$-dimensional. Therefore $\mathscr{H}_{\mA}^{13}\cap\{u_{1}=0\}$
is irreducible by a similar reason to the irreducibility of \textbf{$\mathscr{H}_{\mA}^{13}$},
hence $u_{1}$ is a prime element of $R_{\mathscr{H}}$. 
\end{proof}
\begin{cor}
\label{HFactorial} Let $\mathfrak{\mathscr{H}}_{\mP}^{12}$ be the
weighted projectivization of $\mathfrak{\mathscr{H}}_{\mA}^{13}$
with some positive weights of coordinates. The following assertions
hold:

\begin{enumerate}

\item Any prime Weil divisor on $\mathfrak{\mathscr{H}}_{\mP}^{12}$
is the intersection between $\mathfrak{\mathscr{H}}_{\mP}^{12}$ and
a weighted hypersurface. In particular, $\mathfrak{\mathscr{H}}_{\mP}^{12}$
is $\mQ$-factorial and has Picard number $1$. 

\item Let $X$ be a quasi-smooth threefold such that $X$ is a codimension
$9$ weighted complete intersection in $\mathfrak{\mathscr{H}}_{\mP}^{12}$,
i.e., there exist $9$ weighted homogeneous polynomials $G_{1},\dots,G_{9}$
such that $X=\mathfrak{\mathscr{H}}_{\mP}^{12}\cap\{G_{1}=0\}\cap\dots\cap\{G_{9}=0\}$.
Assume moreover that $\{u_{1}=0\}\cap X$ is a prime divisor. Then
$X$ is irreducible and any prime Weil divisor on $X$ is the intersection
between $X$ and a weighted hypersurface. In particular, $X$ is $\mQ$-factorial
and has Picard number $1$.\end{enumerate}
\end{cor}

\begin{proof}
(1) Let $\frak{\mathfrak{p}}\subset R_{\mathscr{H}}$ be the homogeneous
ideal defining a prime divisor $D$. Then $\mathfrak{p}$ is an ideal
of height $1$. Since $R_{\mathscr{H}}$ is a UFD by Proposition \ref{prop:UFD},
its ideal of height $1$ is principal by \cite[Thm.~20.1]{Matsumura}.
Thus the assertion follows. 

(2) By the assumption, $X$ is normal. Moreover, $X$ is connected
since $X$ is a weighted complete intersection in a normal irreducible
variety $\mathscr{H}_{\mP}^{12}$ by ample divisors. Therefore $X$
is irreducible. 

For the remaining assertion, once we check that $R_{X}:=R_{\mathscr{H}}/(F_{1},\dots,F_{9})$
is a UFD, the proof of (1) works verbatim. Let $R_{X,o}$ be the localization
of $R_{X}$ by the maximal irrelevant ideal. Note that, by \cite[Prop.~7.4]{UFD},
$R_{X}$ is a UFD if and only if so is $R_{X,o}$. Thus we have only
to check the latter. To check $R_{X,o}$ is a UFD, we apply Nagata's
theorem as in the proof of Proposition \ref{prop:UFD}. For this,
we need to check that $u_{1}$ is a prime element of $R_{X,o}$, which
follows from our assumption, and $(R_{X,o})_{u_{1}}$ is a UFD, which
we will check below. We denote by $R_{\mathscr{H},o}$ the localization
of $R_{\mathscr{H}}$ by the maximal irrelevant ideal. Since $\left(R_{\mathscr{H},o}\right)_{u_{1}}$
is a localization of a polynomial ring by the description of the $u_{1}$-chart
as in Subsection \ref{subsec:Descriptrion-of-charts}, we see that
$(R_{X,o})_{u_{1}}$ is a localization of a complete intersection
ring. Therefore $(R_{X,o})_{u_{1}}$ is a complete intersection local
ring of dimension 4. Since $X$ is quasi-smooth, we see that $(R_{X,o})_{u_{1}}$
is a UFD by \cite[Exp.XI, Cor.~3.10 and Thm.~3.13 (ii)]{sga2}. 
\end{proof}
\vspace{0.5cm}

\subsection{Singularities of $\mathscr{H}_{\mA}^{13}$\label{subsec:Singularties-ofH}}
\begin{prop}
\label{prop:HTerm} The variety $\mathscr{H}_{\mA}^{13}$ has only
Gorenstein terminal singularities with the following descriptions: 

\begin{enumerate}[$(1)$]

\item The singular locus of $\mathscr{H}_{\mA}^{13}$ coincides with
$\Delta\cup\mA_{\mathsf{P}}$, where $\Delta$ is defined as in (\ref{eq:Delta}).

\item The singularities along the $6$-dimensional locus $\Delta$
are $c(G(2,5))$-singularities. 

\item There exists a primitive $K$-negative divisorial extraction
$f\colon\widetilde{\mathscr{H}}\to\mathscr{H}_{\mA}^{13}$ such that 

\begin{enumerate}[$(a)$]

\item singularities of $\widetilde{\mathfrak{\mathscr{H}}}$ are
only $c(G(2,5))$-singularities along the strict transforms of the
closure $\overline{\Delta}$ of $\Delta$, and 

\item for the $f$-exceptional divisor $E_{\mathscr{H}}$, the morphism
$f|_{E_{\mathscr{H}}}$ can be identified with the $\mP^{2}\times\mP^{2}$-fibration
$\rho_{\mathfrak{\mathscr{H}}}\colon\widehat{\mathscr{H}}\to\mA_{\mathsf{P}}$
as in Subsection \ref{subsec:-fibration-associated-toP2P2}.

\end{enumerate}

\end{enumerate}
\end{prop}

\begin{proof}
By Proposition \ref{prop:The-singular-locus}, the singular locus
of $\mathscr{H}_{\mA}^{13}$ is contained in $\Delta\cup\mA_{\mathsf{P}}$
and the singularities of $\mathscr{H}_{\mA}^{13}$ along $\Delta$
are $c(G(2,5))$-singularities. It remain to describe the singularities
of $\mathscr{H}_{\mA}^{13}$ along $\mA_{\mathsf{P}}.$\vspace{3pt}

\noindent\textbf{Step 1.} \textit{Let $\widetilde{\mathscr{H}}$
be the variety obtained by blowing up $\mathscr{H}_{\mA}^{13}$ along
$\mA_{\mathsf{P}}$ and taking the reduced structure of the blow-up.
We denote by $f\colon\widetilde{\mathscr{H}}\to\mathscr{H}_{\mA}^{13}$
the natural induced morphism. We show that the assertion (3) (a) holds.}

A crucial fact is that any equation of $\mathscr{H}_{\mA}^{13}$ is
of degree $2$ if we regard $p_{ijk}$ as constants. Therefore, $\widetilde{\mathscr{H}}$
has the same singularities as the $x_{ij}$- and $u_{j}$-charts $(i=1,2,j=1,2,3)$
of $\mathscr{H}_{\mA}^{13}$, and then, by Subsection \ref{subsec:Descriptrion-of-charts},
$\widetilde{\mathscr{H}}$ has also only $c(G(2,5))$-singularities,
which are known to be Gorenstein terminal singularities. Since $\widetilde{\mathscr{H}}$
has only $c(G(2,5))$-singularities, the singular locus of $\widetilde{\mathscr{H}}$
is a smooth $6$-dimensional closed subset. Moreover, it contains
the strict transform of $\overline{\Delta}$ and, by the descriptions
of the charts of $\widetilde{\mathscr{H}}$, the restriction of the
singular locus of $\widetilde{\mathscr{H}}$ to the $f$-exceptional
divisor is $5$-dimensional. Therefore, the singular locus of $\widetilde{\Sigma}$
coincides with the strict transform of $\mathsf{\overline{\Delta}}$.

\vspace{3pt}

\noindent\textbf{Step 2.} \textit{We show that $\mathscr{H}_{\mA}^{13}$
has only Gorenstein terminal singularities also along $\mA_{\mathsf{P}}$.}

Let $E$ be the exceptional divisor of the blow-up $\widetilde{\mA}\to\mA_{\mathscr{H}}^{17}$
along $\mA_{\mathsf{P}}$. Let $X_{ij},U_{j}$ be the projective coordinates
of $\widetilde{\mA}$ corresponding to $x_{ij},u_{j}$. We see that
the $f$-exceptional divisor $E_{\mathscr{H}}=E\cap\widetilde{\mathscr{H}}$
is defined with the equation of $\mathscr{H}_{\mA}^{13}$ by replacing
$x_{ij},u_{j}$ with $X_{ij},U_{j}$. Hence $E_{\mathscr{H}}\to\mA_{\mathsf{P}}$
can be identified with $\rho_{\mathscr{H}}\colon\widehat{\mathscr{H}}\to\mA_{\mathsf{P}}$.
Therefore $E_{\mathscr{H}}$ is irreducible and reduced since so is
$\mathscr{H}_{\mA}^{13}$ by Proposition \ref{prop:9times16 H} (3).
By Proposition \ref{prop:9times16 H} (2) and (4), the variety $\mathscr{H}_{\mA}^{13}$
is normal and Gorenstein, and by Step 1, so is $\widetilde{\mathscr{H}}$.
Thus we may write 
\begin{equation}
K_{\widetilde{\mathscr{H}}}=f^{*}K_{\mathscr{H}_{\mA}^{13}}+aE_{\mathscr{H}}\label{eq:SigmaAdj}
\end{equation}
with some integer $a$. Note that a general fiber of $E_{\mathscr{H}}\to\mA_{\mathsf{P}}$
is $\mP^{2}\times\mP^{2}$ by Proposition \ref{prop:Let--beP2P2fib}.
Therefore we see that $a=2$ restricting (\ref{eq:SigmaAdj}) to a
general fiber. This implies that $\mathscr{H}_{\mA}^{13}$ has only
Gorenstein terminal singularities also along $\mA_{\mathsf{P}}$ since
so does $\widetilde{\mathscr{H}}$ by Step 1.

\noindent\textbf{Step 3.} \textit{Finally we show that $f$ is a
$K$-negative primitive divisorial extraction.}

Note that $-E_{\mathscr{H}}$ is $f$-ample by a general property
of blow-up. Thus $f$ is $K$-negative by (\ref{eq:SigmaAdj}) with
$a=2$. Since $\mathscr{H}_{\mA}^{13}$ is affine and $-E_{\mathscr{H}}$
is $f$-ample, we may find an effective divisor $D$ such that $D$
is linearly equivalent to $n(-E_{\mathscr{H}})$ for sufficiently
big $n\in\mN$ and $(\widetilde{\mathscr{H}},\nicefrac{1}{n}D)$ is
a klt pair since $\widetilde{\mathscr{H}}$ has only terminal singularities.
Note that $\widetilde{\mathscr{H}}$ is locally factorial since so
is a $c(G(2,5))$-singularity. Therefore, by \cite[Cor.~1.4.1 and 1.4.2]{bchm},
we may run the $K$-MMP over $\mathscr{H}_{\mA}^{13}$ with scale
$\nicefrac{1}{n}D$. By the choice of the scale, the strict transforms
of $E_{\mathscr{H}}$ are always negative for the contractions when
we run the $K$-MMP. Therefore the final contraction of the $K$-MMP
is the contraction of the strict transforms of $E_{\mathscr{H}}$.
We denote by $f'$ this contraction, by ${\mathscr{H}}'$ the source
of $f'$ and by ${E}'_{\mathscr{H}}$ the strict transforms of $E_{\mathscr{H}}$
on $\mathscr{H}'$. Since $R_{\mathscr{H}}$ is a UFD by Proposition
\ref{prop:UFD}, we see that the target of $f'$ is equal to $\mathscr{H}_{\mA}^{13}$
and $f'({E}'_{\mathscr{H}})=\mA_{\mathsf{P}}$. Since the $K$-MMP
induces a birational map $\widetilde{\mathscr{H}}\dashrightarrow\mathscr{H}'$
which is isomorphism in codimension $1$, $f$ and $f'$ have the
common target $\mathscr{H}_{\mA}^{13}$, and both $f$ and $f'$ are
$K$-negative contractions, we conclude that the induced map $\widetilde{\mathscr{H}}\dashrightarrow\mathscr{H}'$
is actually an isomorphism by \cite[Lem.~5.5]{Tak2} for example. 
\end{proof}

\section{\textbf{$\mQ$-Fano $3$-fold of type No.$\,$5.4\label{sec:-Fano--fold-of 5.4}}}

This section is devoted to the proof of Theorem \ref{thm:embthm}.
We begin with some comments on the part of the statement of Theorem
\ref{thm:embthm} (1) concerning the $\mQ$-Fano $3$-folds of Type
No. 5.4 with respect to a $\nicefrac{1}{2}(1,1,1)$-singularity. This
claim is about the Sarkisov link constructed from such a $\mQ$-Fano
$3$-fold. Although one can give a short proof by invoking \cite[Subsec.5.4]{Ca},
here we provide a proof based on the Sarkisov link constructed from
the key variety. The construction is carried out in Subsections \ref{subsec:A-toric-Sarkisov}
and \ref{subsec:The-Sarkisov-link H}. Afterwards, we proceed step
by step to the proofs of Theorem \ref{thm:embthm} (1) and (2) in
Subsections \ref{subsec:An-example-of fano 1} and \ref{subsec:Proof-of-Theorem 1.3 (2)}
respectively.

\subsection{The $\mQ$-Fano variety $\mathfrak{\mathscr{H}}_{\mP}^{12}$}

We may verify that $\mathfrak{\mathscr{H}}_{\mathbb{A}}^{13}$ is
positively graded by the following weights of coordinates:

\begin{equation}
\begin{cases}
w(\bm{x}_{i})=\left(\begin{array}{c}
1\\
1
\end{array}\right)\,(i=1,2,3),\\
w(p_{ijk})=1\,(i,j,k=1,2),\,w(u_{i})=2\,(i=1,2,3).
\end{cases}\label{eq:wH}
\end{equation}

In this section, let 
\begin{align*}
 & \mathfrak{\mathscr{H}}_{\mP}^{12}\subset\mP(1^{14},2^{3})
\end{align*}
be the weighted projectivization of $\mathfrak{\mathscr{H}}_{\mathbb{A}}^{13}$
by these weights of coordinates. We denote $\mP(1^{14},2^{3})$ by
$\mP$ for simplicity. 
\begin{prop}
\label{prop:HQFano}$\mathfrak{\mathscr{H}}_{\mP}^{12}$ is a $\mQ$-factorial
$\mQ$-Fano variety with Picard number $1$ satisfying $-K_{\mathfrak{\mathscr{H}}_{\mP}^{12}}=\sO_{\mP}(10).$
\end{prop}

\begin{proof}
By Corollary \ref{HFactorial}, $\mathfrak{\mathscr{H}}_{\mP}^{12}$
is a $\mQ$-factorial variety with Picard number $1$. By Proposition
\ref{prop:HTerm}, $\mathfrak{\mathscr{H}}_{\mA}^{13}$ has only Gorenstein
terminal singularities. Note that $\mathfrak{\mathscr{H}}_{\mP}^{12}$
is the $\mC^{*}$-quotient of $\mathfrak{\mathscr{H}}_{\mA}^{13}\setminus\{\bm{o}\}$
and the $\mC^{*}$-action on $\mathfrak{\mathscr{H}}_{\mA}^{13}\setminus\{\bm{o}\}$
is free outside the inverse images of the $u_{i}$-points ($1\leq i\leq3)$.
Since the $u_{i}$-charts are smooth as we have seen in Subsection
\ref{subsec:Descriptrion-of-charts}, we see that $\mathfrak{\mathscr{H}}_{\mP}^{12}$
has $\nicefrac{1}{2}(1^{12})$-singularities at the $u_{i}$-points
checking the $\mC^{*}$-action near the inverse images of the $u_{i}$-points.
Therefore, $\mathfrak{\mathscr{H}}_{\mP}^{12}$ has only terminal
singularities. By (\ref{eq:omegaP}), we have $-K_{\mathscr{H}_{\mP}^{12}}=\sO_{\mathscr{H}_{\mP}^{12}}(10)$.
Therefore, $\mathfrak{\mathscr{H}}_{\mP}^{12}$ is a $\mQ$-Fano variety.
\end{proof}
Let $\widehat{f}_{\mathscr{H}}\colon\mathscr{H}_{1}\to\mathscr{H}_{\mA}^{13}$
be the blow-up at the $u_{1}$-point. Following \cite{BZ} and \cite{AZ},
we embed $\mathscr{H}_{1}$ in a toric variety $\mathscr{T}_{1}$
of Picard number $2$, and then playing the so-called 2-ray game,
we construct the Sarkisov link starting from $\widehat{f}_{\mathscr{H}}$.
Here we chose the $u_{1}$-point as the center of the blow-up, but
by the natural $\mathfrak{S}_{3}$-action on $\mathscr{H}_{\mP}^{12}$,
we could equally take the $u_{2}$- or $u_{3}$-point.

\subsection{A toric Sarkisov link\label{subsec:A-toric-Sarkisov}}

Recall that $S_{\mathscr{H}}$ is the polynomial ring defined as in
Proposition \ref{prop:9times16 H}. We take the polynomial ring $S_{\mathscr{H}}[v]$
with one more variable $v$. Then let us consider the $(\mC^{*})^{2}$-action
on the polynomial ring $S_{\mathscr{H}}[v]$ with the following weights
on the variables: {\footnotesize{
\begin{equation}
\label{eq:wtmat}
\setlength{\arraycolsep}{6pt}
\begin{array}{@{} c *{8}{c} c @{}}
   & w(v) & w(u_1) & w(\bm{x}_2) & w(\bm{x}_3)
   & w(p_{ijk})\, (i,j,k=1,2) & w(u_2) & w(u_3) & w(\bm{x}_1) & \\[4pt]
  \multirow{2}{*}{$\Bigl(\!$}
   & 0  & 2  & 1  & 1  & 1  & 2  & 2  & 1 & \multirow{2}{*}{$\!\Bigr),$} \\
   & -1 & -1 & 0  & 0  & 0  & 1  & 1  & 1 &
\end{array}
\end{equation}}}where the two rows of the matrix record the weights of the variables
with respect to the two factors of $(\mC^{*})^{2}$; the two entries
under the symbol $w(\bm{x}_{i})$ indicate the common weights of the
two entries of $\bm{x}_{i}$, and the two entries under the symbol
$w(p_{ijk})\,(i,j,k=1,2)$ represent the common weights of all $p_{ijk}$,
independent of the indices. Note that the weights of the variables
except $v$ in the 1st row of the matrix come from (\ref{eq:wH}).

Let $\mathscr{T}_{1}$ be the toric variety whose Cox ring is $S_{\mathscr{H}}[v]$
with the $(\mC^{*})^{2}$-weights on the variables as in (\ref{eq:wtmat})
and whose unstable locus as for the quotient construction is 
\[
\Lambda_{1}:=\{v=u_{1}=0\}\cup\{\bm{x}_{2}=\bm{x}_{3}=\bm{o},p_{ijk}=u_{2}=u_{3}=0,\bm{x}_{1}=\bm{o}\}.
\]
By the general theory on toric GIT (cf. \cite[Sect.4]{BZ}, \cite[Chap.14 and 15]{CLS}),
we can read off from the matrix (\ref{eq:wtmat}) the following toric
Sarkisov link:

\begin{equation}\label{eq:Sarkisov} \xymatrix{& \mathscr{T}_1\ar@{-->}[rr]\ar[dl]_{\widehat{f}_{\mathscr{T}}}\ar[dr] & & \mathscr{T}_2\ar@{-->}[rr]\ar[dl]\ar[dr] & & \mathscr{T}_3\ar[dr]^{\widetilde{f}_{\mathscr{T}}}\ar[dl]\\
\mP& &\overline{\mathscr{T}}_1& &\overline{\mathscr{T}}_2 & & \mP^1,}
\end{equation}where $\mathscr{T}_{i}\,(1=1,2,3)$ are $\mQ$-factorial toric varieties
of Picard number $2$. In what follows, we explain the construction
of this Sarkisov link. We shall not mention it each time, but the
description of the exceptional loci of birational maps is due to \cite[Lem.15.3.11]{CLS}
and \cite[Lem.4.5]{BZ}.

The morphism $\widehat{f}_{\mathscr{T}}\colon\mathscr{T}_{1}\to\mP$
is defined by a linear system corresponding to the ray spanned by
the 2nd column {\scriptsize{$\begin{pmatrix} 2\\-1 \end{pmatrix}$}}
of the matrix (\ref{eq:wtmat}). From the matrix, we can read off
that it is given by
\begin{equation}
(v,u_{1},\bm{x}_{2},\bm{x}_{3},p_{ijk},u_{2},u_{3},\bm{x}_{1})\mapsto(u_{1},v^{\nicefrac{1}{2}}\bm{x}_{2},v^{\nicefrac{1}{2}}\bm{x}_{3},v^{\nicefrac{1}{2}}p_{ijk},v^{2}u_{2},v^{2}u_{3},v^{\nicefrac{3}{2}}\bm{x}_{1}),\label{eq:wtblup}
\end{equation}
where we use the same symbols $(u_{1},\bm{x}_{2},\bm{x}_{3},p_{ijk},u_{2},u_{3},\bm{x}_{1})$
for the homogeneous coordinates of $\mP=\mP(1^{14},2^{3})$ as before;
thus the correspondence of coordinates is to be understood in the
natural way. The exceptional locus of $\mathscr{T}_{1}\to\mP$ coincides
with the divisor $\widehat{E}_{\mathscr{T}}:=\{v=0\}$ and this is
mapped to the $u_{1}$-point.

The morphism $\mathscr{T}_{1}\to\overline{\mathscr{T}}_{1}$ is defined
by a linear system corresponding to the ray spanned by the column
{\scriptsize{$\begin{pmatrix} 1\\0 \end{pmatrix}$}} of the matrix
(\ref{eq:wtmat}). From the matrix, we can read off that it is given
by
\begin{align}
(v,u_{1},\bm{x}_{2},\bm{x}_{3},p_{ijk},u_{2},u_{3},\bm{x}_{1})\mapsto(\bm{x}_{2},\bm{x}_{3},p_{ijk},vu_{2},vu_{3},v\bm{x}_{1},u_{1}u_{2},u_{1}u_{3},u_{1}\bm{x}_{1}),\label{eq:F1F1}
\end{align}
where the ambient weighted projective space of $\overline{\mathscr{T}}_{1}$
is $\mP(1^{14},2^{2},3^{2},4^{2})$, and the coordinates corresponding
to $\bm{x}_{2},\bm{x}_{3},p_{ijk}$ are denoted by the same symbols,
those corresponding to $vu_{2},vu_{3},v\bm{x}_{1}$ are denoted by
$u_{2},u_{3},\bm{x}_{1}$, and those corresponding to $u_{1}u_{2},u_{1}u_{3},u_{1}\bm{x}_{1}$
(with weights $4,4,3,3$) are denoted by $u'_{2},u'_{3},\bm{x}_{1}'$.
The exceptional locus of $\mathscr{T}_{1}\to\overline{\mathscr{T}}_{1}$
coincides with $\{u_{2}=u_{3}=0,\bm{x}_{1}=\bm{o}\}$. By considering
the matrix (\ref{eq:wtmat}) except the right $3$ columns, we see
that the restriction of $\mathscr{T}_{1}\to\overline{\mathscr{T}}_{1}$
to the exceptional locus is a $\mP^{1}$-bundle over $\mP^{11}=\mP(\bm{x}_{2},\bm{x}_{3},p_{ijk})$
by \cite[Def.3.1]{Ah} (in order to match the matrix with the format
of ibid., its columns need to be appropriately reordered and generators
of $(\mC^{*})^{2}$ need to be changed. Similar situations will appear
later, but we shall omit the explanation). 

Let $\mathscr{T}_{2}$ be the toric variety whose Cox ring is $S_{\mathscr{H}}[v]$
with the $(\mC^{*})^{2}$-weights on the coordinates as in (\ref{eq:wtmat})
and whose unstable locus as for the quotient construction is 
\[
\Lambda_{2}:=\{v=u_{1}=0,\bm{x}_{2}=\bm{x}_{3}=\bm{o},p_{ijk}=0\}\cup\{u_{2}=u_{3}=0,\bm{x}_{1}=\bm{o}\}.
\]
The morphism $\mathscr{T}_{2}\to\overline{\mathscr{T}}_{1}$ is given
by (\ref{eq:F1F1}). The exceptional locus of $\mathscr{T}_{2}\to\overline{\mathscr{T}}_{1}$
coincides with $\{v=u_{1}=0\}$. By considering the matrix (\ref{eq:wtmat})
except the left $2$ columns, we see that the restriction of $\mathscr{T}_{2}\to\overline{\mathscr{T}}_{1}$
to the exceptional locus is a $\mP^{3}$-bundle over $\mP^{11}=\mP(\bm{x}_{2},\bm{x}_{3},p_{ijk})$
by \cite[Def.3.1]{Ah}. By the explanation of \cite[p.29]{BZ} for
example, the rational map $\mathscr{F}_{1}\dashrightarrow\mathscr{F}_{2}$
is a $(-1^{2},1^{4})$-flip.

Subtracting twice the 2nd row from the 1st row of the matrix in (\ref{eq:wtmat}),
we obtain the following matrix:

{\footnotesize{
\begin{equation}
\label{eq:wtmat2}
\setlength{\arraycolsep}{6pt}
\begin{array}{@{} c *{8}{c} c @{}}
   & w(v) & w(u_1) & w(\bm{x}_2) & w(\bm{x}_3)
   & w(p_{ijk})\, (i,j,k=1,2) & w(u_2) & w(u_3) & w(\bm{x}_1) & \\[4pt]
  \multirow{2}{*}{$\Bigl(\!$}
   & 2  & 4  & 1  & 1  & 1  & 0  & 0  & -1 & \multirow{2}{*}{$\!\Bigr),$} \\
   & -1 & -1 & 0  & 0  & 0  & 1  & 1  & 1 &
\end{array}
\end{equation}}}

The morphism $\mathscr{T}_{2}\to\overline{\mathscr{T}}_{2}$ is defined
by a linear system corresponding to the ray spanned by the column
{\scriptsize{$\begin{pmatrix} 0\\1 \end{pmatrix}$}} of the matrix
(\ref{eq:wtmat2}). From the matrix, we can read off that it is given
by
\begin{align}
 & (v,u_{1},\bm{x}_{2},\bm{x}_{3},p_{ijk},u_{2},u_{3},\bm{x}_{1})\mapsto\label{eq:F2F2}\\
 & (x_{l1}^{2}v,x_{l1}^{4}u_{1},x_{l1}\bm{x}_{2},x_{l1}\bm{x}_{3},x_{l}p_{ijk},u_{2},u_{3})\,(l=1,2).\nonumber 
\end{align}
The exceptional locus of $\mathscr{T}_{2}\to\overline{\mathscr{T}}_{2}$
coincides with $\{\bm{x}_{1}=\bm{o}\}$, which is a $\mP(1^{12},2,4)$-bundle
over $\mP^{1}=\mP(u_{2},u_{3})$ by \cite[Def.3.1]{Ah}.

Let $\mathscr{T}_{3}$ be the toric variety whose Cox ring is $S_{\mathscr{H}}[v]$
with the $(\mC^{*})^{2}$-weights on the coordinates as in (\ref{eq:wtmat2})
and whose unstable locus as for the quotient construction is 
\[
\Lambda_{3}:=\{v=u_{1}=0,\bm{x}_{2}=\bm{x}_{3}=\bm{o},p_{ijk}=u_{2}=u_{3}=0\}\cup\{\bm{x}_{1}=\bm{o}\}.
\]
The morphism $\mathscr{T}_{3}\to\overline{\mathscr{T}}_{2}$ is given
by (\ref{eq:F2F2}). The exceptional locus of $\mathscr{T}_{3}\to\overline{\mathscr{T}}_{2}$
coincides with $\{v=u_{1}=0,\bm{x}_{2}=\bm{x}_{3}=\bm{o},p_{ijk}=0\}$,
which is a $\mP^{1}$-bundle over $\mP^{1}=\mP(u_{2},u_{3})$ by \cite[Def.3.1]{Ah}.
From the matrix (\ref{eq:wtmat2}), we see the rational map $\mathscr{F}_{1}\dashrightarrow\mathscr{F}_{2}$
is a $(2,4,1^{12},-1^{2})$-flip.

Subtracting the 2nd row from the 1st row of the matrix in (\ref{eq:wtmat}),
and then interchanging the two rows, we obtain the following matrix:

{\footnotesize{
\begin{equation}
\label{eq:wtmat3}
\setlength{\arraycolsep}{6pt}
\begin{array}{@{} c *{8}{c} c @{}}
   & w(v) & w(u_1) & w(\bm{x}_2) & w(\bm{x}_3)
   & w(p_{ijk})\, (i,j,k=1,2) & w(u_2) & w(u_3) & w(\bm{x}_1) & \\[4pt]
  \multirow{2}{*}{$\Bigl(\!$}
   & -1  & -1  & 0  & 0  & 0  & 1  & 1  & 1 & \multirow{2}{*}{$\!\Bigr).$} \\
   & 1 & 3 & 1  & 1  & 1  & 1  & 1  & 0 &
\end{array}
\end{equation}}}The morphism $\mathscr{T}_{3}\to\mP^{1}$ is defined by a linear system
corresponding to the ray spanned by the column {\scriptsize{$\begin{pmatrix} 1\\0 \end{pmatrix}$}}
of the matrix (\ref{eq:wtmat3}). From the matrix, we can read off
that it is given by 
\[
(v,u_{1},\bm{x}_{2},\bm{x}_{3},p_{ijk},u_{2},u_{3},\bm{x}_{1})\mapsto(\bm{x}_{1}).
\]
 By \cite[Def.3.1]{Ah}, the morphism $\mathscr{T}_{3}\to\mP^{1}$
is a $\mP(1^{15},3)$-bundle.

\subsection{The Sarkisov link associated to $\mathscr{H}_{\mP}^{12}$ and the
$u_{1}$-point\label{subsec:The-Sarkisov-link H}}

Let $\mathscr{H}_{i}\,(i=1,2,3)$ be the subvariety of $\mathscr{T}_{i}$
with the same equation as $\mathscr{H}_{\mP}^{12}$ (see (\ref{eq:HypMat})).
We can directly check that this is well-defined by (\ref{eq:HypMat})
and the matrix (\ref{eq:wtmat}). Let $C(\mathscr{H}_{\mA}^{13})$
be the cone over $\mathscr{H}_{\mA}^{13}$ in ${\rm Spec}\,S_{\mathscr{H}}[v]$.
Since $\mathscr{T}_{i}$ is the geometric quotient of an affine open
subset ${\rm Spec}\,S_{\mathscr{H}}[v]\setminus\Lambda_{i}\subset{\rm Spec}\,S_{\mathscr{H}}[v]$
(cf.~\cite[Sect.5.1]{CLS}), $\mathscr{H}_{i}$ is the geometric
quotient of the affine open subset $C(\mathscr{H}_{\mA}^{13})\setminus(C(\mathscr{H}_{\mA}^{13})\cap\Lambda_{i})\subset\mathscr{H}_{\mA}^{13}$.
Since $\mathscr{H}_{\mA}^{13}$ is normal by Proposition \ref{prop:9times16 H}
(4), $\mathscr{H}_{i}$ is also normal. By Proposition \ref{prop:UFD},
the affine coordinate ring of $C(\mathscr{H}_{\mA}^{13})\setminus(C(\mathscr{H}_{\mA}^{13})\cap\Lambda_{i})$
is a UFD. Therefore, once we show that $C(\mathscr{H}_{\mA}^{13})\cap\Lambda_{i}$
is of codimension $\geq2$ in $C(\mathscr{H}_{\mA}^{13})$, we can
conclude that $\mathscr{H}_{i}$ is a Mori dream space by \cite[Cor.2.4]{HK},
and in particular, it is $\mQ$-factorial. 

To check that $C(\mathscr{H}_{\mA}^{13})\cap\Lambda_{i}$ is of codimension
$\geq2$ in $C(\mathscr{H}_{\mA}^{13})$ for $i=1,2,3$, we have only
to show the following:
\begin{lem}
\label{lem:-codim}$C(\mathscr{H}_{\mA}^{13})\cap\{v=u_{1}=0\}$ and
$C(\mathscr{H}_{\mA}^{13})\cap\{\bm{x}_{1}=\bm{o}\}$ are of codimension
$\geq2$ in $C(\mathscr{H}_{\mA}^{13})$.
\end{lem}

\begin{proof}
Since $C(\mathscr{H}_{\mA}^{13})\cap\{v=u_{1}=0\}=\mathscr{H}_{\mA}^{13}\cap\{u_{1}=0\}$
and the latter is of codimension $1$ in $\mathscr{H}_{\mA}^{13}$,
$C(\mathscr{H}_{\mA}^{13})\cap\{v=u_{1}=0\}$ is of codimension $2$
in $C(\mathscr{H}_{\mA}^{13})$. 

Note that $C(\mathscr{H}_{\mA}^{13})\cap\{\bm{x}_{1}=\bm{o}\}$ is
the union of the following three loci:

\begin{align*}
\Gamma_{1} & :=\{u_{2}=u_{3}=0,\bm{x}_{1}=\bm{o},D^{(3)}(x)\bm{x}_{2}=\bm{o},|D^{(2)}(x)|=|D^{(3)}(x)|=0\},\\
\Gamma_{2} & :=\{u_{2}=0,u_{3}\not=0,\bm{x}_{1}=\bm{x}_{3}=\bm{o},u_{1}u_{3}+|D^{(2)}(x)|=0\},\\
\Gamma_{3} & :=\{u_{2}\not=0,u_{3}=0,\bm{x}_{1}=\bm{x}_{2}=\bm{o},u_{1}u_{2}+|D^{(3)}(x)|=0\}.
\end{align*}

The loci $\Gamma_{2}$ and $\Gamma_{3}$ are clearly of codimension
$2$ in $C(\mathscr{H}_{\mA}^{13})$. Since it is easy to check that
$\Gamma_{1}\cap\{\bm{x}_{2}=\bm{o}\}$ is of codimension 3, we may
assume that $\bm{x}_{2}\not=\bm{o}$. Then $|D^{(3)}(x)|=0$ follows
from $D^{(3)}(x)\bm{x}_{2}=\bm{o}.$ Since it holds that $D^{(3)}(x)\bm{x}_{2}=D^{(2)}(x)\bm{x}_{3}$,
we may check that 
\[
\Gamma_{1}\cap\{\bm{x}_{2}\not=\bm{o}\}=\{u_{2}=u_{3}=0,\bm{x}_{1}=\bm{o},\bm{x}_{2}\not=\bm{o}\}\cap\left\{ \rank\left(\begin{array}{cc}
(D^{(2)}(x))^{\dagger} & \bm{x}_{3}\end{array}\right)\leq1\right\} ,
\]
where $(D^{(3)}(x))^{\dagger}$ is the adjoint matrix of $D^{(3)}(x)$.
It is easy to check that this is of codimension $2$ in $C(\mathscr{H}_{\mA}^{13})$
by a property of a determinantal variety. Therefore $C(\mathscr{H}_{\mA}^{13})\cap\{\bm{x}_{1}=\bm{o}\}$
are of codimension $2$ in $C(\mathscr{H}_{\mA}^{13})$.
\end{proof}
We shall explain how restricting (\ref{eq:Sarkisov}) yields the following
diagram:

\begin{equation}\label{eq:SarkisovH} \xymatrix{& \mathscr{H}_1\ar@{-->}[rr]\ar[dl]_{\widehat{f}_{\mathscr{H}}}\ar[dr] & & \mathscr{H}_2\ar@{-->}[rr]\ar[dl]\ar[dr] & & \mathscr{H}_3\ar[dr]^{\widetilde{f}_{\mathscr{H}}}\ar[dl]\\
\mathscr{H}_{\mP}^{12}& &\overline{\mathscr{H}}_1& &\overline{\mathscr{H}}_2 & & \mP^1.}
\end{equation}

We can easily check that $\widehat{f}_{\mathscr{T}}(\mathscr{H}_{1})=\mathscr{H}_{\mP}^{12}$.
We denote by $\widehat{f}_{\mathscr{H}}\colon\mathscr{H}_{1}\to\mathscr{H}_{\mP}^{12}$
the restriction of $\widehat{f}_{\mathscr{T}}$. By the description
of $\widehat{f}_{\mathscr{T}}$, the exceptional locus of $\widehat{f}_{\mathscr{H}}$
is $\widehat{E}_{\mathscr{H}}:=\mathscr{H}_{1}\cap\{v=0\}$, which
is mapped to the $u_{1}$-point. Setting $u_{1}=1$ in (\ref{eq:wtblup}),
we can describe $\widehat{f}_{\mathscr{H}}$ near $u_{1}$-point.
Indeed, by $\bm{G}_{1},\mathsf{G}_{5},\mathsf{G}_{6}$ with setting
$u_{1}=1$, we can erase the coordinates $x_{11},x_{21},u_{2},u_{3}$,
and then $\widehat{f}_{\mathscr{H}}$ is reduced to the morphism defined
by 
\[
(v,\bm{x}_{2},\bm{x}_{3},p_{ijk})\mapsto(v^{\nicefrac{1}{2}}\bm{x}_{2},v^{\nicefrac{1}{2}}\bm{x}_{3},v^{\nicefrac{1}{2}}p_{ijk}),
\]
which is nothing but the weighted blow-up at the $u_{1}$-point (a
$\nicefrac{1}{2}(1^{12})$-singularity) with weight $((\nicefrac{1}{2})^{12})$,
and coincides with the simple blow-up at the $u_{1}$-point. In particular,
since $\mathscr{H}_{\mP}^{12}$ is $\mQ$-factorial and has Picard
number $1$, $\mathscr{H}_{1}$ is $\mQ$-factorial and has Picard
number $2$. Since $\mathscr{H}_{\mP}^{12}$ has only terminal singularities,
so does $\mathscr{H}_{1}$. By the adjunction formula for this weighted
blow-up, we have
\begin{equation}
-K_{\mathscr{H}_{1}}=\widehat{f}_{\mathscr{H}}^{*}\sO_{\mathscr{H}_{\mP}^{12}}(10)-5\widehat{E}_{\mathscr{H}}=5(\widehat{f}_{\mathscr{H}}^{*}\sO_{\mathscr{H}_{\mP}^{12}}(2)-\widehat{E}_{\mathscr{H}}).\label{eq:-KH1}
\end{equation}

Let $\overline{\mathscr{H}}_{1}$ be the image of $\mathscr{H}_{1}$
by the morphism $\mathscr{T}_{1}\to\overline{\mathscr{T}}_{1}$. By
(\ref{eq:F1F1}) and the equation $\bm{G}_{1}=\bm{o}$, $\mathsf{G}_{5}=\mathsf{G}_{6}=0$,
we can consider $\overline{\mathscr{H}}_{1}$ is contained in $\mP(1^{14},2^{2})$
by erasing the coordinates $u_{2}'$, $u'_{3}$ and $\bm{x}_{1}'$
and can check that $\overline{\mathscr{H}}_{1}$ is defined by $\bm{G}_{2}=\bm{G}_{3}=\bm{o}$
and $\mathsf{G}_{4}=0$, equivalently, is defined by the five $4\times4$
Pfaffians of the matrix $M_{{\rm T}}$ as in (\ref{eq:MT}). By the
description of the exceptional locus of $\mathscr{T}_{1}\to\overline{\mathscr{T}}_{1}$,
we see that the exceptional locus $E_{1}$ of $\mathscr{H}_{1}\to\overline{\mathscr{H}}_{1}$
and its image $\overline{E}_{1}$ are defined by the same equation
as $\Gamma_{1}$, and $E_{1}\to\overline{E}_{1}$ is a $\mP^{1}$-bundle.
Similarly to Lemma \ref{lem:-codim}, we can check that $E_{1}$ is
of codimension $2$ in $\mathscr{H}_{1}$. Thus $\mathscr{H}_{1}\to\overline{\mathscr{H}}_{1}$
is a small contraction. Moreover, it is a flopping contraction since,
by (\ref{eq:-KH1}), $-K_{\mathscr{H}_{1}}$ corresponds to the ray
spanned by the column {\scriptsize{$\begin{pmatrix} 1\\0 \end{pmatrix}$}}
of the matrix (\ref{eq:wtmat}). 

We can directly check that the image of $\mathscr{H}_{2}$ by the
morphism $\mathscr{T}_{2}\to\overline{\mathscr{T}}_{1}$ coincides
with $\overline{\mathscr{H}}_{1}$. The exceptional locus of the induced
morphism $\mathscr{H}_{2}\to\overline{\mathscr{H}}_{1}$ is $E_{1}^{+}:=\mathscr{H}_{2}\cap\{v=u_{1}=0\}$.
With some calculations, we can check that the image of $E_{1}^{+}$
is $\overline{E}_{1}$ and $\mathscr{H}_{2}\to\overline{\mathscr{H}}_{1}$
is also a small contraction. It is also a flopping contraction by
the same reason for $\mathscr{H}_{1}\to\overline{\mathscr{H}}_{1}$.
Since we see that nontrivial fibers are not equidimensional by direct
calculations, $\mathscr{H}_{2}\to\overline{\mathscr{H}}_{1}$ is not
isomorphic to $\mathscr{H}_{1}\to\overline{\mathscr{H}}_{1}$ over
$\overline{\mathscr{H}}_{1}$. Therefore $\mathscr{H}_{1}\dashrightarrow\mathscr{H}_{2}$
is the flop. This implies that $\mathscr{H}_{2}$ has only terminal
singularities and has Picard number 2 since so does $\mathscr{H}_{1}$.
We omit detailed calculations of fibers of $\mathscr{H}_{2}\to\overline{\mathscr{H}}_{1}$
but we mention that $2$-dimensional fibers appear exactly over points
in $\mP(p_{ijk}\,(i,j,k=1,2))\simeq\mP^{7}\subset\overline{E}_{1}$,
which we need later.

Let $\overline{\mathscr{H}}_{2}$ be the image of $\mathscr{H}_{2}$
by the morphism $\mathscr{T}_{2}\to\overline{\mathscr{T}}_{2}$. By
the descriptions of the exceptional locus of $\mathscr{T}_{2}\to\overline{\mathscr{T}}_{2}$
and the unstable locus $\Lambda_{2}$ as for the quotient construction,
we see that the exceptional locus of $\mathscr{H}_{2}\to\overline{\mathscr{H}}_{2}$
is $\mathscr{H}_{2}\cap\{\bm{x}_{1}=\bm{o}\}=E_{2}^{(2)}\sqcup E_{2}^{(3)}$,
where $E_{2}^{(2)}$ and $E_{2}^{(3)}$ are respectively defined by
the same equations as those of $\Gamma_{2}$ and $\Gamma_{3}$, and
are mapped to the $u_{3}$- and the $u_{2}$-points. Note that $E_{2}^{(2)}$
and $E_{2}^{(3)}$ are weight 4 hypersurfaces in $\mP(1^{10},2,4)$,
and hence they are isomorphic to $\mP(1^{10},2)$. Hence $\mathscr{H}_{2}\to\overline{\mathscr{H}}_{2}$
is a small contraction. Since $\mathscr{H}_{\mP}^{12}$ is a $\mQ$-Fano
variety, $K_{\mathscr{H}_{2}}$ is not nef. Since $K_{\mathscr{H}_{2}}$
is numerically trivial, $K_{\mathscr{H}_{2}}$ must be negative for
$\mathscr{H}_{2}\to\overline{\mathscr{H}}_{2}$. Therefore $\mathscr{H}_{2}\to\overline{\mathscr{H}}_{2}$
is a ($K$-negative) flipping contraction. 

We can directly check that the image of $\mathscr{H}_{3}$ by the
morphism $\mathscr{T}_{3}\to\overline{\mathscr{T}}_{2}$ coincides
with $\overline{\mathscr{H}}_{2}$. We set 

\begin{align*}
E_{2}^{(2)+} & :=\{v=u_{1}=u_{2}=0,u_{3}\not=0,\bm{x}_{2}=\bm{x}_{3}=\bm{o},p_{ijk}=0\},\\
E_{2}^{(3)+} & :=\{v=u_{1}=u_{3}=0,u_{2}\not=0,\bm{x}_{2}=\bm{x}_{3}=\bm{o},p_{ijk}=0\},
\end{align*}
which are isomorphic to $\mP^{1}$. By the descriptions of the exceptional
locus of $\mathscr{T}_{3}\to\overline{\mathscr{T}}_{2}$ and the unstable
locus $\Lambda_{3}$ as for the quotient construction, we see that
the exceptional locus of $\mathscr{H}_{3}\to\overline{\mathscr{H}}_{2}$
is 
\begin{align*}
\mathscr{H}_{3}\cap\{v=u_{1}=0,\bm{x}_{2}=\bm{x}_{3}=\bm{o},p_{ijk}=0\}=E_{2}^{(2)+}\sqcup E_{2}^{(3)+},
\end{align*}
and $E_{2}^{(2)+}$ and $E_{2}^{(3)+}$ are mapped to the $u_{3}$-
and the $u_{2}$-points, respectively. Hence $\mathscr{H}_{3}\to\overline{\mathscr{H}}_{2}$
is also a small contraction. By the descriptions of the exceptional
loci of $\mathscr{H}_{2}\to\overline{\mathscr{H}}_{2}$ and $\mathscr{H}_{3}\to\overline{\mathscr{H}}_{2}$,
these two morphisms are not isomorphic over $\overline{\mathscr{H}}_{2}$.
Therefore $\mathscr{H}_{2}\dashrightarrow\mathscr{H}_{3}$ is the
flip. This implies that $\mathscr{H}_{3}$ has only terminal singularities
and has Picard number 2 since so does $\mathscr{H}_{2}$. 

Let $\widetilde{f}_{\mathscr{H}}\colon\mathscr{H}_{3}\to\mP^{1}$
be the restriction of $\widetilde{f}_{\mathscr{T}}$. Since $\bm{x}_{1}\not=\bm{o}$
by the description of $\Lambda_{3}$, we can eliminate the coordinate
$u_{1}$ by the equality $\bm{G}_{1}=\bm{o}$. If $x_{11}\not=0$,
we may assume that $x_{11}=1$ and then $\mathscr{H}_{3}$ is defined
by the five $4\times4$ Pfaffians of the matrix $M_{{\rm T}}$ as
in (\ref{eq:MT}). We have a similar description of $\mathscr{H}_{3}$
over $\{x_{21}\not=0\}$. Therefore $\widetilde{f}_{\mathscr{H}}\colon\mathscr{H}_{3}\to\mP^{1}$
is a fibration of the $11$-dimensional cone over ${\rm G}(2,5)$. 

\subsection{Proof of Theorem \ref{thm:embthm} (1)\label{subsec:An-example-of fano 1}}

\begin{proof}[Proof of Theorem \ref{thm:embthm} (1)]
 We set 
\begin{equation}
X:=\mathfrak{\mathscr{H}}_{\mP}^{12}\cap\{L_{1}=\cdots=L_{9}=0\},\label{eq:defX}
\end{equation}
where $L_{1},\dots,L_{9}$ are linear forms of the entries of $\bm{x}_{1},\bm{x}_{2},\bm{x}_{3}$,
and $p_{ijk}$. Using the equations of $\mathfrak{\mathscr{H}}_{\mP}^{12}$,
we may easily verify that ${\rm Bs}|\sO_{\mP}(1)|\cap\mathfrak{\mathscr{H}}_{\mP}^{12}$
consists of the $u_{i}$-points~($i=1,2,3)$. Note that, by Propositions
\ref{prop:The-singular-locus} and \ref{prop:HTerm}, $\dim{\rm Sing\,}\mathfrak{\mathscr{H}}_{\mathbb{A}}^{13}$
is less than $9$, which is the number of $L_{1},\dots,L_{9}$. Therefore,
by the Bertini theorem, we see that a general $X$ is a smooth $3$-fold
outside the $u_{i}$-points. Computing the linear parts of the equations
of a general $X$ at each of the $u_{i}$-points (cf. \cite[Sect.2.8]{Tak7}),
we see that $X$ has a $\nicefrac{1}{2}(1,1,1)$-singularity at each
of the $u_{i}$-points. Therefore we have shown the former assertion
of (1). 

Now assume that $X$ is a (not necessarily general) $3$-fold with
only $\nicefrac{1}{2}(1,1,1)$-singularities. Then the $u_{i}$-points
($1\leq i\leq3)$ are only the singularities of $X$, and they are
$\nicefrac{1}{2}(1,1,1)$-singularities. We show that $X$ is a $\mQ$-Fano
$3$-fold anticanonically embedded of codimension $4$. Indeed, $X$
is of codimension $4$ in the ambient weighted projective space since
so is $\mathscr{H}_{\mathbb{P}}^{12}$ and $\dim X=3$. The condition
that $-K_{X}=\sO_{X}(1)$ can be verified by (\ref{eq:omegaP}) and
(\ref{eq:defX}). The homogeneous coordinate ring of $X$ is Gorenstein
since so is that of $\mathscr{H}_{\mP}^{12}$ by Proposition \ref{prop:9times16 H}
(2) and $X$ has the right dimension (cf.$\,$\cite[Prop. 1.3, equiv. of (ii) and (iii)]{R}).
Therefore $X$ is anticanonically embedded of codimension $4$. We
show that $\rho(X)=1$. By \cite[Thm.4.1]{To}, we have only to show
this claim for a general $X$. By Corollary \ref{HFactorial} (2),
it suffices to check that $X\cap\{u_{1}=0\}$ is a prime divisor.
This follows by the Bertini theorem since $\mathscr{H}_{\mP}^{12}\cap\{u_{1}=0\}$
is irreducible by Proposition \ref{prop:9times16 H} (3), and the
fact that ${\rm Bs}|\sO_{\mP}(1)|\cap\mathfrak{\mathscr{H}}_{\mP}^{12}$
consists of the $u_{i}$-points. 

Finally, based on the investigation in Subsection \ref{subsec:The-Sarkisov-link H},
we show that $X$ is of type No.$\,$5.4 as in \cite{Tak1} for any
$\nicefrac{1}{2}(1,1,1)$-singularity $\mathsf{p}$. 

By \cite[Lem.4.1.13]{BH}, we can read off the Hilbert numerator of
$X$ from the restriction of (\ref{eq:minfree}), and then we have
$(-K_{X})^{3}=\nicefrac{11}{2}$. By the Riemann-Roch theorem, we
have $g(X):=h^{0}(X,\sO_{X}(-K_{X}))-2=3.$ Therefore, by \cite{Tak1},
$X$ is of type No.$\,$4.1 or No$\,$5.4.

Note that no $L_{i}$ is not a linear combination of the entries from
only one of $\bm{x}_{1}$, $\bm{x}_{2}$ and $\bm{x}_{3}$. Indeed,
assume, for a contradiction, that, for some $i$, it holds that $L_{i}=\alpha_{1}x_{11}+\alpha_{2}x_{21}$
($\alpha_{1},\alpha_{2}\in\mC$) (by symmetry, we have only to disprove
this situation). Near the $u_{1}$-point, $\mathscr{H}_{\mP}^{12}$
is quasi-smooth with the entries of $\bm{x}_{2}$ and $\bm{x}_{3}$
and $p_{ijk}$ as the orbifold coordinates as we have seen in Subsection
\ref{subsec:Descriptrion-of-charts}. Then, using $\bm{G}_{1}$ in
Corollary \ref{cor:HsharpN}, we see that, near the $u_{1}$-point,
$\mathscr{H}_{\mP}^{12}\cap\{L_{i}=0\}$ is the hypersurface $\left\{ \left(\begin{array}{cc}
\alpha_{1} & \alpha_{2}\end{array}\right)D^{(3)}(x)\bm{x}_{2}=0\right\} $ in the orbifold space with the entries of $\bm{x}_{2}$ and $\bm{x}_{3}$
and $p_{ijk}$ as orbifold coordinates. Therefore $\mathscr{H}_{\mP}^{12}\cap\{L_{i}=0\}$,
hence $X$ cannot be quasi-smooth at the $u_{1}$-point, a contradiction.

By the natural $\mathfrak{S}_{3}$-action on $\mathscr{H}_{\mP}^{12}$,
we may assume that $\mathsf{p}$ is the $u_{1}$-point. By the above
claim, each $L_{i}$ contains one among $p_{ijk}$, entries of $\bm{x}_{2}$,
or entries of $\bm{x}_{3}$. Therefore the strict transform of $\{L_{i}=0\}$
is linearly equivalent to $\widehat{f}_{\mathscr{T}}^{*}\sO_{\mP}(1)-\nicefrac{1}{2}\widehat{E}_{\mathscr{T}}$,
which corresponds to the ray spanned by {\scriptsize{$\begin{pmatrix} 0\\1 \end{pmatrix}$}}
of the matrix (\ref{eq:wtmat}). It is defined by $L_{i}$ with $x_{j1}\,(j=1,2)$
being replaced by $vx_{j1}$, which we denote by $L'_{i}$. Since
$\mathscr{H}_{\mP}^{12}\cap\{L_{i}=0\}$ is quasi-smooth at the $u_{1}$-point,
$\mathscr{H}_{1}\cap\{L'_{i}=0\}$ is the strict transform of $\mathscr{H}_{\mP}^{12}\cap\{L_{i}=0\}$.

Now we set $X_{1}:=\mathscr{H}_{1}\cap\{L'_{1}=\cdots=L'_{9}=0\}.$
Since $X$ is quasi-smooth at the $u_{1}$-point, $X_{1}\to X$ coincides
with the blow-up at the $u_{1}$-point. By \cite[I, Table 4, No.4.1, Table 5, No.5.4]{Tak1},
$X_{1}$ has a flopping contraction, which we denote by $X_{1}\to\overline{X}_{1}$,
thus this contraction must coincides with the restriction of $\mathscr{H}_{1}\to\overline{\mathscr{H}}_{1}$
since $X_{1}$ has only two nontrivial morphisms. The image of the
exceptional locus of $X_{1}\to\overline{X}_{1}$ is $\overline{E}_{1}\cap X_{1}$,
where we recall that $\overline{E}_{1}$ is the image of the exceptional
locus of $\mathscr{H}_{1}\to\overline{\mathscr{H}}_{1}$. Thus $\overline{E}_{1}\cap X_{1}$
consists of a finite number of points.

Let $X_{1}\dashrightarrow X_{2}$ be the flop for the contraction
$X_{1}\to\overline{X}_{1}$. Note that the image of the exceptional
locus of $X_{2}\to\overline{X}_{1}$ is also $\overline{E}_{1}\cap\overline{X}_{1}$.
We denote $\mP(p_{ijk}\,(i,j,k=1,2))$ simply by $\mP(p_{ijk})$.
If $\overline{X}_{1}\cap\mP(p_{ijk})\not=\emptyset$, then $X\cap\mP(p_{ijk})\not=\emptyset$.
This implies that $X$ has Gorenstein terminal singularities along
$X\cap\mP(p_{ijk})$ by Proposition \ref{prop:HTerm}, which is impossible.
Therefore $\overline{X}_{1}\cap\mP(p_{ijk})=\emptyset$. We show that
$X_{2}$ coincides with $X_{2}':=\mathscr{H}_{2}\cap\{L'_{1}=\cdots=L'_{9}=0\}$.
Indeed, since $\mathscr{H}_{2}\to\overline{\mathscr{H}}_{1}$ has
2-dimensional fibers only over points in $\mP(p_{ijk}),$ $X_{2}'\to\overline{X}_{1}$
has only $1$-dimensional nontrivial fibers, thus a small contraction.
Therefore, $X'_{2}$ is normal by \cite[Cor.5.25]{KM}, and we conclude
that $X_{2}=X_{2}'$.

Let $\overline{X}_{2}$ be the image of $X_{2}$ by the contraction
$\mathscr{H}_{2}\to\overline{\mathscr{H}}_{2}$. Since $-K_{\mathscr{H}_{2}}$
is positive for the contraction $\mathscr{H}_{2}\to\overline{\mathscr{H}}_{2}$
and $\{L_{i}'=0\}$ is proportional for $-K_{\mathscr{H}_{2}}$, $X'_{2}$
intersects the exceptional locus $E_{2}^{(2)}\sqcup E_{2}^{(3)}$
of $\mathscr{H}_{2}\to\overline{\mathscr{H}}_{1}$, thus $X_{2}\to\overline{X}_{2}$
is a nontrivial birational morphism. 

At this point, we can conclude that $X$ is of type No.$\,$5.4 as
for the $u_{1}$-point since it cannot be of type No.$\,$4.1 (if
$X$ is of type No.$\,$4.1, then $X_{2}$ has a conic bundle structure
by \cite[I, Table 4, No.4.1]{Tak1} and does not have a birational
morphism different from the flopping contraction $X_{2}\to\overline{X}_{1}$).
We proceed, however, a bit more since we can complete the construction
of the Sarkisov link starting from the blow-up of $X$ at the $u_{1}$-point.
Since $X$ is of type No.$\,$5.4, $X_{2}\to\overline{X}_{2}$ is
a flipping contraction. Let $X_{3}':=\mathscr{H}_{3}\cap\{L'_{1}=\cdots=L'_{9}=0\}$.
Since $-K_{\mathscr{H}_{3}}$ is negative for the contraction $\mathscr{H}_{3}\to\overline{\mathscr{H}}_{2}$
and $\{L_{i}'=0\}$ is proportional for $-K_{\mathscr{H}_{3}}$, $X'_{3}$
contains the exceptional locus $E_{2}^{(2)+}\sqcup E_{2}^{(3)+}$
($E_{2}^{(2)+}\simeq E_{2}^{(3)+}\simeq\mP^{1}$) of $\mathscr{H}_{3}\to\overline{\mathscr{H}}_{2}$,
thus $X_{3}$ is normal by \cite[Cor.5.25]{KM}, and $X_{2}\dashrightarrow X_{3}$
is the flip for $X_{2}\to\overline{X}_{2}$. 

By \cite[I, Table 5, No.5.4]{Tak1}, $X_{3}$ has a fibration of quintic
del Pezzo surfaces, which is nothing but the restriction of $\mathscr{H}_{3}\to\mP^{1}$,
which is a fibration of cones over ${\rm G(2,5}).$
\end{proof}
\begin{rem}
In \cite{Tak4}, we will construct a prime $\mQ$-Fano threefold of
type No.$\,$4.1 via another key variety $\mathfrak{S}_{\mA}^{8}$
obtained in ibid.
\end{rem}

\subsection{Proof of Theorem \ref{thm:embthm} (2)\label{subsec:Proof-of-Theorem 1.3 (2)}}

This subsection is occupied with the proof of Theorem \ref{thm:embthm}
(2). The proof uses results of \cite{Tak1} freely.
\begin{proof}[Proof of Theorem \ref{thm:embthm} (2)]
Let $X$ be a prime $\mathbb{Q}$-Fano 3-fold of type  No.$\,$5.4
as in \cite{Tak1} and assume that $X$ has only $\nicefrac{1}{2}(1,1,1)$-singularities.
We denote by $\mathsf{p}_{1},\mathsf{p}_{2},\mathsf{p}_{3}$ the $\nicefrac{1}{2}(1,1,1)$-singularities
of $X$. Let $g\colon Z\to X$ be the blow-up at (any) two of $\nicefrac{1}{2}(1,1,1)$-singularities,
say, $\mathsf{p}_{1},\mathsf{p}_{2}$. By \cite[II, Thm.1.0 (2)]{Tak1}
and the fact that $(-K_{Z})^{3}=(-K_{X})^{3}-2\times\nicefrac{1}{2}=\nicefrac{9}{2}>0$,
we see that $Z$ is a weak $\mathbb{Q}$-Fano 3-fold. Let $h\colon Z\to W$
be the morphism defined by $|-mK_{Z}|$ with $m\gg$0 ($W$ is so-called
the anticanonical model of $Z$). Note that $W$ is a $\mQ$-Fano
3-fold with only canonical singularities. 

\vspace{3pt}

In what follows, Steps 1--5 are devoted to substantial preparations.
In Step 6 the variety $\mathscr{H}_{\mA}^{13}$ will appear, and the
proof will be completed in Step 7.

\vspace{3pt}

\noindent\textbf{Step 1.} \textit{We show that $W$ is anticanonically
embedded in $\mathbb{P}(1^{5},2)$ as a weighted compelete intersection
of two weight 3 hypersurfaces.}\vspace{3pt}

We will show this applying \cite[II, Thm.3.0]{Tak1} (see the last
part of the proof of ibid.). It suffices to check the conditions (1)--(5)
of ibid. 

\vspace{3pt}

\noindent(4) We show that $h^{0}(\mathcal{O}(-K_{W}))\geq4$. By
the equality $(-K_{Z})=h^{*}(-K_{W})$, we have $h^{0}(\mathcal{O}(-K_{W}))=h^{0}(\mathcal{O}(-K_{Z}))$.
Since $g$ has the minimal discrepancies at the centers $\mathsf{p}_{1},\mathsf{p}_{2}$
and the discrepancies are positive, we have $h^{0}(\mathcal{O}(-K_{Z}))=h^{0}(\mathcal{O}(-K_{X})).$
Therefore, since $h^{0}(\mathcal{O}(-K_{X}))=5$, we have $h^{0}(\mathcal{O}(-K_{W}))=h^{0}(\mathcal{O}(-K_{X}))\geq4.$
\vspace{5pt}

It is convenient below to keep in mind that there are bijections among
the linear systems $|-K_{W}|$, $|-K_{Z}|$ and $|-K_{X}|$ by the
proof of (4).

\vspace{5pt}

\noindent(3) We show that $|-K_{W}|$ has a member with only canonical
singularities. By \cite[II, Cor.1.2]{Tak1}, there is a member $D_{X}$
of $|-K_{X}|$ with only canonical singularities. Then, since $K_{D_{Z}}\leq(g|_{D_{Z}})^{*}K_{D_{X}}$,
we see that the strict transform $D_{Z}$ belongs to $|-K_{Z}|$ and
$D_{Z}$ has also only canonical singularities. The image $D_{W}$
of $D_{Z}$ belongs to $|-K_{W}|$. Moreover, $D_{W}$ has also only
canonical singularities since so is $D_{Z}$ and $K_{D_{Z}}\sim0$.

\vspace{3pt}

\noindent(5) We check that $W$ has only one index $2$ point and
it is a $\nicefrac{1}{2}(1,1,1)$-singularity. It suffices to show
that there are no $h$-exceptional curves through the unique singularity
$\mathsf{q}$ of $Z$ ($\mathsf{q}$ is the point with $g(\mathsf{q})=\mathsf{p}_{3}$).
If there is such an $h$-exceptional curve $l$, then $l\subset{\rm Bs}|-K_{Z}|$
since $-K_{Z}\cdot l=0$ and $\mathsf{q}\in{\rm Bs}|-K_{Z}|$. Then
it holds that $g(l)\subset{\rm Bs}|-K_{X}|$ and $g(l)$ contains
the $\nicefrac{1}{2}(1,1,1)$-singularity $g(\mathsf{q})$, a contradiction
to \cite[II, Thm.1.0 (2)]{Tak1}. 

\vspace{3pt}

\noindent(2) We show that $|-K_{W}|$ has no base curve passing through
its unique index two point. Assume the contrary. Then $|-K_{Z}|$
has also a base curve passing through its unique index two point $\mathsf{q}$
since $h$ is isomorphic near the unique index $2$ point of $W$
by the proof of (5). Since $g$-exceptional divisors do not contain
$\mathsf{q}$, $|-K_{X}|$ has also a base curve passing through $\mathsf{p}_{3}$,
a contradiction to \cite[Thm.1.0 (2)]{Tak1}. 

\vspace{3pt}

\noindent(1) We check that $W$ is indecomposable, i.e., there are
no effective Weil divisors $A,B$ with $h^{0}(A)\geq2$ and $h^{0}(B)\geq2$
such that $-K_{W}\sim A+B$. Assume, for a contradiction, that there
are such $A$ and $B$. By the check of (5), $h$ is isomorphic near
the unique index $2$ point of $W$. Therefore $-K_{W}$ is a Cartier
divisor near the images of $h$-exceptional divisors if exist. Hence
$h^{*}(A+B)$ is an integral divisor. Then we have $h^{*}(A+B)=A'+B'+\sum E_{i}$,
where $A',B'$ are integral Weil Divisors such that $h'_{*}(A')=A,h^{0}(A')\geq h^{0}(A)\geq2,$
and $h'_{*}(B')=B,h^{0}(B')\geq h^{0}(B)\geq2$, and $E_{i}$ are
$h$-exceptional divisors ($E_{i}=E_{j}$ could happen even if $i\not=j$).
Thus we have $-K_{X}\sim g_{*}(A')+g_{*}(B')+\sum g_{*}(E_{i})$.
Since $h^{0}(g_{*}(A')),h^{0}(g_{*}(B'))\geq2$, we have $g_{*}(A')\not=0$
and $g_{*}(B')\not=0$. This contradicts the assumption that $X$
is prime. 

\vspace{3pt}

Now, applying \cite[II, Thm.3.0]{Tak1} to $W$, we see that $W$
is anticanonically embedded in $\mathbb{P}(1^{5},2)$ as a weighted
compelete intersection of two weight 3 hypersurfaces. We denote by
$u_{3}$ the unique coordinate of weight $2$. 

\vspace{5pt}

Before going to Step 2, let $E_{1}$ and $E_{2}$ be the $g$-exceptional
divisors and $\Pi_{1}$ and $\Pi_{2}$ the images on $W$ of $E_{1}$
and $E_{2}$ respectively. We see that $\Pi_{1}$ and $\Pi_{2}$ do
not contain the $\nicefrac{1}{2}(1,1,1)$-singularity of $W$ since
$E_{1}$ and $E_{2}$ do not contain the $\nicefrac{1}{2}(1,1,1)$-singularity
$\mathsf{q}$ of $Z$ and there are no $h$-exceptional curves through
$z$ as we see in the check of (5) in Step 1. Therefore, by the projection
$\mathbb{P}(1^{5},2)\dashrightarrow\mathbb{P}^{4}$ from the $u_{3}$-point,
$\Pi_{1}$ and $\Pi_{2}$ are mapped to planes, say, $\Pi_{1}'$ and
$\Pi_{2}'$ in $\mathbb{P}^{4}$ since $(-K_{Z})^{2}E_{1}=(-K_{Z})^{2}E_{2}=1$.
Therefore $E_{1}$ and $E_{2}$ are mapped to planes in $\mathbb{P}^{4}$
isomorphically and hence $E_{1}$ and $E_{2}$ are mapped to $\Pi_{1}$
and $\Pi_{2}$ isomorphically. \vspace{3pt}

\noindent\textbf{Step 2.} \textit{We show that $h$ has no exceptional
divisors. In particular, $W$ has only terminal singularities. Moreover
the singular locus of $W$ consists of the $u_{3}$-point and the
locus contained in $\Pi_{1}\cup\Pi_{2}$.}\vspace{3pt}

Since an $h$-exceptional curve is numerically trivial for $-K_{Z}$
while $-K_{X}$ is ample, we see that an $h$-exceptional curve intersects
with $E_{1}$ or $E_{2}$. Therefore the singular locus of $W$ is
contained in $\Pi_{1}\cup\Pi_{2}$ outside of the $u_{3}$-point. 

Assume, for a contradiction, that there exists an $h$-exceptional
divisor $F$. Therefore, by the discussion in the above paragraph,
$h(F)$ is contained in $\Pi_{1}$ or $\Pi_{2}$, say, $\Pi_{1}$.
Since $\Pi_{1}$ is mapped to the plane $\Pi_{1}'$ in $\mathbb{P}^{4}$
isomorphically as we have seen above, there exists a member $D_{W}$
of $|-K_{W}|$ containing $\Pi_{1}$ since $\mathcal{O}(-K_{W})=\mathcal{O}(1)$.
Since there are two independent linear forms to define a plane in
$\mathbb{P}^{4}$, we can write $D_{W}=D_{W}'+\Pi_{1}$ with $h^{0}(D'_{W})=2$.
Let $D_{Z}$ be the total transform of $D_{W}$ on $Z$ and $D_{X}:=g_{*}D_{Z}$.
Then, by a similar argument to the check of (1) in Step 1, we see
that $D_{X}\in|-K_{X}|$ and $D_{X}$ contains at least the strict
transform of some irreducible component of $D_{W}'$ and the image
of $F$. This contradicts the assumption that $X$ is prime. \vspace{3pt}

In the following steps, we also consider the blow-up $g'\colon Z'\to X$
at the three $\nicefrac{1}{2}(1,1,1)$-singularities $\mathsf{p}_{1},\mathsf{p}_{2},\mathsf{p}_{3}$.
By a similar argument to Step 1, we can prove that $Z'$ is a weak
Fano 3-fold. Let $h'\colon Z'\to W'$ be the morphism defined by $|-mK_{Z'}|$
with $m\gg$0. By a similar argument to Step 2, we can prove that
$h'$ has no exceptional divisor.

\vspace{3pt}

\noindent\textbf{Step 3.} \textit{We show that $W'$ is a quartic
3-fold in $\mathbb{P}^{4}$ with only terminal singularities.}\vspace{3pt}

By a similar argument to the one in Step 1, we can check that $W'$
satisfies the conditions (1)--(5) of \cite[II, Thm.3.0]{Tak1}, and
applying \cite[II, Thm.3.0]{Tak1} to $W'$, we see that $W'$ is
a quartic 3-fold in $\mathbb{P}^{4}$ or a double cover of a quadric
3-fold in $\mathbb{P}^{4}$ branched along the intersection with a
quartic 3-fold. Assume, for a contradiction, that the latter case
occurs, i.e., $W'$ is the double cover of a quadric 3-fold $Q$.
Similarly to the discussion just before Step 2, we see that the images
on $\mathbb{P}^{4}$ of $g'$-exceptional divisors are planes. Therefore
$Q$ contains a plane, and then $Q$ is singular. By the geometry
of a singular quadric 3-fold, there is a hyperplane section $H_{Q}$
which is the union of two planes $P_{1}$ and $P_{2}$, and we may
choose $P_{1}$ and $P_{2}$ such that their strict transforms are
not $g'$-exceptional divisors. Then, pulling back $H_{Q}$ to $Z'$
and pushing down to $X$, we obtain a reducible anticanonical divisor
of $X$, which contradicts the assumption that $X$ is prime (cf.
the argument in Step 2). Therefore $W'$ is a quartic in $\mathbb{P}^{4}$.
The assertion as for the singularities of $W'$ follows similarly
to Step 2. \vspace{3pt}

\noindent\textbf{Step 4.} \textit{We show that $\Pi_{1}'\cap\Pi_{2}'$
consists of one point.} \vspace{3pt}

First we show that $\Pi_{1}'\cap\Pi_{2}'\not=\emptyset$. To see this,
let $f\colon Y\to X$ be the blow-up at $\mathsf{p}_{1}$ and $E$
the $f$-exceptional divisor. Then, since X is of type No.$\,$5.4,
there is a flop $Y\dashrightarrow Y_{1}$ and a flip $Y_{1}\dashrightarrow Y'$
by \cite[I, Table 5]{Tak1}. Let $E_{Y_{1}}$ be the strict transform
of $E$ on $Y_{1}$. By ibid., there is a flipping curve $m$ on $Y_{1}$
with $-K_{Y_{1}}\cdot m=\nicefrac{1}{2}$ and $E_{Y_{1}}\cdot m>0$
through the $\nicefrac{1}{2}(1,1,1)$-singularity corresponding to
$\mathsf{p}_{2}$. Since the flop $Y\dashrightarrow Y_{1}$ is crepant,
we see that $m$ is the strict transform of a curve $m_{Y}$ on $Y$
with $-K_{Y}\cdot m_{Y}=\nicefrac{1}{2}$ passing through the $\nicefrac{1}{2}(1,1,1)$-singularity
corresponding to $\mathsf{p}_{2}$. Moreover, since a flopping curve
is positive for $E$, we have $E\cdot m_{Y}>0$, or there exists a
connected chain of flopping curves intersecting both $E$ and $m_{Y}$.
Regarding $Z$ as the blow-up of $Y$ at the $\nicefrac{1}{2}(1,1,1)$-singularity
corresponding to $\mathsf{p}_{2}$ (then $E_{1}$ is the strict transform
of the $f$-exceptional divisor $E$), we see that the strict transform
$m_{Z}$ of $m_{Y}$ on $Z$ satisfies $-K_{Z}\cdot m_{Z}=0$. Therefore
$m_{Z}$ is contracted to a point on $\mathbb{P}^{4}$. Moreover $m_{Z}$
intersects both with $E_{1}$ and $E_{2}$, or $m_{Z}$ intersects
with $E_{2}$ and there exists a connected chain of flopping curves
intersecting $E_{1}$ and $m_{Z}$. In any case, we conclude that
$\Pi_{1}'\cap\Pi_{2}'$ contains at least the image of $m_{Z}$. 

Assume, for a contradiction, that $\Pi_{1}'\cap\Pi_{2}'$ does not
consists of one point, namely, $\Pi_{1}'\cap\Pi_{2}'$ is a line.
We may easily check that the closure of the image of $W$ by the projection
$\mathbb{P}(1^{5},2)\dashrightarrow\mathbb{P}^{4}$ is nothing but
$W'$, which is shown to be a quartic 3-fold in $\mathbb{P}^{4}$
in Step 3. Then, since $\Pi_{1}'\cap\Pi_{2}'$ is a line and $h'$
is birational, there exists an $h'$-exceptional divisor, which contradicts
the argument just above Step 3. Thus we have shown that $\Pi_{1}'\cap\Pi_{2}'$
consists of one point. 

\vspace{3pt}

Let $x_{11},x_{21},x_{12},x_{22},x_{13},u_{3}$ be the coordinates
of $\mathbb{P}(1^{5},2)$ such that the weight of $u_{3}$ is 2 and
the weights of the other coordinates are 1 (this deliberate choice
of notation for coordinates will convenient for relating the construction
below with the equation of $\mathscr{H}_{\mA}^{13}$ as in Corollary
\ref{cor:HsharpN}).

\vspace{3pt}

\noindent\textbf{Step 5.} \textit{We prove that we may assume $\Pi_{1}=\{x_{11}=x_{21}=u_{3}=0\}$
and $\Pi_{2}=\{x_{12}=x_{22}=u_{3}=0\}$ by coordinate changes.}\vspace{3pt}

Since $\Pi_{1}'\cap\Pi_{2}'$ consists of one point by Step 4, we
may assume that $\Pi_{1}'=\{x_{11}=x_{21}=0\}$ and $\Pi_{2}'=\{x_{12}=x_{22}=0\}$.
Therefore, since the degrees of $\Pi_{1}$ and $\Pi_{2}$ are one
with respect to $\mathcal{O}(1)$, we have $\Pi_{1}=\{x_{11}=x_{21}=f_{1}(x_{12},x_{22},x_{13},u_{3})=0\}$
and $\Pi_{2}=\{x_{12}=x_{22}=f_{2}(x_{11},x_{21},x_{13},u_{3})=0\}$,
where $f_{1}(x_{12},x_{22},x_{13},u_{3})$ and $f_{2}(x_{11},x_{21},x_{13},u_{3})$
are of weight $2$. Since $\Pi_{1},\Pi_{2}$ do not contain the $\nicefrac{1}{2}(1^{5})$-singularity
of $\mathbb{P}(1^{5},2)$ by the argument just before Step 2, we see
that the polynomials $f_{1},f_{2}$ contains $u_{3}$. By a coordinate
change, we may assume that $f_{1}=u_{3}$ and $f_{2}=u_{3}+g_{2}(x_{11},x_{21},x_{13})$.
Note that $g_{2}$ does not contain $x_{13}^{2}$ since otherwise
$\Pi_{1}\cap\Pi_{2}$ must be empty. Now we set $u_{3}'=u_{3}+g_{2}(x_{11},x_{21},x_{13}).$
Then $\Pi_{2}=\{x_{12}=x_{22}=u_{3}'=0\}$ and $\Pi_{1}=\{x_{11}=x_{21}=u_{3}'-g_{2}(x_{11},x_{21},x_{13})=0\}=\{x_{11}=x_{21}=u_{3}'=0\}$.
Denoting $u_{3}'$ by $u_{3}$, we obtain the desired expression of
$\Pi_{1}$ and $\Pi_{2}$.

\vspace{3pt}

\noindent\textbf{Step 6. }\textit{We extend $W$ to a $12$-dimensional
weighted compelete intersection $\overline{\mathscr{H}}$ of two weight
3 hypersurfaces in $\mathbb{P}(1^{14},2)$.}\vspace{3pt}

We may write equations of $W$ as follows since it contains $\Pi_{1}$
and $\Pi_{2}${\small :
\begin{equation}
\begin{cases}
u_{3}l_{13}-\begin{pmatrix}x_{12} & x_{22}\end{pmatrix}\begin{pmatrix}l_{221} & -l_{121}\\
-l_{211} & l_{111}
\end{pmatrix}\left(\begin{array}{c}
x_{11}\\
x_{21}
\end{array}\right)=0,\\
u_{3}l_{23}-\begin{pmatrix}x_{12} & x_{22}\end{pmatrix}\begin{pmatrix}l_{222} & -l_{122}\\
-l_{212} & l_{112}
\end{pmatrix}\left(\begin{array}{c}
x_{11}\\
x_{21}
\end{array}\right)=0,
\end{cases}\label{eq:eqW}
\end{equation}
}where $l_{ijk}$'s, $l_{st}$ are linear forms. 

We show that $l_{13}$ or $l_{23}$ contains $x_{13}$. Assume, for
a contradiction, that neither $l_{13}$ nor $l_{23}$ contains $x_{13}$.
Then both equations of (\ref{eq:eqW}) are singular at the $x_{13}$-point.
Therefore $W$ has a non-hypersurface, hence non-terminal singularity
at the $x_{13}$-point. This contradicts Step 2.

By changing the order of the two equations of (\ref{eq:eqW}), we
may assume that $l_{13}$ contains $x_{13}$. Moreover, changing the
coordinates, we may also assume that $l_{13}=x_{13}.$ Now (\ref{eq:eqW})
can be written as follows: {\small
\begin{equation}
u_{3}\left(\begin{array}{c}
x_{13}\\
l_{23}
\end{array}\right)=\begin{pmatrix}-\begin{vmatrix}l_{211} & x_{12}\\
l_{221} & x_{22}
\end{vmatrix} & \begin{vmatrix}l_{111} & x_{12}\\
l_{121} & x_{22}
\end{vmatrix}\\
-\begin{vmatrix}l_{212} & x_{12}\\
l_{222} & x_{22}
\end{vmatrix} & \begin{vmatrix}l_{112} & x_{12}\\
l_{122} & x_{22}
\end{vmatrix}
\end{pmatrix}\bm{x}_{1},\label{eq:eqWvect}
\end{equation}
}where we set $\bm{x}_{1}:=\left(\begin{array}{c}
x_{11}\\
x_{21}
\end{array}\right)$. We also set $\bm{x}_{2}:=\left(\begin{array}{c}
x_{12}\\
x_{22}
\end{array}\right),$ $\bm{x}_{3}:=\left(\begin{array}{c}
x_{13}\\
x_{23}
\end{array}\right)$, and denote by $\mathbb{P}(\bm{x}_{1},\bm{x}_{2},\bm{x}_{3},p_{ijk},u_{3})$
($i,j,k=1,2$) the weighted projective space with the entries of $\bm{x}_{i}\,(i=1,2,3)$
and $p_{ijk}$ being weight $1$ coordinates, and $u_{3}$ being a
weight $2$ coordinate. Then $W$ is equal to 
\[
\left\{ u_{3}\bm{x}_{3}=D^{(2)}\bm{x}_{1}\right\} \cap\{p_{ijk}=l_{ijk},x_{23}=l_{23}\}\subset\mathbb{P}(\bm{x}_{1},\bm{x}_{2},\bm{x}_{3},p_{ijk},u_{3}),
\]
where $D^{(2)}$ is defined as in Corollary \ref{cor:HsharpN}.\vspace{5pt}

Finally, we set 
\[
\overline{\mathscr{H}}:=\left\{ u_{3}\bm{x}_{3}=D^{(2)}\bm{x}_{1}\right\} \subset\mathbb{P}(\bm{x}_{1},\bm{x}_{2},\bm{x}_{3},p_{ijk},u_{3}).
\]

\vspace{3pt}

Let $\mathscr{H}_{\mathbb{P}}^{12}$ be the weighted projectivization
of $\mathscr{H}_{\mathbb{A}}^{13}$ such that the weights of $u_{1},u_{2},u_{3}$
is 2 and the weights of the other coordinates are 1. We set 
\[
X':=\mathscr{H}_{\mathbb{P}}^{12}\cap\{p_{ijk}=l_{ijk},x_{23}=l_{23}\}.
\]
Note that the projection $\mathbb{P}(\bm{x}_{1},\bm{x}_{2},\bm{x}_{3},p_{ijk},u_{1},u_{2},u_{3})\dashrightarrow\mathbb{P}(\bm{x}_{1},\bm{x}_{2},\bm{x}_{3},p_{ijk},u_{3})$
induces rational maps $\mathscr{H}_{\mathbb{P}}^{12}\dashrightarrow\overline{\mathscr{H}}$
and $X'\dashrightarrow W$.

\vspace{3pt}

\noindent\textbf{Step 7. }\textit{We show that $X$ is isomorphic
to $X'$, which finishes the proof of Theorem \ref{thm:embthm} (2).}\vspace{3pt}

By \cite[Lem.5.5]{Tak2}, it suffices to show that 

\noindent(i) $X'$ is a normal 3-fold,

\noindent(ii) $-K_{X'}$ is $\mathbb{Q}$-Cartier and ample, and 

\noindent(iii) the rational map $X\dashrightarrow X'$ which is the
composite of $X\dashrightarrow W$ and the inverse of $X'\dashrightarrow W$
is isomorphic in codimension 1. 

Note that (ii) follows from (i) and the proof of Theorem \ref{thm:embthm}
(1). Thus it suffices to show (i) and (iii). For this, we give preparatory
argument.

It is easy to check that the rational map $\mathscr{H}_{\mathbb{P}}^{12}\dashrightarrow\overline{\mathscr{H}}$
is defined outside the $u_{1}$- and the $u_{2}$-points. Moreover,
using the equations (\ref{eq:HypMat}), we see that, outside $\{\bm{x}_{1}=\bm{o},u_{3}=0\}\cup\{\bm{x}_{2}=\bm{o},u_{3}=0\}$,
the coordinates $u_{1}$ and $u_{2}$ are recovered by this rational
map with the equations $\bm{G}_{1}=\bm{o}$, $\bm{G}_{2}=\bm{o}$,
$\bm{G}_{3}=\bm{o}$, $\mathsf{G}_{4}=0$, $\mathsf{G}_{5}=0$. Therefore
it holds that 
\begin{align*}
\mathscr{H}_{\mathbb{P}}^{12}\setminus\Gamma\simeq\overline{\mathscr{H}}\setminus\left((\{\bm{x}_{1}=\bm{o},u_{3}=0\}\cup\{\bm{x}_{2}=\bm{o},u_{3}=0\})\right),
\end{align*}
where we set 
\[
\Gamma:=(\{\bm{x}_{1}=\bm{o},u_{3}=0\}\cup\{\bm{x}_{2}=\bm{o},u_{3}=0\})\cap\mathscr{H}_{\mathbb{P}}^{12}.
\]
This implies that 
\begin{equation}
X'\setminus\left(\Gamma\cap X'\right)\simeq W\setminus(\Pi_{1}\cup\Pi_{2}).\label{eq:X'=00003DW}
\end{equation}
 By direct computations, we can verify that

{\small
\begin{align*}
\{\bm{x}_{1}=\bm{o},u_{3}=0\}\cap\mathscr{H}_{\mathbb{P}}^{12}= & \{\bm{x}_{1}=\bm{o},u_{2}=u_{3}=0,D^{(3)}\bm{x}_{2}=\bm{o},|D^{(2)}|=|D^{(3)}|=0\}\\
\cup & \{\bm{x}_{1}=\bm{x}_{2}=\bm{o},u_{3}=0,u_{1}u_{2}+|D^{(3)}|=0\},
\end{align*}
} and it is mapped by $\mathscr{H}_{\mP}^{12}\dashrightarrow\overline{\mathscr{H}}$
to the locus
\[
S_{1}:=\{\bm{x}_{1}=\bm{o},u_{3}=0,D^{(3)}\bm{x}_{2}=\bm{o},|D^{(2)}|=0\}
\]
with $1$-dimensional fibers. We note that $D^{(1)}\bm{x}_{3}=-D^{(3)\dagger}\bm{x}_{1}$,
where $D^{(3)\dagger}$ is the adjoint of $D^{(3)}.$ Then, setting
\[
S_{2}:=\{\bm{x}_{2}=\bm{o},u_{3}=0,=\bm{o},D^{(3)\dagger}\bm{x}_{1}=\bm{o},|D^{(1)}|=0\},
\]
we can similarly check that $\{\bm{x}_{2}=\bm{o},u_{3}=0\}\cap\mathscr{H}_{\mathbb{P}}^{12}$
is mapped by $\mathscr{H}_{\mP}^{12}\dashrightarrow\overline{\mathscr{H}}$
to $S_{2}$ with $1$-dimensional fibers. Therefore, $\Gamma\cap X'$
is mapped to the image $S_{1}|_{W}\cup S_{2}|_{W}$ with 1-dimensional
fibers. On the other hand, using (\ref{eq:eqWvect}), we can directly
check that ${\rm Sing}\,W=\{\text{the}\,u_{3}\text{-point}\}\cup(S_{1}|_{W}\cup S_{2}|_{W})$.
Therefore, by Step 2, $\dim(S_{1}|_{W}\cup S_{2}|_{W})=0$, and then
$\dim\Gamma\cap X'=1.$

Now we show (i). By the above discussion, $X'$ is 3-dimentional since
$X'$ is the union of a 1-dimensional closed subset $\Gamma\cap X'$
and a $3$-dimensional open subset isomorphic to $W\setminus(\Pi_{1}\cup\Pi_{2})$.
Since then the affine cone $X_{\mathbb{A}}'$ of $X'$ is of codimension
$9$ in $\mathscr{H}_{\mathbb{A}}^{13}$ and the ideal of $X_{\mathbb{A}}'$
in $\mathscr{H}_{\mathbb{A}}^{13}$ is generated by 9 elements, $X'_{\mA}$
is Gorenstein and is of pure dimension 4 by unmixedness of the Gorenstein
variety $\mathscr{H}_{\mathbb{A}}^{13}$. Thus $X'$ is an irreducible
3-fold and satisfies $S_{2}$-condition. Since $W\setminus(\Pi_{1}\cup\Pi_{2})$
has only one $\nicefrac{1}{2}(1,1,1)$-singularity by Step 2, $X'$
satisfies the $R_{1}$-condition. Consequently, $X'$ is a normal
3-fold as desired.

Finally we show (iii). Since $W$ has only terminal singularities
by Step 2, the morphism $h$ is a small contraction and the images
of $h$-exceptional curves are contained in $\Pi_{1}\cup\Pi_{2}$.
Thus $W\setminus(\Pi_{1}\cup\Pi_{2})$ is isomorphic to an open subset
of $X$ obtained by removing a closed subset of codimension 2. Therefore,
$X$ and $X'$ are isomorphic in codimension $1$ as desired by (\ref{eq:X'=00003DW})
and the fact that $\dim\Gamma\cap X'=1$.
\end{proof}

\section{\textbf{Relatives of $\mathscr{H}_{\mA}^{13}$--affine varieties
$\mathscr{M}_{\mathbb{A}}^{8}$ and $\mathscr{S}_{\mA}^{6}$--\label{sec:Relatives-of-=002013AffineMS}}}

We define two special subvarieties $\mathscr{M}_{\mathbb{A}}^{8}$
and $\mathscr{S}_{\mA}^{6}$ of $\mathscr{H}_{\mA}^{13}$ with good
group actions. In Section \ref{sec:The--cluster-variety}, it turns
out that their several weighted projectivizations are key varieties
for prime $\mQ$-Fano 3-folds.

\subsection{Definition of $\mathscr{M}_{\mathbb{A}}^{8}$ }
\begin{defn}
\label{def:MS}

We define $\mathscr{M}_{\mathbb{A}}^{8}$ to be the affine scheme
in the $12$-dimensional affine space $\mathbb{A}_{\mathscr{M}}$
with the coordinates implemented in 
\[
\mathsf{X}_{1}:=\left(\begin{array}{cc}
x_{1} & x_{2}\\
x_{3} & -x_{1}
\end{array}\right),\mathsf{X_{2}}:=\left(\begin{array}{cc}
x_{4} & x_{5}\\
x_{6} & x_{7}
\end{array}\right),W:=\left(\begin{array}{cc}
w_{1} & w_{2}\\
w_{3} & -w_{1}
\end{array}\right),y,z,
\]
defined by the following equations:
\begin{equation}
yz=\det\mathsf{X}_{2},\ y\mathsf{W}=\mathsf{X}_{2}^{\dagger}\mathsf{X}_{1}\mathsf{X}_{2},\ \mathsf{X}_{2}\mathsf{W}=z\mathsf{X}_{1}\mathsf{X}_{2},\ \det\mathsf{W}=z^{2}\det\mathsf{X}_{1},\label{eq:M8}
\end{equation}
where $\mathsf{X}_{2}^{\dagger}$ is the adjoint matrix of $\mathsf{X}_{2}$. 
\end{defn}

We can show the following proposition by a direct calculation:
\begin{prop}
The affine scheme $\mathscr{H}_{\mA}^{13}\cap\{p_{221}=-1,p_{222}=0,p_{211}=0,p_{212}=1,p_{121}=p_{112}\}$
coincides with $\mathscr{M}_{\mA}^{8}$ by changing the notation of
the coordinates as follows:

\begin{align*}
\left(\begin{array}{cc}
x_{1} & x_{2}\\
x_{3} & -x_{1}
\end{array}\right) & =\left(\begin{array}{cc}
p_{112} & -p_{111}\\
p_{122} & -p_{112}
\end{array}\right),\,\left(\begin{array}{cc}
x_{4} & x_{5}\\
x_{6} & x_{7}
\end{array}\right)=\left(\begin{array}{cc}
x_{13} & x_{12}\\
x_{23} & x_{22}
\end{array}\right),\\
\left(\begin{array}{cc}
w_{1} & w_{2}\\
w_{3} & -w_{1}
\end{array}\right) & =\left(\begin{array}{cc}
x_{11} & -u_{3}\\
-u_{2} & -x_{11}
\end{array}\right),\,y=-u_{1},\,z=x_{21},
\end{align*}
where the l.h.s.~are the coordinates of $\mathscr{M}_{\mA}^{8}$,
and the r.h.s.~are those of $\mathscr{H}_{\mA}^{13}$. 
\end{prop}

\subsection{$({\rm GL}_{2})^{2}$-action on $\mathscr{M}_{\mathbb{A}}^{8}$\label{subsec:Group-action-H}}

By the equations (\ref{eq:M8}) of $\mathscr{M}_{\mathbb{A}}^{8}$,
we can directly check the following:
\begin{prop}
The group $({\rm GL}_{2})^{2}$ acts on $\mathscr{M}_{\mathbb{A}}^{8}$
by the following rules: for any $(g_{1},g_{2})\in({\rm GL}_{2})^{2}$,
we set

\begin{align*}
 & \mathsf{X}_{1}\mapsto g_{1}\mathsf{X}_{1}g_{1}^{\dagger},\text{\ensuremath{\mathsf{X}_{2}\mapsto g_{1}\mathsf{X}_{2}g_{2}^{\dagger}}},\mathsf{W}\mapsto g_{2}\mathsf{W}g_{2}^{\dagger},\\
 & y\mapsto(\det g_{1})^{2}y,z\mapsto(\det g_{2})(\det g_{1})^{-1}z
\end{align*}

where $g_{1}^{\dagger},g_{2}^{\dagger}$ are the adjoint matrices
of $g_{1},g_{2}$ respectively.
\end{prop}

\subsection{$\mP^{2}\times\mP^{2}$-fibration related to $\mathscr{M}_{\mA}^{8}$\protect 
}Let $\mA_{\mathsf{X}_{1}}$ be the affine 3-space with the coordinates
$x_{1},x_{2},x_{3}$. As in the case of $\text{\ensuremath{\mathscr{H}_{\mathbb{A}}^{13}}}$,
considering the coordinates in the equations of $\text{\text{\ensuremath{\mathscr{M}_{\mathbb{A}}^{8}}}}$
except those of $\mA_{\mathsf{X}_{1}}$ as projective coordinates,
we obtain a quasi-projective variety with the same equations as $\text{\ensuremath{\mathscr{M}_{\mathbb{A}}^{8}}}$,
which we denote by $\widehat{\mathscr{M}}$. We also denote by $\rho_{\mathscr{M}}\colon\widehat{\mathscr{M}}\to\mA_{\mathsf{X}_{1}}$
the natural projection. Using the group action on $\mathscr{M}_{\mA}^{8}$,
we can show the following in a similar way to Proposition \ref{prop:Let--beP2P2fib}:
\begin{prop}
\label{prop:Let--beP2P2fibM8} The $\rho_{\mathscr{M}}$-fibers over
points of the complement of $\{\det\mathsf{X}_{1}=0\}$ are isomorphic
to $\mathbb{P}^{2}\times\mathbb{P}^{2}$, and the $\rho_{\mathscr{M}}$-fibers
over points of are $\{\det\mathsf{X}_{1}=0\}\setminus\{o\}$ isomorphic
to $\mathbb{P}^{2,2}$. Any $\rho_{\mathscr{M}}$-fibers are $4$-dimensional.
In particular, $\mathscr{M}_{\mA}^{8}$ is irreducible and $8$-dimensional.
\end{prop}

\subsection{Definition of $\mathscr{S}_{\mA}^{6}$\protect 
}\begin{defn}
We define $\mathscr{S}_{\mathbb{A}}^{6}$ to be the affine scheme
in the $10$-dimensional affine space $\mathbb{A}_{\mathscr{S}}$
with the coordinates implemented in 

\[
\mathsf{S}:=\left(\begin{array}{ccc}
s_{11} & s_{12} & s_{13}\\
s_{12} & s_{22} & s_{23}\\
s_{13} & s_{23} & s_{33}
\end{array}\right),\,\bm{\sigma}:=\left(\begin{array}{c}
\sigma_{1}\\
\sigma_{2}\\
\sigma_{3}
\end{array}\right),t
\]
defined by the following equations: 

\begin{equation}
\mathsf{S}\bm{\bm{\sigma}}=\bm{0},\mathsf{S}^{\dagger}=t\bm{\bm{\sigma}}(\empty^{t}\!\bm{\bm{\sigma}}),\label{eq:S}
\end{equation}
where $\mathsf{S}^{\dagger}$ is the adjoint matrix of $\mathsf{S}$. 
\end{defn}

\begin{rem}
The equations (\ref{eq:S}) appears in \cite[Ex.5.4]{GN}.
\end{rem}

We can show the following proposition by a direct calculation:
\begin{prop}
The affine scheme $\mathscr{H}_{\mA}^{13}\cap\{p_{221}=p_{212}=p_{122}=1,p_{222}=p_{211}=p_{121}=p_{112}=0\}$
coincides with $\mathscr{S}_{\mA}^{6}$ by changing the notation of
the coordinates as follows:

\begin{align*}
\left(\begin{array}{ccc}
s_{11} & s_{12} & s_{13}\\
 & s_{22} & s_{23}\\
 &  & s_{33}
\end{array}\right) & =\left(\begin{array}{ccc}
u_{1} & x_{13} & x_{12}\\
 & u_{2} & x_{11}\\
 &  & u_{3}
\end{array}\right),\,\left(\begin{array}{c}
\sigma_{1}\\
\sigma_{2}\\
\sigma_{3}
\end{array}\right)=\left(\begin{array}{c}
x_{21}\\
x_{22}\\
x_{23}
\end{array}\right),\,t=-p_{111},\\
\end{align*}
where the l.h.s.~are the coordinates of $\mathscr{S}_{\mA}^{6}$,
and the r.h.s.~are those of $\mathscr{H}_{\mA}^{13}$.
\end{prop}

\subsection{${\rm GL}{}_{3}$-action on $\mathscr{S}_{\mA}^{6}$}

By the equations (\ref{eq:S}) of $\mathscr{S}_{\mathbb{A}}^{6}$,
we can directly check the following:
\begin{prop}
The group ${\rm GL}{}_{3}$ acts on $\mathscr{S}_{\mathbb{A}}^{6}$
by the following rules: for any $g\in{\rm GL}_{3}$, we set

\begin{align*}
 & \mathsf{S}\mapsto g\mathsf{S(\empty^{t}}g),\bm{\bm{\sigma}}\mapsto\mathsf{(\empty^{t}}g)^{\dagger}\bm{\bm{\sigma}},t\mapsto t
\end{align*}

where $g^{\dagger}$ is the adjoint matrix of $g$ respectively.
\end{prop}

\subsection{$\mP^{2}\times\mP^{2}$-fibration related to $\mathscr{S}_{\mA}^{6}$}

Let $\mA_{t}$ be the affine line with the coordinates $t$. As in
the case of $\text{\ensuremath{\mathscr{H}_{\mathbb{A}}^{13}}}$,
considering the coordinates in the equations of $\text{\text{\ensuremath{\mathscr{S}_{\mathbb{A}}^{6}}}}$
except $t$ as projective coordinates, we obtain a quasi-projective
variety with the same equations as $\text{\ensuremath{\mathscr{S}_{\mathbb{A}}^{6}}}$,
which we denote by $\widehat{\mathscr{S}}$. We also denote by $\rho_{\mathscr{S}}\colon\widehat{\mathscr{S}}\to\mA_{t}$
the natural projection. Using the group action on $\mathscr{S}_{\mA}^{6}$,
we can show the following in a similar way to Proposition \ref{prop:Let--beP2P2fib}:
\begin{prop}
\label{prop:Let--beP2P2fibS6} The $\rho_{\mathscr{M}}$-fibers except
the one over the origin are isomorphic to $\mathbb{P}^{2}\times\mathbb{P}^{2}$,
and the $\rho_{\mathscr{M}}$-fiber over the origin is isomorphic
to $\mathbb{P}^{2,2}$. In particular, $\mathscr{S}_{\mA}^{6}$ is
irreducible and $6$-dimensional.
\end{prop}

\vspace{0.5cm}

As for $\mathscr{M}_{\mA}^{8}$ and $\mathscr{S}_{\mA}^{6}$, we obtain
similar properties to those of $\mathscr{H}_{\mA}^{13}$ as in Subsections
\ref{subsec:Descriptrion-of-charts} -- \ref{subsec:Singularties-ofH}.
We omit deriving these since the method is similar to that in the
case of $\mathscr{H}_{\mA}^{13}$. Actually, we need those properties
of $\mathscr{M}_{\mA}^{8}$ and $\mathscr{S}_{\mA}^{6}$ if we directly
construct prime $\mQ$-Fano 3-folds from their weighted projectivizations.
We can, however, take a shortcut using the result of \cite{CD1} in
Section \ref{sec:The--cluster-variety}.

\section{\textbf{The $C_{2}$-cluster variety as a specialization of $\mathscr{H}_{\mA}^{13}$\label{sec:The--cluster-variety}}}

\subsection{The cluster variety $X_{\mathsf{C}_{2}}$ of type $C_{2}$}

In the paper \cite{CD1}, inspired by mirror symmetry of log Calabi-Yau
surfaces after Gross, Hacking and Keel \cite{GHK}, Coughlan and Ducat
define\textit{ }the cluster variety of type $C_{2}$ as follows and
produces, as weighted complete intersections of their weighted projectivizations,
several examples of prime $\mathbb{Q}$-Fano 3-folds $X$ anticanonically
embedded of codimension $4$.
\begin{defn}
\textit{\label{def:C2}The cluster variety of type $C_{2}$, denoted
by $X_{\mathsf{C}_{2}}$}, is the affine variety in the $13$-affine
space $\mathbb{A}_{\mathsf{C}_{2}}$ with coordinates implemented
as in 

\[
\left(\begin{array}{cccccc}
\theta_{12} & \theta_{23} & \theta_{31} & \theta_{1} & \theta_{2} & \theta_{3}\\
A_{12} & A_{23} & A_{31} & A_{1} & A_{2} & A_{3}
\end{array}\right),\lambda,
\]
defined by the following $9$ equations:
\begin{equation}
\begin{cases}
\theta_{i}\theta_{j}=A_{ij}\theta_{ij}+A_{jk}A_{k}A_{ki},\\
\theta_{ki}\theta_{ij}=A_{i}\theta_{i}^{2}+\lambda A_{jk}\theta_{i}+A_{j}A_{jk}^{2}A_{k},\\
\theta_{i}\theta_{jk}=A_{ij}A_{j}\theta_{j}+\lambda A_{ki}A_{ij}+A_{k}A_{ki}\theta_{k},
\end{cases}\label{eq:ClusterEq}
\end{equation}
where $(i,j,k)=(1,2,3),(2,3,1),(3,1,2)$.
\end{defn}

\subsection{$X_{\mathsf{C}_{2}}$ as a subvariety of $\mathscr{H}_{\mathbb{A}}^{13}$ }

By direct computations, we may check the following:
\begin{prop}
\label{prop:C2}The affine variety $\mathscr{H}_{\mathbb{A}}^{13}\cap\{p_{121}=p_{112}=p_{211}=0,\,p_{111}=1\}$
coincides with $X_{\mathsf{C}_{2}}$ by changing the notation of coordinates
as follows:
\begin{align*}
u_{1} & =\theta_{12},\,u_{2}=\theta_{23},\,u_{3}=\theta_{31},\\
\bm{x}_{1} & =\left(\begin{array}{c}
A_{12}\\
\theta_{3}
\end{array}\right),\,\bm{x}_{2}=\left(\begin{array}{c}
A_{23}\\
\theta_{1}
\end{array}\right),\,\bm{x}_{3}=\left(\begin{array}{c}
A_{31}\\
\theta_{2}
\end{array}\right),\\
p_{212} & =-A_{1},\,p_{221}=-A_{2},\,p_{122}=-A_{3},\,p_{222}=\lambda,
\end{align*}
where the l.h.s. are the coordinates of $\mathscr{H}_{\mathbb{A}}^{13}$
and the r.h.s. are those of $X_{\mathsf{C}_{2}}$.
\end{prop}

By this proposition, we often consider $X_{\mathsf{C}_{2}}$ as a
subvariety of $\mathscr{H}_{\mA}^{13}$.

\subsection{Subvarieties $\mathscr{H}_{\mathbb{A}}^{11}$ and $\mathscr{H}_{\mathbb{A}}^{12}$
of $\mathscr{H}_{\mathbb{A}}^{13}$ related to $X_{\mathsf{C}_{2}}$}
\begin{defn}
We set
\begin{align*}
\mathscr{H}_{\mathbb{A}}^{12}:=\mathscr{H}_{\mathbb{A}}^{13}\cap\{p_{111}=1\},\ \mathscr{H}_{\mathbb{A}}^{11}:=\mathscr{H}_{\mathbb{A}}^{12}\cap\{p_{121}=1\}.
\end{align*}
Actually, $\mathscr{H}_{\mathbb{A}}^{12}\cap\{p_{121}=1\}$, $\mathscr{H}_{\mathbb{A}}^{12}\cap\{p_{112}=1\}$,
$\mathscr{H}_{\mathbb{A}}^{12}\cap\{p_{221}=1\}$ are transformed
each other by permuting the coordinates $p_{121},p_{112},p_{211}$
by the $\mathfrak{S}_{3}$-action. Therefore, abusing notation, we
also denote $\mathscr{H}_{\mathbb{A}}^{12}\cap\{p_{112}=1\}$, $\mathscr{H}_{\mathbb{A}}^{12}\cap\{p_{221}=1\}$
by $\mathscr{H}_{\mathbb{A}}^{11}$. 
\end{defn}

Note that $X_{\mathsf{C}_{2}}$ is a subvariety of $\mathscr{H}_{\mA}^{12}$.
Moreover, we can regard $X_{\mathsf{C}_{2}}$ as a specialization
of a subvariety of $\mathscr{H}_{\mathbb{A}}^{11}$. Indeed, we can
verify $\mathscr{H}_{\mathbb{A}}^{12}\cap\{p_{121}=1\}$ is isomorphic
to $\mathscr{H}_{\mathbb{A}}^{12}\cap\{p_{121}=\alpha\}\,(\alpha\not=0)$.
Then $X_{\mathsf{C}_{2}}$ is a specialization of $\mathscr{H}_{\mathbb{A}}^{12}\cap\{p_{121}=\alpha,p_{112}=p_{211}=0\}$
by $\alpha\to0$.

The reason why we consider $\mathscr{H}_{\mathbb{A}}^{11}$ and $\mathscr{H}_{\mathbb{A}}^{12}$
is that, using these, we can generalize some examples of prime $\mQ$-Fano
$3$-folds obtained from $X_{\mathsf{C}_{2}}$; let us consider the
following types of prime $\mQ$-Fano $3$-folds of \cite{GRDB}:

\begin{enumerate}[(A)]

\item No.$\,$1397, 1405, 2427, 2511, 5052, 5305, 5410, 5516, 5963,
11222 (10 classes).

\item No.$\,$5530, 5970, 11437, 11455, 16339 (5 classes).

\end{enumerate}

\noindent We can show directly that positive weights of coordinates
are associated to $\mathscr{H}_{\mA}^{12}$ (resp. $\mathscr{H}_{\mA}^{11}$)
by 
\begin{itemize}
\item the data of \cite{CD2} for the classes in the group (A) (resp. (B))
,
\item the transformations of the coordinates as in Proposition \ref{prop:C2},
and
\item the set-up
\end{itemize}
{\small
\begin{equation}
\begin{cases}
w(p_{121}):=w(p_{122})+w(p_{221})-w(p_{222}),\\
w(p_{112}):=w(p_{122})+w(p_{212})-w(p_{222}),\\
w(p_{211}):=w(p_{212})+w(p_{221})-w(p_{222}).
\end{cases}\label{eq:p121}
\end{equation}
}(we note that one of $w(p_{121})$, $w(p_{112})$ or $w(p_{211})$
is zero for the group (B), so (\ref{eq:p121}) gives suitable positive
weights of coordinates of $\mathscr{H}_{\mA}^{11}$) . We may also
directly verify that all the equations of $\mathscr{H}_{\mathbb{P}}^{11}$
and $\mathscr{H}_{\mathbb{P}}^{10}$ in each case are homogeneous
with these weights of coordinates. Therefore, by \cite{CD1}, we can
construct examples of prime $\mathbb{Q}$-Fano $3$-folds $X$ for
the classes in the group (A) (resp. (B)) as weighted complete intersection
of $\mathscr{H}_{\mathbb{P}}^{11}$(resp. $\mathscr{H}_{\mathbb{P}}^{10}$).
Actually, more general examples are obtained since $X_{\mathsf{C}_{2}}$
is a specialization of a subvariety of $\mathscr{H}_{\mathbb{A}}^{11}$
or $\mathscr{H}_{\mathbb{A}}^{12}$ as we see above.

We describe an example of a prime $\mathbb{Q}$-Fano $3$-fold obtained
from $\mathscr{H}_{\mathbb{\mP}}^{11}$.
\begin{example}
Let $X$ be a prime $\mathbb{Q}$-Fano $3$-fold $X\subset\mathbb{P}(1^{2},3,4^{2},5^{3})$
of type No.$\,$5052 with three $\nicefrac{1}{5}\,(1,1,4)$-singularities.
Based on the data of weights of coordinates of $X_{\mathsf{C}_{2}}$
as in \cite{CD2} and (\ref{eq:p121}), we put the following weights
of coordinates for $\mathscr{H}_{\mathbb{A}}^{12}$ :

\begin{align*}
w(\bm{x}_{1}) & =w(\bm{x}_{2})=w(\bm{x}_{3})=\begin{pmatrix}3\\
4
\end{pmatrix},\\
w(u_{1}) & =w(u_{2})=w(u_{3})=5,\\
w(p_{112}) & =w(p_{121})=w(p_{211})=1,\\
w(p_{122}) & =w(p_{212})=w(p_{221})=2,\\
w(p_{222}) & =3.
\end{align*}
Then we obtain a weighted projectivization $\mathscr{H}_{\mathbb{P}}^{11}$
of $\mathscr{H}_{\mathbb{A}}^{12}$ in $\mathbb{P}^{15}(1^{3},2^{3},3^{4},4^{3},5^{3})$.
A general 
\[
X:=\mathscr{H}_{\mathbb{P}}^{11}\cap(1)\cap(2)^{3}\cap(3)^{3}\cap(4)
\]
 is an example of a prime $\mathbb{Q}$-Fano $3$-fold $X$ of type
No.$\,$5052.

By \cite{CD1}, an example of such an $X$ is obtained from $X_{\mathsf{C}_{2}}$
but two weight 1 coordinates are missing for $X_{\mathsf{C}_{2}}$,
thus it is necessary to take the cone over $X_{\mathsf{C}_{2}}$ with
two weight 1 coordinates added. Therefore the example obtained here
is more general than one in \cite{CD1}.
\end{example}

\subsection{$\mathscr{M}_{\mathbb{A}}^{8}$ and $\mathscr{S}_{\mathbb{A}}^{6}$
as subvarieties of $X_{\mathsf{C}_{2}}$\label{subsec:MS}}

In the following proposition, we use the notation in Section \ref{sec:Relatives-of-=002013AffineMS}
freely:
\begin{prop}
\label{prop:MSC2}\begin{enumerate}[$(1)$] 

\item By the isomorphism $\mathbb{A}_{\mathsf{C}_{2}}\cap\{A_{3}=1\}\to\mathbb{A}_{\mathscr{M}}$
defined by 

{\small
\begin{align}
{\sf{X}}_{1} & =\begin{pmatrix}\nicefrac{1}{2}\lambda & A_{2}\\
-A_{1} & -\nicefrac{1}{2}\lambda
\end{pmatrix},\,{\sf{X}}_{2}=\begin{pmatrix}\theta_{1} & A_{31}\\
A_{23} & \theta_{2}
\end{pmatrix},\nonumber \\
{\sf{W}} & =\begin{pmatrix}\theta_{3}+\nicefrac{1}{2}\lambda A_{12} & \theta_{23}\\
-\theta_{31} & -\theta_{3}-\nicefrac{1}{2}\lambda A_{12}
\end{pmatrix},\label{eq:MC2}\\
y & =\theta_{12},\,z=A_{12},\nonumber 
\end{align}
}the affine variety $\mathscr{M}_{\mathbb{A}}^{8}$ is the isomorphic
image of $X_{\mathsf{C}_{2}}^{I}:=X_{\mathsf{C}_{2}}\cap\{A_{3}=1\}$.

If we put the positive weights to the coordinates of $\mathbb{A}_{\mathsf{C}_{2}}\cap\{A_{3}=1\}$
such that all the equations of $X_{\mathsf{C}_{2}}^{I}$ are homogeneous,
then we can define the positive weights to the coordinates of $\mathbb{A}_{\mathscr{M}}$
by the equalities $(\ref{eq:MC2})$ and the isomorphism $X_{\mathsf{C}_{2}}^{I}\to\mathscr{M}_{\mathbb{A}}^{8}$
becomes homogeneous.

\item Let $\mathbb{A}_{\mathscr{S}}'$ be the affine space obtained
from $\mathbb{A}_{\mathscr{S}}$ by adding one coordinate $\sigma_{0}$.
By the isomorphism $\mathbb{A}_{\mathsf{C}_{2}}\cap\{A_{3}=1,A_{1}=-1\}\to\mathbb{A}_{\mathscr{S}}'$
defined by {\small
\begin{align}
\mathsf{S} & =\left(\begin{array}{ccc}
\theta_{12} & -\theta_{1}+\nicefrac{1}{2}\lambda A_{23} & A_{31}-\nicefrac{1}{2}\lambda\theta_{2}\\
 & -\theta_{31} & \theta_{3}+\nicefrac{1}{2}\lambda A_{12}\\
 &  & -\theta_{23}
\end{array}\right),\bm{\bm{\sigma}}=\left(\begin{array}{c}
A_{12}\\
\theta_{2}\\
A_{23}
\end{array}\right),\label{eq:SC2}\\
 & t=-\nicefrac{1}{4}\lambda^{2}-A_{2},\sigma_{0}=\lambda,\nonumber 
\end{align}
}the cone $\mathscr{S}_{\mathbb{A}}^{7}$ over $\mathscr{S}_{\mathbb{A}}^{6}$
with the coordinate $\sigma_{0}$ added is the isomorphic image of
$X_{\mathsf{C}_{2}}^{I\!I}:=X_{\mathsf{C}_{2}}\cap\{A_{3}=1,A_{1}=-1\}$. 

If we put the positive weights to the coordinates of $\mathbb{A}_{\mathsf{C}_{2}}\cap\{A_{3}=1,A_{1}=-1\}$
such that all the equations of $X_{\mathsf{C}_{2}}^{I\!I}$ are homogeneous,
then we can define the positive weights to the coordinates of $\mathbb{A}_{\mathscr{S}}$
by the equality $(\ref{eq:SC2})$ and the isomorphism $X_{\mathsf{C}_{2}}^{I\!I}\to\mathscr{S}_{\mathbb{A}}^{7}$
becomes homogeneous.

\end{enumerate}
\end{prop}

\begin{proof}
By direct calculations, we see that the conditions of the weights
of coordinates so that the equations (\ref{eq:ClusterEq}) are homogeneous
are reduced to the following if all the coordinates are nonzero:{\small
\begin{align*}
w(\theta_{1}) & =w(A_{23})+w(A_{3})+w(A_{31})-w(t_{2}),\,w(\theta_{3})=w(A_{1})+w(A_{12})+w(A_{31})-w(t_{2}),\\
w(A_{2}) & =w(A_{1})+w(A_{3})+2w(A_{31})-2w(t_{2}),\\
w(\theta_{12}) & =-w(A_{12})+w(A_{23})+w(A_{3})+w(A_{31}),\\
w(\theta_{23}) & =w(A_{1})+w(A_{12})-w(A_{23})+w(A_{31}),\\
w(\theta_{31}) & =w(A_{1})+w(A_{12})+w(A_{23})+w(A_{3})+w(A_{31})-2w(t_{2}),\\
w(\lambda) & =w(A_{1})+w(A_{3})+w(A_{31})-w(t_{2}).
\end{align*}
}{\small\par}

For the assertion (1), we have only to check that $\theta_{3}+\nicefrac{1}{2}\lambda A_{12}$
is weighted homogeneous, namely, it holds that $w(\theta_{3})=w(\lambda)+w(A_{12})$.
This is easily seen by the above conditions of the weights since $w(A_{3})=0$. 

For the assertion (2), we have only to check that $-\theta_{1}+\nicefrac{1}{2}\lambda A_{23}$,
$A_{31}-\nicefrac{1}{2}\lambda\theta_{2}$, $\theta_{3}+\nicefrac{1}{2}\lambda A_{12}$
and $-\nicefrac{1}{4}\lambda^{2}-A_{2}$ are weighted homogeneous,
namely, it holds that $w(\theta_{1})=w(\lambda)+w(A_{23})$, $w(\lambda)+w(\theta_{2})=w(A_{31}),$
$w(\theta_{3})=w(\lambda)+w(A_{12})$, and $2w(\lambda)=w(A_{2})$.
These are again easily seen by the above conditions of the weights
since $w(A_{1})=w(A_{3})=0$.
\end{proof}

\subsection{$\mathscr{H}_{\mathbb{A}}^{12}$,~$\mathscr{H}_{\mathbb{A}}^{11}$,~$\mathscr{M}_{\mathbb{A}}^{8}$
or $\mathscr{S}_{\mathbb{A}}^{6}$ as key varieties }

By the above considerations in Section \ref{sec:The--cluster-variety},
the result of \cite{CD1,CD2}, and Theorem \ref{thm:embthm}, we have
the following:
\begin{cor}
\label{cor:example}There exist quasi-smooth prime $\mathbb{Q}$-Fano
$3$-folds $X$ anticanonically embedded of codimension $4$ for the
$109$ numerical data of \cite{GRDB} such that they are obtained
as weighted complete intersections of weighted projectivizations $\mathscr{H}_{\mathbb{P}}^{12}$,
$\mathscr{H}_{\mathbb{P}}^{11}$,~$\mathscr{H}_{\mathbb{P}}^{10}$,~$\mathscr{M}_{\mathbb{P}}^{7}$
or $\mathscr{S}_{\mathbb{P}}^{5}$ of $\mathscr{H}_{\mathbb{A}}^{13}$,
$\mathscr{H}_{\mathbb{A}}^{12}$,~$\mathscr{H}_{\mathbb{A}}^{11}$,~$\mathscr{M}_{\mathbb{A}}^{8}$
or $\mathscr{S}_{\mathbb{A}}^{6}$ respectively or suitable weighted
cones of them. Precise information about which $X$ can be constructed
from which key varieties is as in Table \ref{tab:prime--Fano-3-folds}
below:

\begin{table}[H]
\begin{center}

\begin{tabular}{|c|c|}
\hline 
${\rm \text{key variety}}$ & $\#$ ${\rm \text{of classes}}$\tabularnewline
\hline 
\hline 
$\ensuremath{\mathscr{H}_{\mathbb{P}}^{12}}$ & $1$\tabularnewline
\hline 
$\ensuremath{\mathscr{H}_{\mathbb{P}}^{11}}$ & $10$\tabularnewline
\hline 
$\ensuremath{\mathscr{H}_{\mathbb{P}}^{10}}$ & $5$\tabularnewline
\hline 
$\ensuremath{\mathscr{M}_{\mathbb{P}}^{8}}$ & $32$\tabularnewline
\hline 
$\ensuremath{\mathscr{S}_{\mathbb{P}}^{5}}$ or $\ensuremath{\mathscr{S}_{\mathbb{P}}^{6}}$  & $61$\tabularnewline
\hline 
\end{tabular}

\end{center}

\caption{\label{tab:prime--Fano-3-folds}prime $\mQ$-Fano 3-folds obtained
from $\mathscr{H}_{\mA}^{13}$ }
\end{table}
\end{cor}

\end{document}